\documentclass[final]{siamltex}

\usepackage{color,amsmath,amssymb,amssymb,mathrsfs}
\usepackage{color}
\usepackage{graphicx}
\usepackage{showkeys,cite}
\usepackage[normalem]{ulem}
\usepackage[pdfborder={0 0 0}, colorlinks=true, linkcolor=blue]{hyperref}
\hypersetup{pdfborder=0 0 0}

\newcommand{\vv}{\ensuremath{\vec{v}}}
\newcommand{\ff}{\ensuremath{\vec{f}}}

\renewcommand{\matrix}[1] {\ensuremath{\boldsymbol{#1}}}
\renewcommand{\vec}[1] {\ensuremath{\boldsymbol{#1}}}

\newcommand{\M}{{\ensuremath{\matrix{M}}}}
\newcommand{\K}{{\ensuremath{\matrix{K}}}}
\newcommand{\A}{{\ensuremath{\matrix{A}}}}
\newcommand{\B}{{\ensuremath{\matrix{B}}}}
\newcommand{\V}{{\ensuremath{\matrix{V}}}}

\renewcommand{\L}{{\ensuremath{\matrix{L}}}}
\newcommand{\F}{{\ensuremath{\matrix{F}}}}
\newcommand{\iF}{{\ensuremath{\vec{F}}}}
\newcommand{\g}{\ensuremath{\boldsymbol{g}}}

\newcommand{\nn}{\ensuremath{\boldsymbol{n}}}

\renewcommand{\u}{\ensuremath{\boldsymbol{u}}}
\newcommand{\w}{\ensuremath{\boldsymbol{w}}}
\newcommand{\HT}{\matrix{\tilde{H}}_{\text{misfit}}}

\newcommand{\obs}{\vec{y}}  
\newcommand{\ipar}{m}
\newcommand{\dpar}{\vec{m}}
\newcommand{\iparmap}{\ipar_\text{MAP}}
\newcommand{\dparmap}{\dpar_\text{MAP}}

\newcommand{\xx}{\ensuremath{\boldsymbol x}}

\newcommand{\x}{\xx}
\newcommand{\y}{\ensuremath{\boldsymbol y}}

\newcommand{\seclab}[1]{\label{sec:#1}}
\newcommand{\secref}[1]{Section~\ref{sec:#1}}
\newcommand{\figlab}[1]{\label{fig:#1}}
\newcommand{\figref}[1]{Figure~\ref{fig:#1}}

\newcommand{\eqnlab}[1]{\label{eq:#1}}
\newcommand{\eqnref}[1]{\eqref{eq:#1}}

\definecolor{pacificorange}{cmyk}{.15,.45,1,0} 
\definecolor{pacificgray}{cmyk}{0,.15,.35,.60}
\definecolor{pacificlgray}{cmyk}{0,0,.2,.4}
\definecolor{pacificcream}{cmyk}{.05,.05,.15,0}
\definecolor{deepyellow}{cmyk}{0,.17,.80,0}
\definecolor{lightblue}{cmyk}{.49,.01,0,0}
\definecolor{lightbrown}{cmyk}{.09,.15,.34,0}
\definecolor{deepviolet}{cmyk}{.79,1,0,.15}
\definecolor{deeporange}{cmyk}{0,.59,1,.18}
\definecolor{dustyred}{cmyk}{0,.7,.45,.4}
\definecolor{grassgreen}{RGB}{92,135,39}
\definecolor{pacificblue}{RGB}{59,110,143}
\definecolor{pacificgreen}{cmyk}{.15,0,.45,.30}
\definecolor{deepblue}{cmyk}{1,.57,0,.2}
\definecolor{turquoise}{cmyk}{.43,0,.24,0}
\definecolor{green}{rgb}{0,0.65,0}

\newcommand{\bs}[1]{\boldsymbol{#1}}

\newcommand{\CC}{\ensuremath{c}}
\newcommand{\tCC}{\ensuremath{\tilde{c}}}

\newcommand{\equaldef}{{:=}}
\newcommand{\mb}[1]{\mathbf{#1}}
\newcommand{\mc}[1]{\mathcal{#1}}
\newcommand{\nor}[1]{\left\| #1 \right\|}

\renewcommand{\u}{p}
\newcommand{\yobs}{\obs^{\text{obs}}}
\newcommand{\R}{\mathbb{R}}
\newcommand{\half} {\ensuremath{\frac{1}{2}}}
\newcommand{\D}{\Omega}
\newcommand{\Acal}{\mc{A}}

\newcommand{\X}{X}

\newcommand{\C}{\mc{C}}
\newcommand{\Grad} {\ensuremath{\nabla}}
\newcommand{\DD}[2] {\ensuremath{\frac{d {#1}}{d {#2}}}}

\newcommand{\E}{E}
\newcommand{\LRp}[1]{\left( #1 \right)}
\newcommand{\LRs}[1]{\left[ #1 \right]}

\newcommand{\LRc}[1]{\left\{ #1 \right\}}
\newcommand{\GM}[2]{\mc{N}\left( #1, #2 \right)}

\newcommand{\adjMacroMM}[1]{{#1}^*}
\newcommand{\adjMacroME}[1]{{#1}^\natural}
\newcommand{\adjMacroEM}[1]{{#1}^\diamond}

\newcommand{\Badj}{\adjMacroMM{\B}}

\newcommand{\iFadj}{\adjMacroME{\iF}}
\newcommand{\Fadj}{\adjMacroME{\F}}

\newcommand{\Vadj}{\adjMacroEM{\V}}
\newcommand{\Ladj}{\adjMacroEM{\L}}

\def \addressices{Institute for Computational Engineering \& Sciences, The
  University of Texas at Austin, Austin, TX, USA}

\def \addressgeo{Department of Geological Sciences, The University of
  Texas at Austin, Austin, TX, USA}

\def \addressmech{Department of Mechanical Engineering, The
  University of Texas at Austin, Austin, TX, USA}

\begin{document}

\author{Tan Bui-Thanh\footnotemark[4] 
\and Omar Ghattas\footnotemark[4]
  \footnotemark[2] \footnotemark[3] \and James Martin\footnotemark[4]
  \and Georg Stadler\footnotemark[4] }

\renewcommand{\thefootnote}{\fnsymbol{footnote}}
\footnotetext[4]{\addressices}
\footnotetext[2]{\addressmech}
\footnotetext[3]{\addressgeo}
\renewcommand{\thefootnote}{\arabic{footnote}}

\title{A computational framework for 
  \\ infinite-dimensional Bayesian inverse problems.\\ Part I: The
  linearized case, with application to global seismic
  inversion\thanks{This research was partially supported by AFOSR
    grant FA9550-09-1-0608; DOE grants DE-FC02-11ER26052,
    DE-FG02-09ER25914, DE-FG02-08ER25860, DE-FC52-08NA28615, and DE-SC0009286; and
    NSF grants CMS-1028889, OPP-0941678, and DMS-0724746. We are
    grateful for this support.}}

\maketitle

\begin{abstract}
We present a computational framework for estimating the uncertainty in
the numerical solution of linearized infinite-dimensional statistical
inverse problems. We adopt the Bayesian inference formulation: given
observational data and their uncertainty, the governing forward
problem and its uncertainty, and a prior probability distribution
describing uncertainty in the parameter field, find the posterior
probability distribution over the parameter field. The prior must be
chosen appropriately in order to guarantee well-posedness of the
infinite-dimensional inverse problem and facilitate computation of the
posterior. Furthermore, straightforward discretizations may not lead
to convergent approximations of the infinite-dimensional problem. And
finally, solution of the discretized inverse problem via explicit
construction of the covariance matrix is prohibitive due to the need
to solve the forward problem as many times as there are parameters.

Our computational framework builds on the infinite-dimensional
formulation proposed by Stuart \cite{Stuart10}, and incorporates a
number of components aimed at ensuring a convergent discretization of
the underlying infinite-dimensional inverse problem. The framework
additionally incorporates algorithms for manipulating the prior,
constructing a low rank approximation of the data-informed component
of the posterior covariance operator, and exploring the posterior that
together ensure scalability of the entire framework to very high
parameter dimensions.  We demonstrate this computational framework on
the Bayesian solution of an inverse problem in 3D global seismic wave
propagation with hundreds of thousands of parameters.

\end{abstract}



\pagestyle{myheadings}
\thispagestyle{plain}
\markboth{T. Bui-Thanh, O. Ghattas, J. Martin, G. Stadler}
{COMPUT. FRAMEWORK FOR INF.-DIM.\ BAYESIAN INVERSION}

\section{Introduction}
\label{introduction}

We present a scalable computational framework for the quantification
of uncertainty in large scale {\em inverse problems}; that is, we seek
to estimate probability densities for uncertain
parameters,\footnote{We use the term {\em parameters} broadly to
  describe general model 
inputs that may be subject to uncertainty, which might include model
  parameters, boundary conditions, initial conditions, sources,
  geometry, and so on.} given noisy observations or measurements and a
model that maps parameters to output observables. The forward
problem---which, without loss of generality, we take to be governed by
PDEs---is usually well-posed (the solution exists, is unique, and is
stable to perturbations in inputs), causal (later-time solutions
depend only on earlier time solutions), and local (the forward
operator includes derivatives that couple nearby solutions in space
and time). The inverse problem, on the other hand, reverses this
relationship by seeking to estimate uncertain parameters from
measurements or observations. The great challenge of solving inverse
problems lies in the fact that they are usually ill-posed, non-causal,
and non-local: many different sets of parameter values may be
consistent with the data, and the inverse operator couples solution
values across space and time.

Non-uniqueness stems in part from the sparsity of data and the
uncertainty in both measurements and the PDE model itself, and in part
from non-convexity of the parameter-to-observable map. The popular
approach to obtaining a unique ``solution'' to the inverse problem is
to formulate it as an optimization problem: minimize the misfit
between observed and predicted outputs in an appropriate norm while
also minimizing a {\em regularization} term that penalizes unwanted
features of the parameters.
Estimation of parameters using the regularization approach to inverse
problems as described above will yield an estimate of the ``best''
parameter values that simultaneously fit the data and minimize the
regularization penalty term. However, we are interested not just in
point estimates of the best-fit parameters, but a {\em complete
  statistical description} of the parameters values that is consistent
with the data. The {\em Bayesian} approach
\cite{KaipioSomersalo05,Tarantola05} does this by reformulating the
inverse problem as a problem in {\em statistical inference},
incorporating uncertainties in the observations, the
parameter-to-observable map, and prior information on the
parameters. The solution of this inverse problem is the {\em
  posterior} probability distribution of the parameters, which
reflects the degree of confidence in their values. Thus we are able to
quantify the resulting uncertainty in the parameters, taking into
account uncertainties in the data, model, and prior information.

The inverse problems we target here are characterized by {\em infinite
  dimensional} parameter fields. This presents multiple difficulties,
including proper choice of prior to guarantee well-posedness of the
infinite-dimensional inverse problem, proper discretization to assure
convergence to solutions of the infinite-dimensional problem, and
algorithms for constructing and manipulating the posterior covariance
matrix that insure scalability to very large parameter dimensions.
The approach we adopt in this paper follows \cite{Stuart10}, which
seeks to first fully specify the statistical inverse problem on the
infinite-dimensional parameter space.  In order to accomplish this
goal, we postulate the prior distribution as a Gaussian random field
with covariance operator given by the square of the inverse of an
elliptic PDE.  This choice ensures that samples of the parameter field
are (almost surely) continuous as functions, and that the statistical
inverse problem is well-posed. To achieve a finite-dimensional
approximation to the infinite-dimensional solution, we carefully
construct a function-space-aware discretization of the parameter
space. 

The remaining challenge presented by infinite-dimensional statistical
inverse problems is in computing statistics of the (discretized)
posterior distribution. This is notoriously challenging for inverse
problems governed by expensive-to-solve forward problems and
high-dimensional parameter spaces (as in our application to global
seismic wave propagation in Section \ref{seismic}). The difficulty
stems from the fact that evaluation of the probability of each point
in parameter space requires solution of the forward PDE problem (which
can take many hours on a large supercomputer), and many such
evaluations (millions or more) are required to adequately sample the
(discretized) posterior density in high dimensions by conventional
Markov-chain Monte Carlo (MCMC) methods.
In complementary work \cite{MartinWilcoxBursteddeEtAl12}, we are
developing methods that accelerate MCMC sampling of the posterior by
employing a local Gaussian approximation of the posterior as a
proposal density, which is computed from the Hessian of the negative
log posterior. Here, as an alternative, we consider the case of the
linearized inverse problem; by linearization we mean that the
parameter-to-observable map is linearized about the point that
maximizes the posterior, which is known as the maximum {\em a
  posteriori} (MAP) point. With this linearization, the posterior
becomes Gaussian, and its mean is given by the MAP point; this can be
found by solving an appropriately weighted regularized nonlinear least
squares optimization problem. Furthermore, the posterior covariance
matrix is given by the Hessian of the negative log posterior evaluated
at the MAP point. 

Unfortunately, straightforward computation of the---nominally
dense---Hessian is prohibitive, requiring as many forward-like solves
as there are uncertain parameters (which in our example problem in
Section \ref{seismic}, is hundreds of thousands). However, the data
are typically informative about a low dimensional subspace of the
parameter field: that is, the Hessian of the data misfit term is
 a compact operator that is sparse with respect to some
basis. We exploit this fact to construct a low rank approximation of
the (prior preconditioned) data misfit Hessian using matrix-free
Lanczos iterations \cite{MartinWilcoxBursteddeEtAl12,
  FlathWilcoxAkcelikEtAl11}, which we observe to require a
dimension-independent number of iterations. Each iteration requires a
Hessian-vector product, which amounts to just a pair of
forward/adjoint PDE solves, as well as a prior covariance operator
application.  Since we take the prior covariance in the form of the
inverse of an elliptic differential operator, its application can be
computed scalably via multigrid. The Sherman-Morrison-Woodbury formula
is then invoked to express the covariance of the posterior. Finally,
we show that the resulting expressions necessary for visualization and
interrogation of the posterior distribution require just elliptic PDE
solves and vector sums and inner products. In particular, the
corresponding dense operators are never formed or stored.
Solving the statistical inverse problem thus reduces to solving a
fixed number of forward and adjoint PDE problems as well as an
elliptic PDE representing the action of the prior. Thus, when the
forward  PDE problem can be solved in a scalable manner (as it is for
our seismic wave propagation example in Section \ref{seismic}), the
entire computational framework is scalable with respect to forward
problem dimension, uncertain parameter field dimension, and data
dimension.

The computational framework presented here is applied to a sequence of
realistic large-scale 3D Bayesian inverse problems in global
seismology, in which the acoustic wavespeed of an unknown
heterogeneous medium is to be inferred from noisy waveforms recorded
at sparsely located receivers. Numerical results are presented for
several problems with the number of unknown parameters up to
431,000. We have employed a similar approach for problems with more than
one million parameters in related work 
\cite{Bui-ThanhBursteddeGhattasEtAl12}.

In the following sections, we provide an overview of the framework for
infinite-dimensional Bayesian inverse problems following
\cite{Stuart10} (Section \ref{sec:bayesian}), present a consistent
discretization scheme (Section \ref{sec:finite_bayesian}) for the
infinite-dimensional problem, summarize a method for computing the MAP
point (Section \ref{sec:findMAP}), describe our low rank-based
covariance approximation (Section \ref{lowrank}), and present results
of the application of our framework to the Bayesian solution of an
inverse problem in 3D global seismic wave propagation (Section
\ref{seismic}).

\section{Bayesian framework for infinite-dimensional inverse problems}
\label{sec:bayesian}

\subsection{Overview}
In the Bayesian formulation, we state the inverse problem as a problem
of \emph{statistical inference} over the space of parameters. The
solution of the resulting statistical inverse problem is a posterior
probability distribution that reflects our degree of confidence that
any set of candidate parameters might contain the actual values that
gave rise to the data via the model and were consistent with the prior
information. Bayes' formula, presented in its infinite dimensional
form in Section~\ref{sec:InfiniteBayes}, defines this posterior
probability distribution by combining a prior probability distribution
with a likelihood model.

The inversion parameter is a function assumed to be defined over an
open, bounded, and sufficiently regular set $\D \subset \R^3$.
The statistical inverse problem is therefore naturally posed in an
appropriate function space setting.
Here, we adopt the infinite-dimensional framework developed in
\cite{Stuart10}. In particular, we choose a prior that ensures
sufficient regularity of the parameter as required for the statistical
inverse problem to be well-posed.
We will represent the prior as a Gaussian random field whose
covariance operator is the inverse of an elliptic differential
operator. For certain problems, non-Gaussian priors can be important,
but the use of
non-Gaussian priors in statistical inverse problems is still subject
to active research, in particular for infinite-dimensional parameters.
Thus, here we restrict ourselves to priors given by Gaussian random
fields.
Let us motivate the choice of the covariance operator as inverse
of an elliptic differential operator by considering two
alternatives.  A common choice for covariance operators in statistical
inverse problems with a moderate number of parameters is to specify
the covariance function, which gives the covariance of the parameter
field between any two points.  This necessitates either construction
and ``inversion'' of a dense covariance matrix or expansion in a
truncated Karhunen-Lo\'eve (KL) basis. In the large-scale setting,
inversion of a dense covariance matrix is clearly intractable, and the
truncated KL approach can be impractical since it may require many
terms to prevent biasing of the solution toward the strong prior
modes.
On the contrary, specifying the covariance as the inverse of an
elliptic differential operator enables us to build on existing fast
solvers for elliptic operators
without constructing the dense operator.  Discretizations of elliptic
operators often satisfy a conditional independence property, which
relates them to Gaussian Markov random fields and allows for
statistical interpretation \cite{RueHeld05,Bardsley13}.  Even if a
Gaussian Markov random field is not based on an elliptic differential
operator, this Markov property permits the use of fast,
sparsity-exploiting algorithms for instance for taking samples from
the distribution, \cite{Rue01}.
Our implementation employs multigrid as solver for the discretized
elliptic systems.

A useful prior distribution must have bounded variance and have
meaningful realizations. In our infinite-dimensional setting, we
require samples to be pointwise well-defined, for instance,
continuous. Furthermore, it is convenient to have the ability to apply
the square root of the covariance operator, e.g., this is used to
compute samples from a Gaussian distribution.  We consider a Gaussian
random field $m$ on a domain $\Omega\subset\mathbb{R}^3$ with mean
$\ipar_0$ and covariance function $c(\x,\y)$ describing the covariance
between $\ipar(\x)$ and $\ipar(\y)$
\begin{equation}
\eqnlab{CovarianceFunctionPrior}
c(\x,\y) = \mathbb{E} \LRs{ (\ipar(\x) - \ipar_0(\x)) (\ipar(\y) -
  \ipar_0(\y))  } \:\text{ for } \x,\y \in \Omega.
\end{equation}
The corresponding covariance operator $\mc{C}_0$ is
\begin{equation}\eqnlab{covarianceC}
\LRp{\C_0\phi}\LRp{\x} = \int_\D c\LRp{\x,\y} \phi\LRp{\y}\, d \y
\end{equation}
for sufficiently regular functions $\phi$ defined over $\Omega$.
Thus, if the covariance operator is given by the solution operator of
an elliptic PDE, the covariance function is the corresponding Green's
function. Thus, Green's function properties have direct implications
for properties of the random field $m$. For instance, since Green's
functions of the Laplacian in one spatial dimension are bounded, the
random field with the Laplacian as covariance operator is of bounded
variance. However, in two and three space dimensions, Green's
functions $c(\x,\y)$ of the Laplacian are singular along the diagonal,
and thus the corresponding distribution has unbounded variance.  Thus,
intuitively the PDE solution operator used as covariance operator
$\C_0$ has to be sufficiently smoothing and have bounded Green's
functions. Indeed, this is necessary for the well-posedness of the
infinite-dimensional Bayesian formulation \cite{Stuart10}. The
biharmonic operator, for example, has bounded Green's functions in two
and three space dimensions.  We choose $\C_0=\mc{A}^{-2}$, where
$\mc{A}$ is a Laplacian-like operator specified in
Section~\ref{subsec:prior}.   This provides the desired
simple and fast-to-apply square root operator $\C_0^{1/2}=\mc{A}^{-1}$
and allows a straightforward discretization.

An approach to extract information from the posterior distribution is
to find the maximum a posterior (MAP) point, which amounts to the
solution of an optimization problem as summarized in
Section~\ref{sec:infdim_map_point}. Finally, in
Section~\ref{sec:linearized}, we introduce a linearization of the
parameter-to-observable map. This results in a Gaussian approximation
of the posterior, which is the main focus of this paper.

\subsection{Bayes' formula in infinite dimensions}
\label{sec:InfiniteBayes}
To define Bayes' formula, we require a likelihood function that
defines, for a
given parameter field $\ipar$, the distribution of observations
$\yobs$. Here, we assume a finite-dimensional vector $\yobs \in
\mathbb{R}^q$ of such observations. We introduce the
\emph{parameter-to-observable map} $\ff: X:=L^2(\Omega) \to \mathbb{R}^q$
as a deterministic function mapping a parameter field $\ipar$ to
so-called observables $\obs \in \mathbb{R}^q$, which are predictions
of the observations.  For the problems targeted here, an evaluation of
$\ff(\ipar)$ requires a PDE solve followed by the application of an
observation operator to extract $\obs$ from the PDE solution. Even
when the parameter $\ipar$ coincides with the ``true'' parameter, the
observables $\obs$ may still differ from the measurements $\yobs$ due
to measurement noise and inadequacy (i.e., the lack of fidelity of the
governing PDEs with respect to reality) of the parameter-to-observable
map $\ff$. As is common practice, we assume the discrepancy between
$\obs$ and $\yobs$ to be described by a Gaussian additive noise
$\vec\eta \sim \mu_{\text{noise}} = \GM{0}{\bs{\Gamma}_\text{noise}}$,
independent of $\ipar$.  In particular, we have
\begin{equation}\label{equ:Nobservation}
\yobs = \bs f\LRp{\ipar} + \vec\eta,
\end{equation}
which allows us to write the likelihood probability density function (pdf) as
\begin{equation}
\eqnlab{likelihoodInf}
\pi_{\text{like}}(\yobs | \ipar)
  \propto \exp \LRp{ -\frac12 (\ff\LRp{\ipar} -
  \yobs)^T \bs{\Gamma}_\text{noise}^{-1} (\ff\LRp{\ipar} - \yobs) }.
\end{equation}
The Bayesian solution to the infinite-dimensional inverse problem is
then defined as follows: given the likelihood $\pi_{\text{like}}$ and
  the prior measure $\mu_0$, find the conditional measure $\mu^y$ of
  $\ipar$ that satisfies the Bayes' formula
\begin{equation}\label{eq:BayesianSolution}
\DD{\mu^y}{\mu_0} = \frac{1}{Z} \pi_{\text{like}}(\yobs | \ipar),
\end{equation}
where $Z = \int_{\X}\pi_{\text{like}}(\yobs|\ipar) \, d\mu_0$ is a
normalization constant.  The formula \eqnref{BayesianSolution} is
understood as the Radon-Nikodym derivative of the posterior
probability measure $\mu^y$ with respect to the prior measure $\mu_0$.
In order for \eqref{eq:BayesianSolution} to be well defined,
$\ff:\X \to \mathbb{R}^q$ is assumed to be
locally Lipschitz and quadratically bounded in the sense of Assumption
2.7 in \cite{Stuart10}.
While the
Bayes' formula \eqref{eq:BayesianSolution} is
valid in finite and infinite dimensions, a more intuitive form of Bayes'
formula that uses Lebesgue measures and thus only holds in finite
dimensions is given in Section \ref{sec:finitePosterior}.

\subsection{Parameter space and the prior}
\label{subsec:prior}

As discussed in the introduction of Section~\ref{sec:bayesian}, we use
a squared inverse elliptic operator as covariance operator $\mc{C}_0$
in \eqref{eq:CovarianceFunctionPrior}, i.e., $\C_0=\mc{A}^{-2}$.
We first specify the
elliptic PDE corresponding to $\mc A$ in weak
form. For $s\in L^2(\Omega)$, the solution $\ipar=\mc A^{-1}s$ satisfies
\begin{equation}
\eqnlab{Wspace}
    \alpha\int_\Omega (\bs{\Theta} \nabla \ipar) \cdot \nabla p + \ipar p \,d\xx = \int_\Omega
    sp\,d\xx \text{ for all } p \in H^1(\Omega),
\end{equation}
with $\alpha > 0$, and $\bs{\Theta}(\xx) \in \mathbb{R}^{3 \times 3}$
is symmetric, uniformly bounded, and positive definite.  Note that for
$s\in L^2(\Omega)$, there exists a unique solution $\ipar \in
H^1(\Omega)$ by the Lax-Milgram theorem.
Since $s \in L^2(\D)$ in \eqref{eq:Wspace},
regularity results, e.g. \cite{ArbogastBona08,Evans98}, show that in fact $\ipar \in
H^2(\D)$ provided $\partial\D$ is sufficiently smooth, e.g., $\D$ is a
$C^{1,1}$ domain.  In this case, $(\ipar,s)$
satisfies the elliptic differential equation
\begin{subequations}\label{eq:Wstrong}
\begin{align}
-\alpha \nabla \cdot (\bs{\Theta} \nabla \ipar) + \alpha m &= s \quad \text{in } \D,\label{eq:Wstrong1}\\
\alpha (\bs{\Theta} \nabla \ipar) \cdot \nn &= 0 \quad \text{on } \partial\D,\label{eq:Wstrong2}
\end{align}
\end{subequations}
where $\nn$ denotes the outward unit normal on $\partial \Omega$.

Let us denote by $\Acal$ the differential operator together with its
domain of definition specified by \eqnref{Wstrong}; hence
$\Acal$ is a densely defined operator on $L^2(\D)$
with the following domain
\[
D\LRp{\Acal} \equaldef \LRc{
m \in H^2\LRp{\D}: \alpha \bs{\Theta} \Grad m \cdot \mb{n} = 0}.
\]
The operator $\Acal$ is assumed to be ``Laplacian-like'' in the sense
of Assumption 2.9 in \cite{Stuart10}.  In brief, this assumption
requires that $\Acal$ be positive definite, self-adjoint, invertible,
and have eigenfunctions that form an orthonormal basis of
$L^2(\Omega)$. Additionally, certain growth conditions on the
eigenvalues and $L^\infty(\Omega)$ norms of the
eigenfunctions are enforced\footnote{We note that this growth condition on the
  eigenfunctions may not be straightforward to demonstrate (or may not
  even hold) for a non-rectangular domain $\Omega$ and nonconstant
  coefficient $\bs{\Theta}$.  In these cases, we expect that
  alternative proofs of the results in \cite{Stuart10} can
  be accessed via regularity properties of the covariance function for
  the prior distribution.  See for example \cite{Abrahamsen97,
    LindgrenRueLindstroem11}. }.

To summarize, we consider $\ipar$ as a Gaussian random field whose
distribution law is a Gaussian measure $\mu_0 \equaldef
\GM{\ipar_0}{\mc{C}_0}$ on $\L^2(\D)$, with mean $\ipar_0 \in
D\LRp{\Acal}$ and covariance operator $\mc{C}_0 \equaldef \Acal^{-2}$.
The definition of the Gaussian prior measure is meaningful since
$\Acal^{-2}$ is a trace class operator on $L^2\LRp{\D}$
\cite{Stuart10}, which guarantees bounded variance and almost surely
pointwise well-defined samples since $\mu_0(X)=1$ holds, where
$X := C(\Omega)$ denotes the space of continuous functions defined on
$\Omega$ (see \cite[Lemma 6.25]{Stuart10}).

\subsection{The MAP point}
\label{sec:infdim_map_point}

As a first step in exploring the solution of the statistical inverse
problem, we determine the maximum a posteriori (MAP) estimate of the
posterior measure.  In a finite-dimensional setting, the MAP estimate
is the point
in parameter space that maximizes the posterior
probability density function.  This notion does not generalize
directly to the infinite-dimensional setting, but we can still define
the MAP estimate $\iparmap$ as the point $\ipar$ in parameter space
that asymptotically maximizes the measure of a ball with radius
$\varepsilon$ centered at $\ipar$, in the limit as $\varepsilon \to
0$.
We recall that the Cameron-Martin space $\E$ equipped with the inner
product $(\cdot\,,\cdot)_\E :=
(\mc{C}_0^{-1/2}\cdot,\mc{C}_0^{-1/2}\cdot)$ associated with
$\mc{C}_0$ is the range of $\mc{C}_0^{1/2}$ \cite{Hairer09}, and hence
coincides with $D\LRp{\Acal}$.
Using variational arguments, it can be shown (see \cite{Stuart10})
that  $\iparmap$ is given by solving the optimization problem
\begin{equation}
\label{eq:opt}
\min_{\ipar \in \E} 
\mc{J}\LRp{\ipar},
\end{equation}
where
\begin{equation}\label{eq:optfun}
\mc{J}\LRp{\ipar} \equaldef
\frac 12 \nor{\ff\LRp{\ipar} - \yobs}_{\bs{\Gamma}^{-1}_\text{noise}}^2 +
\half \nor{\mc{A}(\ipar-\ipar_0)}_{L^2\LRp{\D}}^2.
\end{equation}
The
well-posedness of the optimization problem \eqref{eq:opt} is guaranteed by
the  assumptions on $\bs f(\ipar)$ in Section~\ref{sec:InfiniteBayes}.

\subsection{A linearized Bayesian formulation}
\label{sec:linearized}

Once we have obtained the MAP estimate $\iparmap$, we approximate the
parameter-to-observable map $\bs f(\ipar)$ by its linearization about
$\iparmap$, which ultimately results in a Gaussian approximation to
the posterior distribution, as shown below. When the
parameter-to-observable map is nearly linear this is a reasonable
approximation; moreover, there are other scenarios in which the
linearization, and the resulting Gaussian approximation, may be
useful. Of particular interest here are the limits of small data noise
and many observations.  In the small noise case, the
parameter-to-obvervable map can be nearly linear as a mapping into the
subset of the observable space on which the likelihood distribution is
non-negligible---even when $\bs f(\ipar)$ is significantly nonlinear.
The asymptotic normality discussions in
\cite{GelmanCarlinSternEtAl03, LeCam86} suggest that under certain
conditions, the many observations case can lead to a Gaussian
posterior.
Finally, even if this approximation fails to describe the posterior
distribution adequately,
the linearization is still useful in building an initial step for the
rejection sampling approach or a Gaussian proposal distribution for
the Metropolis-Hastings algorithm \cite{RobertCasella04,
  MartinWilcoxBursteddeEtAl12}. These methods are
related to the sampling algorithm in
\cite{GirolamiCalderhead11}, which also employs derivative information to
respect the local structure of the parameter space.

Assuming that the parameter-to-observable map $\ff$ is Fr\'echet
differentiable, we linearize the right hand side of
\ref{equ:Nobservation} around $\iparmap$ to obtain
\[
\yobs \approx \ff\LRp{\iparmap} + \iF \LRp{\ipar - \iparmap} +
\vec\eta 
\]
where $\iF$ is the Fr\'echet derivative of $\ff\LRp{\ipar}$ evaluated at
$\iparmap$.  Consequently, the posterior distribution $\mu^y$ of
$\ipar$ conditional on $\yobs$ is a Gaussian measure
$\GM{\iparmap}{\C_\text{post}}$ with  mean $\iparmap$ and
covariance operator $\C_\text{post}$ defined by \cite{Stuart10}:
\begin{align}
\eqnlab{InfPostCov}
\C_\text{post} &= (\iFadj \bs{\Gamma}_\text{noise}^{-1} \iF + \C_0^{-1})^{-1},
\end{align}
with $\iFadj$ denoting the adjoint of $\iF$, an operator from the
space of observations $\mathbb{R}^q$ to $L^2\LRp{\Omega}$.  In
principle, a local Gaussian approximation of the posterior at the MAP
point can also be found for non-Gaussian priors and when the noise in
the observables is not additive and Gaussian as in
\eqref{equ:Nobservation}. In these cases,
however, even for a linear parameter-to-observable map
the local Gaussian approximation might
only be reasonable approximation to the true posterior distribution in
a small neighborhood around the MAP point.

\section{Discretization of the Bayesian inverse problem}
\label{sec:finite_bayesian}

\subsection{Overview}
Next, we present a numerical discretization of the
infinite-dimensional Bayesian statistical inverse problem described in
Section~\ref{sec:InfiniteBayes}.
The discretized parameter space is inherently high-dimensional (with
dimension dependent upon the mesh size).  If discretization is not
performed carefully at each step, it is unlikely that the discrete
solutions will converge to the desired infinite-dimensional solution
in a meaningful way \cite{Stuart10, LassasSaksmanSiltanen09}.

In the following, and particularly in Section
\ref{sec:finite_inner_product}, we choose a mass matrix-weighted
vector product instead of the standard Euclidean vector product. While
this is a natural choice in finite element discretizations \cite{BeskosRobertsStuart09,Voss12}, this does lead
to a few complications, for instance, the use of covariance operators
that are not symmetric in the conventional sense (they are
self-adjoint however). This choice is much better suited for proper
discretization of the infinite-dimensional expressions given in this
paper, and the resulting numerical expressions for computation will
more closely resemble their infinite-dimensional counterparts in
\secref{bayesian}.  By contrast, the correct corresponding expressions
in the Euclidean inner product are significantly less intuitive in our
opinion, and ultimately more cumbersome to manipulate
and interpret than the development we give here.

We provide finite-dimensional approximations of the prior and the
posterior distributions in Sections~\ref{sec:finite_prior} and
\ref{sec:finitePosterior}, respectively.
To study and visualize the uncertainty in Gaussian random fields, such
as the prior and posterior distributions, we generate realizations
(i.e., samples) and compute pointwise variance fields. This must be
done carefully in light of the mass-weighted inner products due to the
finite element discretization introduced in
Section~\ref{sec:finite_inner_product}. We present explicit
expressions for computing these quantities for the prior in the
Sections~\ref{sec:samples} and \ref{sec:variance}. The fast generation
of samples and the pointwise variance field from the Gaussian
approximation of the posterior exploits the low rank ideas presented
in Section~\ref{lowrank}. Thus, the presentation of the corresponding
expressions is postponed to Section~\ref{subsec:samplegeneration}.

\subsection{Finite-dimensional parameter space}
\label{sec:finite_function_space}

We consider a finite-dimensional subspace $V_h$ of $L^2(\Omega)$
originating from a finite element discretization with
continuous Lagrange basis functions
$\LRc{\phi_j}_{j=1}^n$, which correspond to the nodal points
$\LRc{\x_j}_{j=1}^n$, such that
\[
\phi_j(\x_i) = \delta_{ij}, \qquad \text{for } i,j \in \{1, \ldots, n\}.
\]
Instead of statistically inferring parameter functions $\ipar \in L^2(\Omega)$,
we perform this task on the approximation $\ipar_h = \sum_{j=1}^nm_j\phi_j
\in V_h$. Consequently, the
coefficients $\LRp{m_1,\hdots,m_n}^T \in \mathbb R^n$ are the actual
parameters to be inferred. For simplicity of notation, we
shall use the boldface symbol $\bs m = \LRp{m_1,\hdots,m_n}^T$ to denote the
nodal vector of a function $\ipar_h$ in $V_h$.

\subsection{Discrete inner product}
\label{sec:finite_inner_product}

Since we postulate the prior Gaussian measure on $L^2\LRp{\D}$, the
finite-dimensional space $V_h$ inherits the $L^2$-inner product. Thus,
inner products between nodal coefficient vectors must be weighted by a
mass matrix $\M$ to approximate the infinite-dimensional $L^2$-inner
product. We denote this weighted inner product by $( \cdot\,, \cdot
)_{\M}$ and assume that $\M \in \R^{n\times n}$ is symmetric and
positive definite.  To distinguish $\R^n$ with the $\M$-weighted inner
product from the usual Euclidean space $\R^n$, we denote it by
$\mathbb{R}^n_{\M}$.  For any $m_1,m_2\in L^2(\Omega)$, observe
that $(\ipar_1, \ipar_2)_{L^2(\Omega)} \approx (\ipar_{1h},
\ipar_{2h})_{L^2(\Omega)} = (\dpar_1, \dpar_2)_\M = \dpar_1^T \M
\dpar_2$, which motivates the choice of $\M$ as the finite element
mass matrix defined by
\begin{equation}
M_{ij} = \int_\Omega \phi_i(\x) \phi_j(\x) d\x ~, \quad i,j \in
\LRc{1,\hdots,n}.
\end{equation}

When using the $\M$-inner product, there is a critical distinction
that must be made between the matrix adjoint and the matrix
transpose. For an operator $\B : \R^n_{\M} \rightarrow \R^n_{\M}$,
we denote the matrix transpose by $\B^T$ with entries $[B^T]_{ij} =
B_{ji}$.  The adjoint $\Badj$ of $\B$,
however, must satisfy
\begin{equation}
\label{eq:implicit_innerprod_def}
\LRp{ \Badj \dpar_1, \dpar_2}_\M = \LRp{ \dpar_1, \B \dpar_2 }_\M
\: \text{ for all } \dpar_1, \dpar_2 \in \R^n_{\M}.
\end{equation}
This implies that
\begin{align}
\label{eq:MM_adj}
\Badj &= \M^{-1} \B^T \M.
\end{align}
In the following, we also need the adjoints $\Fadj$ of $\F : \R^n_{\M}
\rightarrow \R^q$ and $\Vadj$ of $\V : \R^r \rightarrow \R^n_{\M}$
(for some $r$), where $\R^q$ and $ \R^r$ are endowed with the
Euclidean inner product. The desired adjoints can be 
can be expressed as
\begin{align}
\label{eq:ME_adj}
\Fadj &= \M^{-1} \F^T,
\\
\label{eq:EM_adj}
\Vadj &= \V^T \M.
\end{align}
Next, let $P_h$ be the projection from $L^2\LRp{\D}$ to
$V_h$. Then, the matrix representation  $\B : \R^n_{\M} \rightarrow
\R^n_{\M}$ for the operator  $\mc{B}_h \equaldef P_h\mathcal B P_h'$, where 
$\mathcal{B}: L^2(\D) \rightarrow L^2(\D)$ and $P_h': V_h \rightarrow L^2(\D)$,
is implicitly given with respect to the Lagrange basis
$\LRc{\phi_i}_{i=1}^n$  in $V_h$ by
\begin{equation*}
\int_\Omega\phi_i\mathcal B\phi_j\,dx = (\bs e_i,\B\bs e_j)_\M,
\end{equation*}
where $\bs e_i$ is the coordinate vector corresponding to the basis
function $\phi_i$. As a result, one can write $\B$ explicitly as
\begin{equation}
\eqnlab{Arepresentation}
\B = \M^{-1}\K,
\end{equation}
where $\K$ is given by
\[
K_{ij} = \int_\Omega\phi_i\mathcal B\phi_j\,dx, \quad i,j \in \LRc{1,\hdots,n}.
\]

\subsection{Finite-dimensional approximation of the prior}
\label{sec:finite_prior}
Next, we derive the finite-dimensional representation of the prior.
The matrix representation of the operator $\Acal$ defined in
Section~\ref{subsec:prior} is given by the stiffness matrix
$\K$ with entries
\[
K_{ij}
= \alpha\int_\Omega
\left(\bs{\Theta}(\x)\nabla\phi_i(\x)\right)\cdot \nabla\phi_j(\x) +
\phi_i(\x)\phi_j(\x) \, d\x, \quad i,j \in \LRc{1,\hdots,n}.
\]
It follows that both $\A = \M^{-1}\K$ and
$\A^{-1}=\K^{-1}\M$ are self-adjoint operators in the sense of \eqref{eq:MM_adj}.

We are now in a position to define the finite-dimensional Gaussian prior
measure $\mu_0^h$ specified by the following density (with respect to
the Lebesgue measure):
\begin{equation}
\label{eq:prior_pdf}
\pi_{\text{prior}}(\dpar) \:\propto\:
 \exp\left[
- \frac 12 \left\| \A (\dpar - \dpar_{0}) \right\|^2_{\M}
\right].
\end{equation}
This definition implies that $\matrix{\Gamma}_{\text{prior}} = \A^{-2}$.

\subsection{Finite-dimensional approximation of the posterior}
\label{sec:finitePosterior}

In infinite dimensions, the Bayes' formula \eqref{eq:BayesianSolution}
has to be expressed in terms of the Radon--Nikodym derivative
since the prior and posterior distributions do not have density
functions with respect to the Lebesgue measure.
Since we approximate the prior
measure $\mu_0$ by $\mu_0^h$, it is natural to define a finite-dimensional approximation  $\mu^{y,h}$ of the posterior measure
$\mu^y$ such that
\[
\DD{\mu^{y,h}}{\mu_0^h} = \frac{1}{Z^h} \pi_{\text{like}}(\yobs | \ipar_h),
\]
where $Z^h = \int_{\X}\pi_{\text{like}}(\yobs|\dpar) \, d\mu_0^h$, and
$\pi_{\text{like}}$ is the likelihood \eqnref{likelihoodInf} evaluated
at $m_h$.  If we define $\pi_{\text{post}}(\dpar | \yobs)$ as the
density of $\mu^{y,h}$, again with respect to the Lebesgue measure, we
recover the familiar finite-dimensional Bayes' formula
\begin{align}
\label{eq:new_finite_bayes}
\pi_{\text{post}}(\dpar | \yobs) \:\propto\:
\pi_{\text{prior}}(\dpar) \pi_{\text{like}}(\yobs | \ipar_h),
\end{align}
where
the normalization constant
$1/Z^h$ , which does not depend on $\dpar$, is omitted.
Finally, we can express the posterior pdf explicitly as
\begin{equation}
\label{eq:explicit_finite_posterior}
\pi_{\text{post}}(\dpar)
\:\propto\:
\exp \left(
- \frac 12 \left\| \ff(\ipar_h) - \vec{y}^{\text{obs}} \right\|^2_{\matrix{\Gamma}^{-1}_{\text{noise}}}
- \frac 12 \left\| \A (\dpar - \dpar_{0}) \right\|^2_{\M}
\right),
\end{equation}
where, to recall our notation, $\ipar_h = \sum_{j=1}^n{\ipar_j\phi_j}
\in V_h$ and $\dpar = \LRp{\ipar_1,\hdots,\ipar_n}^T$.  We observe
that the negative log of the right side of
\eqnref{explicit_finite_posterior} is the finite-dimensional
approximation of the objective functional 
in \eqnref{opt}.

As a finite-dimensional counterpart of \secref{linearized}, we
linearize the parameter-to-observable map $\ff$ at the MAP point, but
now considering it as a function of the coefficient vector
$\dpar$. Let $\matrix{\Gamma}_{\text{post}}$ be the posterior
covariance matrix in the $\M$-inner product. Using
\eqnref{InfPostCov}, we obtain
\begin{equation}
\matrix{\Gamma}_{\text{post}} =
\left(
\Fadj\matrix{\Gamma}_{\text{noise}}^{-1}\F+\matrix{\Gamma}^{-1}_{\text{prior}}\right)^{-1},
\label{eqn:posterior}
\end{equation}
with $\Fadj = \M^{-1}\F^T$ as defined in \eqref{eq:ME_adj}. Note that
$\matrix{\Gamma}_{\text{post}}$ is self-adjoint, i.e.,
$\matrix{\Gamma}_{\text{post}} = \matrix{\Gamma}_{\text{post}}^*$ in
the sense of \eqref{eq:MM_adj}.

Since the posterior covariance matrix $\matrix{\Gamma}_{\text{post}}$ is
typically dense, we wish to avoid  explicitly storing it, especially when the
parameter dimension $n$ is large. Even if we are able to do so, it is
prohibitively expensive to construct. The reason is that the
Jacobian of the parameter-to-observable map, $\F$, is generally a
dense matrix, and its construction typically requires $n$ forward PDE
solves. This is clearly intractable when $n$ is large and solving the
PDEs is expensive. However, one can exploit the structure of the
inverse problem,
to approximate the posterior covariance matrix with desired
accuracy, as we shall show in Section~\ref{lowrank}.

\subsection{Sample generation in a finite element discretization}
\label{sec:samples}

We begin by developing expressions for a
general Gaussian distribution with mean $\bar{\dpar}$ and
covariance matrix $\matrix{\Gamma}$. Then, they are specified for the
Gaussian prior with $(\dpar_0, \matrix{\Gamma}_{\text{prior}})$.
Realizations of a finite-dimensional Gaussian random variable with
mean $\bar{\dpar}$ and covariance matrix $\matrix{\Gamma}$ can be
found by choosing a vector $\bs n$ containing independent and
identically distributed ({\em i.i.d.}) standard normal random values and computing
\begin{equation}
\dpar = \bar{\dpar} + \L\bs{n},
\end{equation}
where $\L$ is a linear map from $\mathbb{R}^n$ to $\R^n_{\M}$ such
that $\matrix{\Gamma} = \L \adjMacroEM{\L}$, in which the adjoint
$\adjMacroEM{\L} = \L^T\M$ (see also \eqnref{EM_adj}). Note that
$\M^{-1/2}\bs n$ is a sample from $\mathcal N(\bs 0,\bs I)$ in the
mass-weighted inner product.

In particular, for $\bar\dpar = \dpar_0$ and $\bs{\Gamma} =
\matrix{\Gamma}_{\text{prior}}$, we have $\L_{\text{prior}} = \K^{-1}
\M \M^{-1/2} = \K^{-1}\M^{1/2}$ (see the Appendix) and samples from the
prior are computed as $\dpar = \dpar_0 +
\K^{-1}\M^{1/2}\bs{n}$. Samples from the Gaussian approximation to the
posterior use the low-rank representation introduced in
Section~\ref{lowrank} and the corresponding expressions are given in
\eqref{eq:Lmatrix} and \eqref{eq:postsample}.


\subsection{The pointwise variance field in a finite element discretization}
\label{sec:variance}
Let us approximate the covariance function in $V_h$ for a generic
Gaussian measure with covariance operator $\mc{C}$.
Recall from Section \ref{subsec:prior} that
the covariance function $c\LRp{\x,\y}$ corresponding to the covariance
operator $\mc{C}$ is the Green's function of 
$\mc{C}^{-1}$, i.e.,
\[
\mc{C}^{-1}c\LRp{\x,\y} \equaldef \delta_{\y}(\x) \text{ for } \x\in \Omega,
\]
where $\delta_{\y}$ denotes the Dirac delta function concentrated at
$\y\in \Omega$.  We approximate $c\LRp{\x,\y}$ in the finite element
space $V_h$ by $c_h\LRp{\x,\y} = \sum_{i=1}^nc_i\LRp{\y} \phi_i\LRp{\x}$
with coefficient vector $\bs{c}\LRp{\y} = \LRs{c_1\LRp{\y}, \hdots,
  c_n\LRp{\y}}^T$. Using the Galerkin finite element method
to obtain a finite element approximation of the preceding equation
results in
\[
\mb{C}^{-1} \bs{c}\LRp{\y} = \bs{\Phi}\LRp{\y} \:\:\text{ with } \:\:
\bs{\Phi}\LRp{\y} = [ \phi_1(\y), \ldots, \phi_n(\y) ]^T
\]
and the entries of the matrix $\mb{C}^{-1}$ are given by
${C}_{ij}^{-1} = \LRp{\phi_i,\mc{C}^{-1}\phi_j}_{L^2(\Omega)}$. It follows
that $\bs{c}\LRp{\y} = \mb{C} \bs{\Phi}\LRp{\y}$ and
\[
c_h\LRp{\x,\y} = \sum_{i=1}^n c_i(\y)\Phi_i(\x) = \bs{\Phi}\LRp{\x}^T\mb{C} \bs{\Phi}\LRp{\y} \:\text{ for } \x,\y\in \Omega.
\]
Let us denote by $\bs{\Gamma}^{-1}$ the representation of
$P_h\mc{C}^{-1}P_h'$ in $V_h$; then, using \eqnref{Arepresentation} yields that
$
\mb{C} = \bs{\Gamma} \M^{-1}.
$
Consequently, the discretized covariance function for the covariance operator $\mc{C}$ now becomes
\begin{equation}
\eqnlab{genericCovariance}
c_h\LRp{\x,\y} = \bs{\Phi}\LRp{\x}^T\bs{\Gamma}\M^{-1} \bs{\Phi}\LRp{\y}.
\end{equation}
Let us now apply \eqref{eq:genericCovariance} to compute the prior
variance field. As discussed in Section~\ref{sec:finite_prior},
$\bs{\Gamma}_\text{prior} = \A^{-2} = \K^{-1}\M
\K^{-1}\M$. This results in the discretized prior covariance function
\[
c_h^{\text{prior}}\LRp{\x,\y} = \bs{\Phi}\LRp{\x}^T\K^{-1}\M \K^{-1}\bs{\Phi}\LRp{\y},
\]
By taking $\y = \x$, the prior variance field at an
arbitrary point $\x\in \Omega$ reads
\[
c_h^{\text{prior}}\LRp{\x,\x} =
\bs{\Phi}\LRp{\x}^T \K^{-1}\M \K^{-1} \bs{\Phi}\LRp{\x}.
\]
The pointwise variance field of the posterior distribution builds on
the low-rank representation introduced in Section~\ref{lowrank}.  The
resulting expression, which requires the prior variance field, is
given in \eqref{eq:postvarfield}.

\section{Finding the MAP point}
\label{sec:findMAP}

Section \ref{sec:infdim_map_point} introduced the idea of the MAP
point as a first step in exploring the solution of the statistical
inverse problem. To find the MAP point, one needs to solve a discrete
approximation (using the discretizations of Section
\ref{sec:finite_bayesian}) of the optimization problem \eqref{eq:opt},
which amounts to a large-scale nonlinear least squares numerical
optimization problem. In this section, we provide just a brief summary
of a scalable method we use for solving this problem, and refer the
reader to our earlier work for details, in particular the work on
inverse wave propagation \cite{AkcelikBielakBirosEtAl03,
  AkcelikBirosGhattas02, EpanomeritakisAkccelikGhattasEtAl08}.
We use an inexact matrix-free Newton-conjugate gradient (CG) method in
which only Hessian-vector products are required. These Hessian-vector
products are computed by solving a linearized forward-like and an
adjoint-like PDE problems, and thus the Hessian matrix is never
constructed explicitly. Inner CG iterations are terminated prematurely
when sufficient reduction is made in the norm of the gradient, or when
a direction of negative curvature is encountered. The prior operator
is used to precondition the CG iterations. Globalization is through an
Armijo backtracking line search.

Because the major components of the method can be expressed as
solving PDE-like systems,
the method inherits the scalability (with respect to problem
dimension) of the forward PDE solve.  The remaining ingredient for
overall scalability is that the optimization algorithm itself be
scalable with increasing problem size. This is indeed the case: for a
wide class of nonlinear inverse problems, the outer Newton iterations
and the inner CG iterations are independent of the mesh size, as is
found to be the case for instance for inverse wave propagation
\cite{AkcelikBirosGhattas02,
  EpanomeritakisAkccelikGhattasEtAl08}. This is a consequence of the
use of a Newton solver, of the compactness of the Hessian of the data
misfit term (i.e., the first term) in \eqref{eq:optfun}, and of the use
of preconditioning by $\matrix{\Gamma}_{\text{prior}}$, so that
the resulting preconditioned Hessian is a compact perturbation of the
identity, for which CG exhibits mesh-independent iterations.

\section{Low rank approximation of the Hessian matrix}
\label{lowrank}

\subsection{Overview}
As discussed in Section \ref{sec:linearized}, linearizing the
parameter-to-observable map results in the posterior covariance matrix
being given by the inverse of the Hessian of the negative log
posterior.  Explicitly computing this Hessian matrix requires a
(linearized) forward PDE problem for each of its columns, and thus as
many (linearized) forward PDE solves are required as there are
parameters.  For inverse problems in which one seeks to infer an
unknown parameter field, discretization results in a very large number
of parameters; explicitly computing the Hessian---and hence the
covariance matrix---is thus out of the question.  As a remedy, we
exploit the structure of the problem to find an approximation of the
Hessian that can be constructed and dealt with efficiently.

When the linearized parameter-to-observable map is used in
$\mathcal{J}(m)$ (as defined in~\eqref{eq:optfun}) and second
derivatives of the resulting functional are computed, one obtains the
Gauss-Newton portion of the Hessian of $\mathcal{J}(m)$. Both, the
full Hessian matrix as well as its Gauss Newton portion are positive
definite at the MAP point and they only differ in terms that involve the
adjoint variable. Since the adjoint system is driven only by the data
misfit (see, for instance, the adjoint wave
equation~\ref{eq:adjoint}), the adjoint variable is expected to be
small when the data misfit is small, which occurs provided the model
and observational errors are not too large. The Gauss-Newton portion
of the Hessian is thus often a good approximation of the full
Hessian of $\mathcal{J}(m)$.

For conciseness and convenience of the notation, we focus on computing
a low rank approximation of the Gauss-Newton portion of the (misfit)
Hessian in Section~\ref{subsec:lowrank}. The same approach also
applies to the computation of a low rank approximation of the full
Hessian, whose inverse might be a better approximation for the
covariance matrix if the data is very noisy and the data misfit at the
MAP point cannot be neglected.  The low rank construction of the
misfit Hessian is based on the Lanczos method and thus only requires
Hessian-vector products.  Using the Sherman-Morrison-Woodbury formula,
this approximation translates into an approximation of the posterior
covariance matrix.

In Section~\ref{subsec:samplegeneration}, we present low rank-exploiting
methods for sample generation from the Gaussian approximation of the
posterior, as well as methods for the efficient computation of the
pointwise variance field. Finally, in
Section~\ref{subsec:scalability}, we discuss the overall scalability
of our approach.

\subsection{Low rank covariance approximation}\label{subsec:lowrank}
For many ill-posed inverse problems, the Gauss-Newton
portion of the Hessian matrix (called the Gauss-Newton Hessian for
short) of the data misfit term in \eqref{eq:optfun} evaluated at
any $\mb{m}$,
\begin{equation}
\matrix{H}_{\text{misfit}} :=
   \Fadj\matrix{\Gamma}_{\text{noise}}^{-1}\F,
\label{Hmisfit}
\end{equation}
behaves like (the discretization of) a compact operator (see, e.g.,
\cite[p.17]{Vogel02}). The intuitive reason for this is that only
parameter modes that strongly influence the observations through the
linearized parameter-to-observable map $\F$ will be present in the
dominant spectrum of the Hessian \eqref{Hmisfit}.  For typical inverse
problems, observations are sparse, and hence the dimension of the
observable space is much smaller than that of the parameter
space. Furthermore, highly oscillatory perturbations in the parameter
field often have negligible effect on the output of the
parameter-to-observable map. In~\cite{Bui-ThanhGhattas12a,
  Bui-ThanhGhattas12}, we have shown that the Gauss-Newton
Hessian of the data misfit is a compact operator, and that for smooth media
its eigenvalues decay exponentially to zero. Thus, the range space of
the Gauss-Newton Hessian is effectively finite-dimensional even before
discretization, i.e., it is independent of the mesh.  We can exploit the
compact nature of the data misfit Hessian to
construct scalable algorithms
for approximating the inverse of the Hessian 
\cite{FlathWilcoxAkcelikEtAl11, MartinWilcoxBursteddeEtAl12}.

A simple manipulation of \eqref{eqn:posterior} yields the following
expression for the posterior covariance matrix, which will prove
convenient for our purposes:
\begin{align}
\label{eqn:gamma_post}
\matrix{\Gamma}_{\text{post}} =
\matrix{\Gamma}_{\text{prior}}^{1/2}\left(\matrix{\Gamma}_{\text{prior}}^{1/2}
\matrix{H}_{\text{misfit}}
\matrix{\Gamma}_{\text{prior}}^{1/2}+
\matrix{I}\right)^{-1}\matrix{\Gamma}_{\text{prior}}^{1/2}.
\end{align}
We now present a fast method for approximating
$\matrix{\Gamma}_{\text{post}}$ with controllable accuracy by making a
low rank approximation of the so-called {\em prior-preconditioned Hessian
of the data misfit}, namely,
\begin{equation}\label{eq:HT}
\HT{:=}
\matrix{\Gamma}_{\text{prior}}^{1/2}
\matrix{H}_{\text{misfit}}
\matrix{\Gamma}_{\text{prior}}^{1/2}.
\end{equation}
Let $\LRp{\lambda_i, \vv_i}, i = 1,\hdots,n$ be the eigenpairs of
$\HT$, and
$\matrix{\Lambda} = \diag(\lambda_1,\hdots,\lambda_n) \in \mathbb{R}^{n\times n}$.
Define $\V \in \mathbb{R}^{n\times
   n}$ such that its columns are the eigenvectors
 $\vv_i$ of $\HT$. Replacing
 $\HT$  by its spectral decomposition (recall that $\Vadj$ is the
adjoint of $\V$ as defined in \eqref{eq:EM_adj}),
\begin{equation*}
  \left(\matrix{\Gamma}_{\text{prior}}^{1/2}
\matrix{H}_{\text{misfit}}
  \matrix{\Gamma}_{\text{prior}}^{1/2}+ \matrix{I}\right)^{-1} =
  (\V\matrix{\Lambda} \Vadj + \matrix{I})^{-1} .
\end{equation*}
When the eigenvalues of $\HT$ decay rapidly we can construct a low-rank
approximation of $\HT$ by computing only the $r$ largest eigenvalues, i.e.,
\begin{equation*}
  \matrix{\Gamma}_{\text{prior}}^{1/2}
\matrix{H}_{\text{misfit}}
  \matrix{\Gamma}_{\text{prior}}^{1/2} =
  \matrix{V}_r\matrix{\Lambda}_r \Vadj_r +
  \mathcal{O}\left(\sum_{i=r+1}^{n} \lambda_i\right),
\end{equation*}
where $\matrix{V}_r \in \mathbb{R}^{n\times r}$ contains $r$
eigenvectors of $\HT$ corresponding to the $r$ largest eigenvalues
$\lambda_i, i=1,\hdots, r$, and $\matrix{\Lambda}_r = \diag
(\lambda_1, \hdots, \lambda_r) \in \mathbb{R}^{r \times r}$.  We can
then invert the approximate Hessian using the
Sherman-Morrison-Woodbury formula to obtain
\begin{equation}
\label{eqn:eigen_error}
\left(\matrix{\Gamma}_{\text{prior}}^{1/2}
\matrix{H}_{\text{misfit}}
\matrix{\Gamma}_{\text{prior}}^{1/2}+ \matrix{I}\right)^{-1}
=
\matrix{I}-\V_r \matrix{D}_r \Vadj_r +
\mathcal{O}\left(\sum_{i=r+1}^{n} \frac{\lambda_i}{\lambda_i + 1}\right),
 \end{equation}
where
$\matrix{D}_r
:=\diag(\lambda_1/(\lambda_1+1), \hdots, \lambda_r/(\lambda_r+1))
\in \mathbb{R}^{r\times r}$. Equation \eqref{eqn:eigen_error}
shows the truncation error due to the low-rank approximation based on
the first $r$ eigenvalues. To obtain an accurate approximation of
$\matrix{\Gamma}_{\text{post}}$, only eigenvectors corresponding to
eigenvalues that are small compared to $1$ can be neglected.
With such a low-rank
approximation, the final expression for the approximate posterior
covariance
is  given by
\begin{equation}
\label{eqn:reduced_post}
\matrix{\Gamma}_{\text{post}} \approx
\matrix{\Gamma}_{\text{prior}}
-\matrix{\Gamma}_{\text{prior}}^{1/2} \V_r \matrix{D}_r \Vadj_r
\matrix{\Gamma}_{\text{prior}}^{1/2}.
\end{equation}
Note that \eqref{eqn:reduced_post} expresses the posterior uncertainty
(in terms of the covariance matrix) as the prior uncertainty less any
information gained from the data. Due to the square root of the prior
in the rightmost term in \eqref{eqn:reduced_post}, the information
gained from the data is filtered through the prior, i.e., only
information consistent with the prior can reduce the posterior
uncertainty.

\subsection{Fast generation of samples and the pointwise variance
  field}\label{subsec:samplegeneration}
Properties of the last term in \eqref{eqn:reduced_post}, such as its
diagonal (which provides the reduction in variance due to the
knowledge acquired from the data) can be obtained numerically through
just $r$ applications of the square root of the prior covariance
matrix to $r$ columns of $\V_r$. Let us define
\[
  \tilde{\V_r} = \matrix{\Gamma}_{\text{prior}}^{1/2}
  \V_r ,
\]
then \eqref{eqn:reduced_post} becomes
\begin{equation}
\label{eqn:reduced_postfinal}
\matrix{\Gamma}_{\text{post}} \approx
\matrix{\Gamma}_{\text{prior}}
-\tilde{\V_r} \matrix{D}_r \tilde{\Vadj_r},
\end{equation}
with $\tilde{\Vadj_r} = \Vadj_r
\matrix{\Gamma}_{\text{prior}}^{1/2}$.

The linearized posterior is a Gaussian distribution with known
mean, namely the MAP point, and low rank-based covariance
\eqref{eqn:reduced_postfinal}.  Thus, the pointwise variance field and samples
can be generated as in Section \ref{sec:variance} and \ref{sec:samples}, respectively.
The variance field can be computed as
\begin{equation}\label{eq:postvarfield}
c^{\text{post}}_h(\x,\x) = c^{\text{prior}}_h(\x,\x)
- \sum_{k=1}^r d_k \LRp{ \tilde{\vv}_{kh}(\x) }^2,
\end{equation}
where $\tilde{\vv}_{kh}(\x) = \bs{\Phi}\LRp{\x}^T \tilde{\vv}_k$, with
$\tilde{\vv}_k$ denoting the $k$th column of $\tilde{\V_r}$, is the
function in $V_h$ corresponding to the nodal vector $\tilde \vv_k$.

Now, we can compute samples from the posterior provided that we have a
factorization $\matrix{\Gamma}_{\text{post}} = \L \Ladj$.  One
possibility for $\L $ (see the Appendix for the detailed derivation) reads
\begin{equation}
\label{eq:Lmatrix}
\L :=
\matrix{\Gamma}^{1/2}_{\text{prior}}
\LRp{
  \V_r \matrix{P}_r
  \Vadj_r + \matrix{I}
}\matrix{M}^{-1/2}
\end{equation}
with $\matrix{P}_r = \diag\LRp{1/\sqrt{\lambda_1+1}-1, \hdots,
  1/\sqrt{\lambda_r+1}-1} \in \R^{r\times r}$, $\L$ as a linear map
from $\R^n$ to $\R^n_M$, and $\matrix{I}$ as the identity map in both
$\R^n$ and $\R^n_\M$.  As discussed in Section
\ref{sec:samples}, samples can be then computed as
\begin{equation}\label{eq:postsample}
\bs{\nu}^{\text{post}} = \dparmap + \matrix{L} \bs{n},
\end{equation}
where $\bs{n}$ is an {\em i.i.d.} standard normal random vector.

\subsection{Scalability}\label{subsec:scalability}
We now discuss the scalability of the above low rank construction of
the posterior covariance matrix in \eqref{eqn:reduced_post}. The
dominant task is the computation of the dominant spectrum of the prior
preconditioned Hessian of the data misfit, $\HT$, given by
\eqref{eq:HT}. Computing the spectrum by a matrix-free eigensolver
such as Lanczos means that we need only form actions of $\HT$ with a
vector. As argued at the end of \secref{finitePosterior}, the
linearized parameter-to-observable map $\F$ is too costly to be
constructed explicitly since it requires $n$ linearized forward PDE
solves. However, its action on a vector can be computed by solving a
single linearized forward PDE (which we term the {\em incremental
  forward problem}), regardless of the number of parameters $n$ and
observations $q$. Similarly, the action of $\Fadj$ on a vector can be
found by solving a single linearized adjoint PDE (which we term the
{\em incremental adjoint problem}). Explicit expressions for the
incremental forward and incremental adjoint PDEs in the context of
inverse acoustic wave propagation will be given in Section
\ref{seismic}. Solvers for the incremental forward and adjoint
problems of course inherit the scalability of the forward PDE solver.
The other major cost in computing the action of $\HT$ on a vector is
the application of the square root of the prior,
$\matrix{\Gamma}_{\text{prior}}^{1/2}$, to a vector. As discussed in
Section \ref{subsec:prior}, this amounts to solving a Laplacian-like
problem.  Using a scalable elliptic solver such as multigrid renders
this component scalable as well. Therefore, the scalability of
the application of $\HT$ to a vector---which is the basic operation of
a matrix-free eigenvalue solver such as Lanczos---is assured, and the
cost is independent of the parameter dimension.

The remaining requirement for independence of parameter dimension in
the construction of the low rank-based representation of the posterior
covariance in \eqref{eqn:reduced_post} is that the number of dominant
eigenvalues of $\matrix{H}_{\text{misfit}}$ be independent of the
dimension of the discretized parameter. This is the case when
$\matrix{H}_{\text{misfit}}$ and $\matrix{\Gamma}_{\text{prior}}$ in
\eqref{Hmisfit} are discretizations of a compact and a continuous
operator, respectively. The continuity of $\mc{C}_0$ is a direct
consequence of the prior Gaussian measure $\mu_0$; in fact,
$\mc{C}_0$, the infinite-dimensional counterpart of
$\matrix{\Gamma}_{\text{prior}}$, is also a compact operator.
Compactness of the data misfit Hessian $\matrix{H}_{\text{misfit}}$
for inverse wave propagation problems has long been observed (e.g.,
\cite{ColtonKress98}) and, as mentioned above, has been proved for
frequency-domain acoustic inverse scattering for both continuous and
pointwise observation operators \cite{Bui-ThanhGhattas12,
  Bui-ThanhGhattas12a}. Specifically, we have shown that the data
misfit Hessian is a compact operator at any point in the parameter
domain.  We also quantify the decay of the data misfit Hessian eigenvalues
in terms of the smoothness of the medium, i.e., the smoother it is the
faster the decay rate. For an analytic target medium, the rate can be shown
to be exponential. That is, the data misfit Hessian can be
approximated well with a small number of its dominant eigenvectors
and eigenvalues. 

As a result, the number of Lanczos iterations required to obtain a low
rank approximation of $\HT$ is independent of the dimension of the
discretized parameter field.  Once the low-rank approximation of $\HT$
is constructed, no additional forward or adjoint PDE solves are
required.  Any action  of $\matrix{\Gamma}_{\text{post}}$
in \eqref{eqn:reduced_post} on a vector (which is required to generate
samples from the posterior distribution and to compute the
diagonal of the covariance) 
is now dominated by the 
action of $\matrix{\Gamma}_{\text{prior}}$ on a vector. But as
discussed above, this amounts to an elliptic solve and can be readily carried
out in a scalable manner. 
Since $r$ is independent of the dimension of the discretized parameter
field, estimating the posterior covariance matrix requires a constant
number of forward/adjoint PDE solves, independent of the number of
parameters, observations, and state variables. Moreover, since the
dominant cost is that of solving forward and adjoint PDEs as well as
elliptic problems representing the prior, scalability of the overall
uncertainty quantification method follows when the forward and adjoint
PDE solvers are scalable.

\section{Application to global seismic statistical inversion}
\label{seismic}

In this section, we apply the computational framework developed in the
previous sections to the statistical inverse problem of global seismic
inversion, in which we seek to reconstruct the heterogeneous
compressional (acoustic) wave speed from observed seismograms, i.e.,
seismic waveforms recorded at points on earth's surface. With the
rapid advances in observational capabilities, exponential growth in
supercomputing, and maturation of forward seismic wave propagation
solvers, there is great interest in solving the global seismic inverse
problem governed by the full acoustic or elastic wave equations
\cite{Fichtner11,PeterKomatitschLuoEtAl11}. Already, successful
deterministic inversions have been carried out at regional scales; for
example, see 
\cite{TapeLiuMaggiEtAl09, FichtnerIgelBungeEtAl09,
  FichtnerKennettIgelEtAl09, ZhuBozdaugPeterEtAl12, LekicRomanowicz11}.

We consider global seismic model
problems in which the seismic source is taken as a simple point source.
Sections~\ref{sec:wave_parameter} and \ref{sec:wave_prior} define the
prior mean and covariance operator for the wave speed and its
discretization. Section \ref{sec:wave_map} presents the
parameter-to-observable map $\bs f\LRp{\ipar}$ (which involves
solution of the acoustic wave equation) and the likelihood model. We
next provide the expressions for the gradient and application of the
Hessian  of the negative log-likelihood in
Section~\ref{sec:wave_gradhessian}. Then, we discuss the
discretization of the forward and adjoint wave equations and
implementation details in Section \ref{sec:wave_implementation}.
Section \ref{sec:wave_problemsetup} provides the inverse problem
setup, while numerical results and discussion are provided in Sections
\ref{sec:wave_lowrank} and \ref{sec:wave_resultsUQ}.

\subsection{Parameter space for seismic inversion}
\label{sec:wave_parameter}

\begin{figure}\centering
\hfill\hfill
\includegraphics[height=.37\columnwidth]{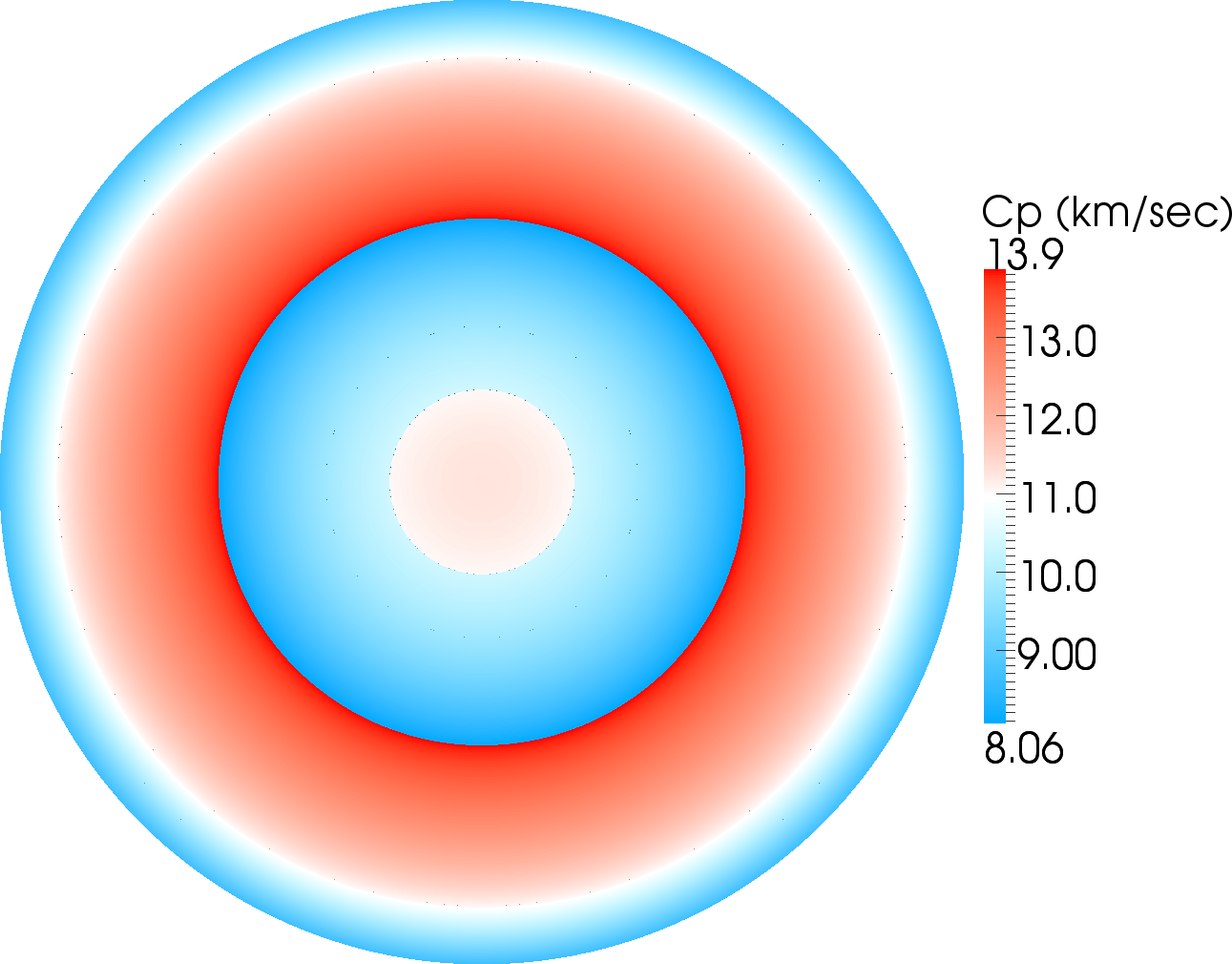}\hfill
\includegraphics[height=.37\columnwidth]{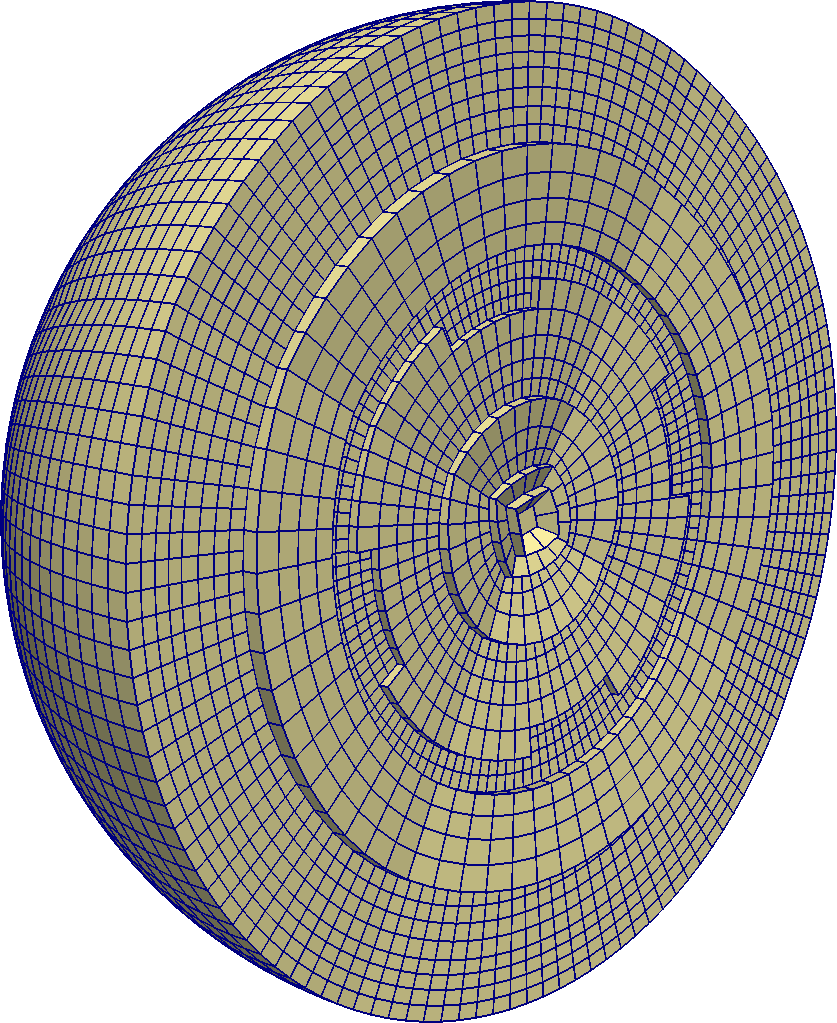}\hfill\hfill
\caption{Left image: Cross section through the spherically symmetric PREM
  earth model, which is the prior mean in the inversion.
  Right image: Mesh used for both wave speed parameters (discretized with $N=1$)
  and wave propagation unknowns ($N=3$). The mesh is tailored to the local wave
  lengths.}
\label{fig:s20rts}
\end{figure}

We are interested in inferring the heterogenous compressional acoustic
wavespeed in the earth.  In order to do this, we represent the earth as a
sphere of radius 6,371km.
We employ two earth models, i.e.\., two 
representations of the compressional wave speed and density in the
earth. We suppose that our current knowledge of the earth is given by
the spherically symmetric Preliminary Reference Earth Model (PREM)
\cite{DziewonskiAnderson81}, which is depicted in Figure
\ref{fig:s20rts}. 

The PREM is used as the mean of the prior distribution, and as the
starting point for the
determination of the MAP 
point by the optimization solver. 
Then, we presume that the real earth behaves according to the S20RTS
velocity model \cite{RitsemaVanHeijst02}, which superposes lateral
wave speed variations on the laterally-homogeneous PREM. S20RTS is
used to generate waveforms used as synthetic observational data for
inversion; we refer to it as the ``ground truth'' earth model.  The inverse
problem then aims to reconstruct the S20RTS ground truth model from
the (noisy) synthetic data and from prior knowledge of the PREM, and
quantify the uncertainty in doing so.

The parameter field $\ipar$ of interest for the 
inverse problem is the deviation or anomaly from the
PREM, and hence it is sensible to choose the zero funtion as the prior mean. 
Owing to the prior covariance operator specified in Section~\ref{subsec:prior}, the deviation is smooth; in fact it is continuous almost surely.
The wave speed parameter space is discretized using continuous
isoparametric trilinear finite elements on a hexahedral octree-based
mesh. To generate the mesh, we partition the earth into 3 layers
described by 13 mapped cubes. The first layer consists of a single
cube surrounded by two layers of six mapped cubes. 
The resulting mesh is aligned with the interface
between the outer core and the mantle, where the wave speed has a
significant discontinuity (see Figure \ref{fig:s20rts}). It is also aligned with
several known weaker discontinuities between layers.

The parameter mesh coincides with the mesh used to solve the wave
equation described in Section \ref{sec:wave_map}. 
The mesh is 
locally refined to resolve the local seismic wavelength resulting from
a given frequency of interest for the PREM. We choose a conservative
number of grid points per wavelength to permit the same mesh to be
used for anticipated variations in the earth model across the
iterations needed to determine the MAP point. 
For the parallel mesh generation and its distributed storage, we use
fast forest-of-octree algorithms for scalable adaptive mesh refinement
from the \texttt{p4est} library \cite{BursteddeWilcoxGhattas11,
  BursteddeGhattasGurnisEtAl10}. 

\subsection{The choice of prior}
\label{sec:wave_prior}

Since the prior is a Gaussian measure, it is completely specified by a
mean function and a covariance operator.  As discussed in Section
\ref{sec:wave_parameter}, the prior distribution for the anomaly (the
deviation of the acoustic wavespeed from that described by the PREM
model) is naturally chosen to have zero mean.
The choice of covariance operator for the prior distribution has to
encode several important features. Recall that we specify the
covariance operator via the precision operator $\mc{A}$ in Section
\ref{subsec:prior}. Therefore, the size of the variance about the zero
mean is set by $\alpha$, while the product $\alpha \Theta$ determines
the correlation length of the prior Gaussian random field. We next
specify the scalar $\alpha$ and the tensor $\Theta$ based on the
following observations of models for the local wave speeds in the
earth.
\begin{itemize}
\item Smoothness.
  The parameter field describes the
  \emph{effective} local wave-speed, which, for a finite
  source frequency, depends on the local average of material parameters
  within a neighborhood of each point in space. This makes the
  effective wave speed mostly a smooth field.
  Note that the S20RTS-based target wave speed model (see
  \cite{RitsemaVanHeijst02}) is smooth.
\item Prior variance. The deviation in this effective wave
  speed from the PREM model is believed to be within a few
  percent. Thus, we select $\alpha$ such that the prior standard
  deviation is about 3.5\%. The S20RTS target model has a maximal deviation
  from PREM of 7\%.
\item Anisotropy in the mantle. We further incorporate the prior
  belief that the compressional wave speed has a stronger variation
  in depth than in the lateral directions.  We encode this anisotropic
  variation through $\Theta$. In particular, we select $\Theta$ such
  that the anisotropy is strongest near the surface, and gradually
  becomes weaker with higher overall correlation length at larger
  depths. We observe that the S20RTS target model also obeys a similar
  anisotropy.
\end{itemize}

\begin{figure}
  \begin{center}
    \includegraphics[width=.22\columnwidth]{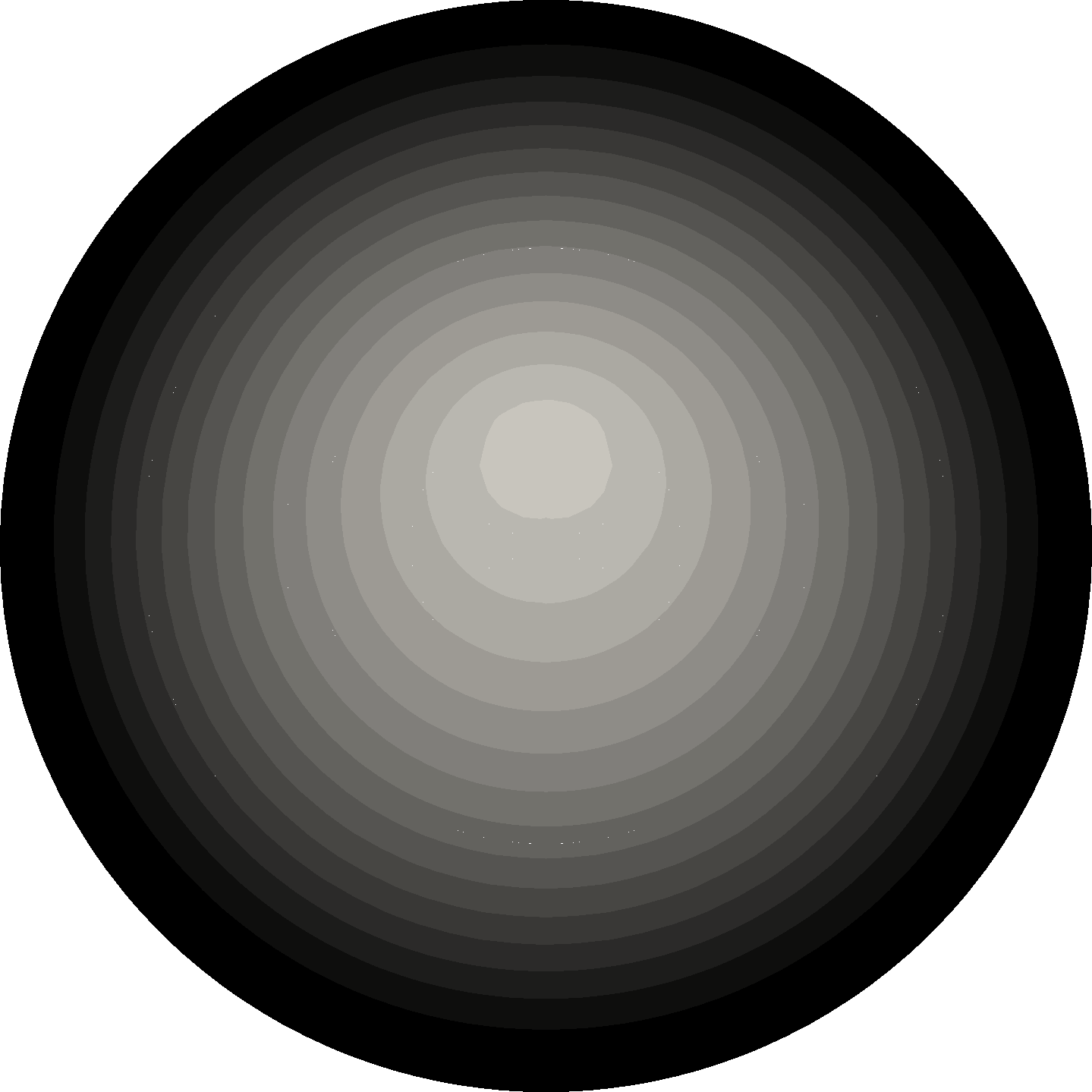} 
    \includegraphics[width=.22\columnwidth]{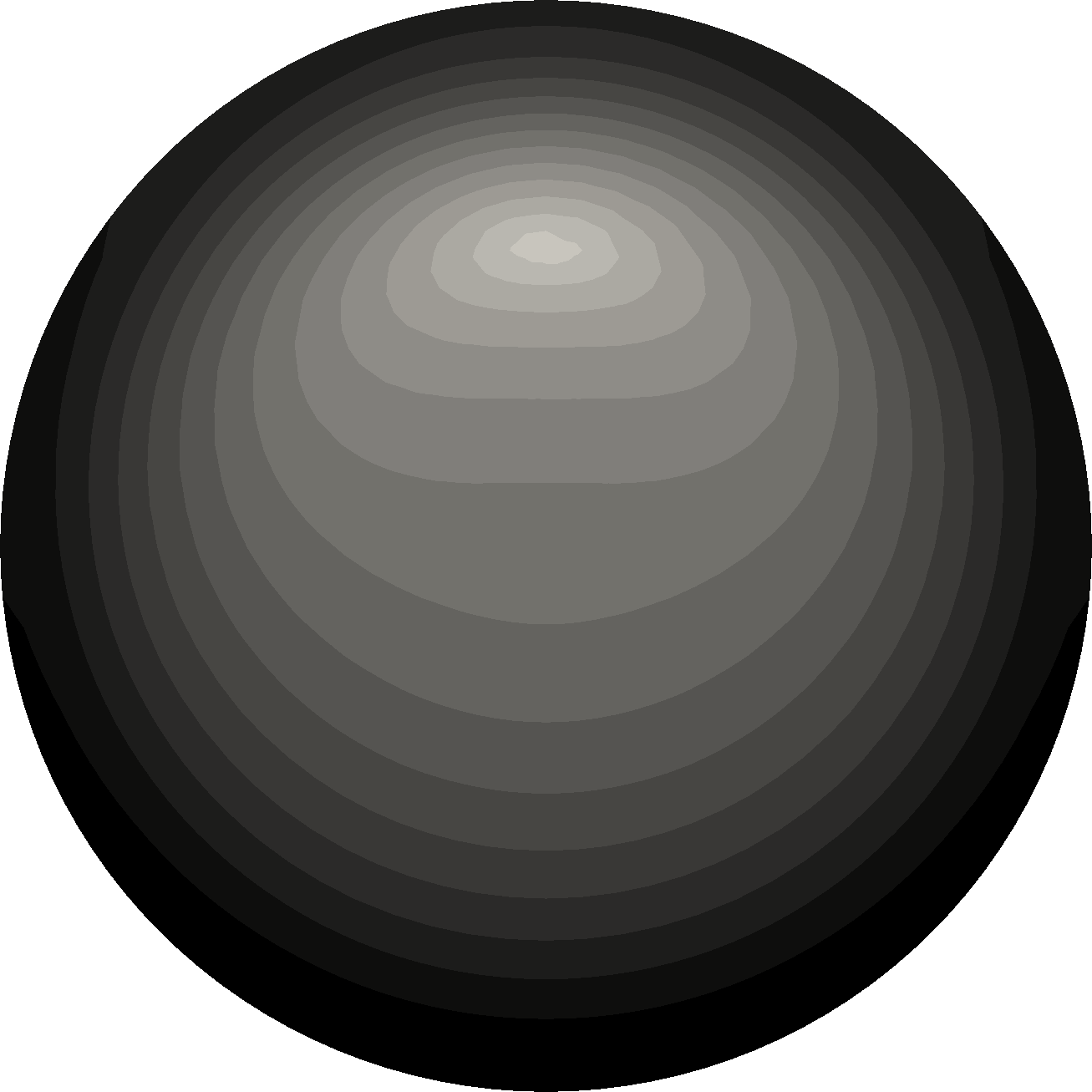} 
    \includegraphics[width=.22\columnwidth]{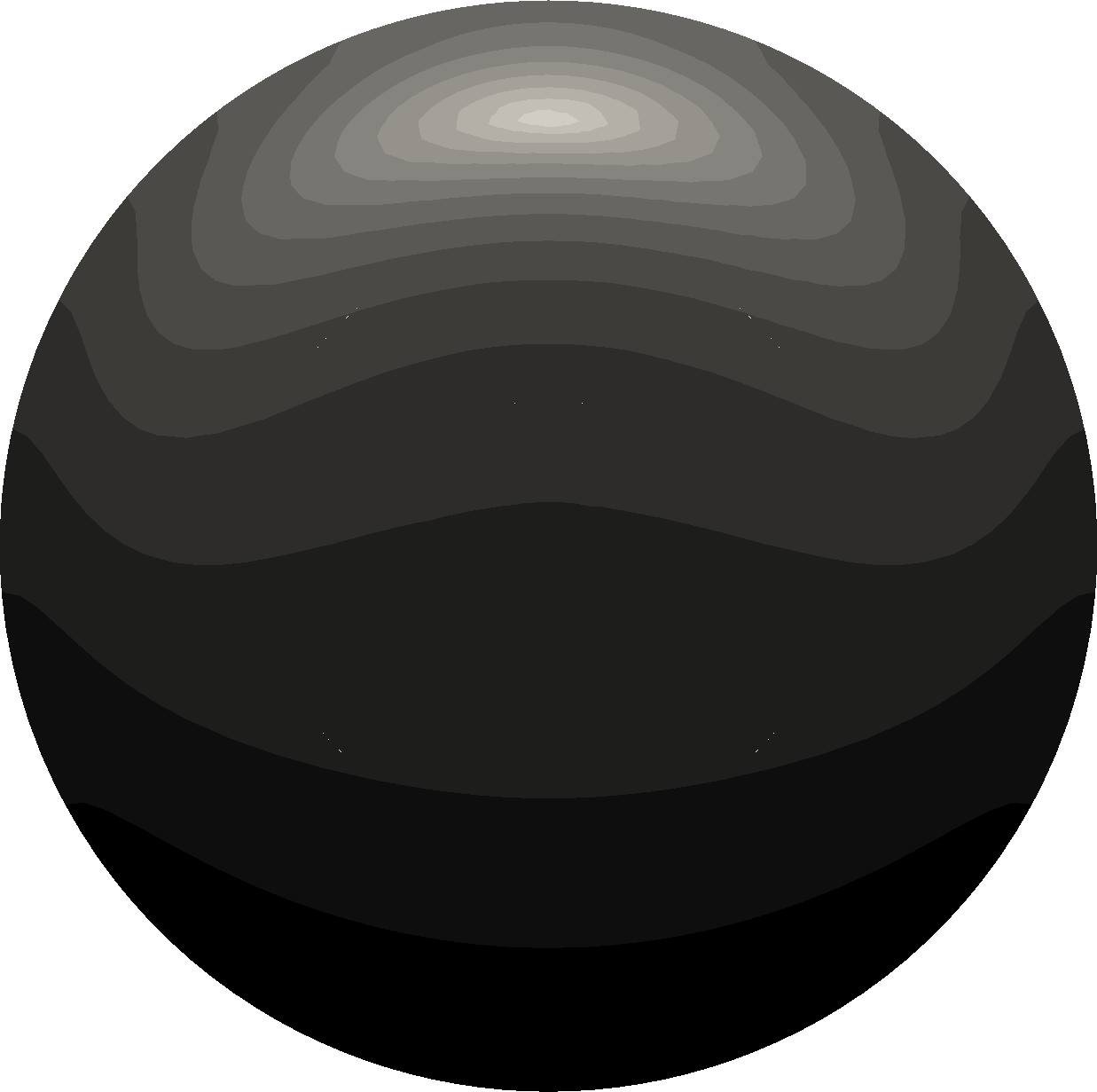} 
    \includegraphics[width=.22\columnwidth]{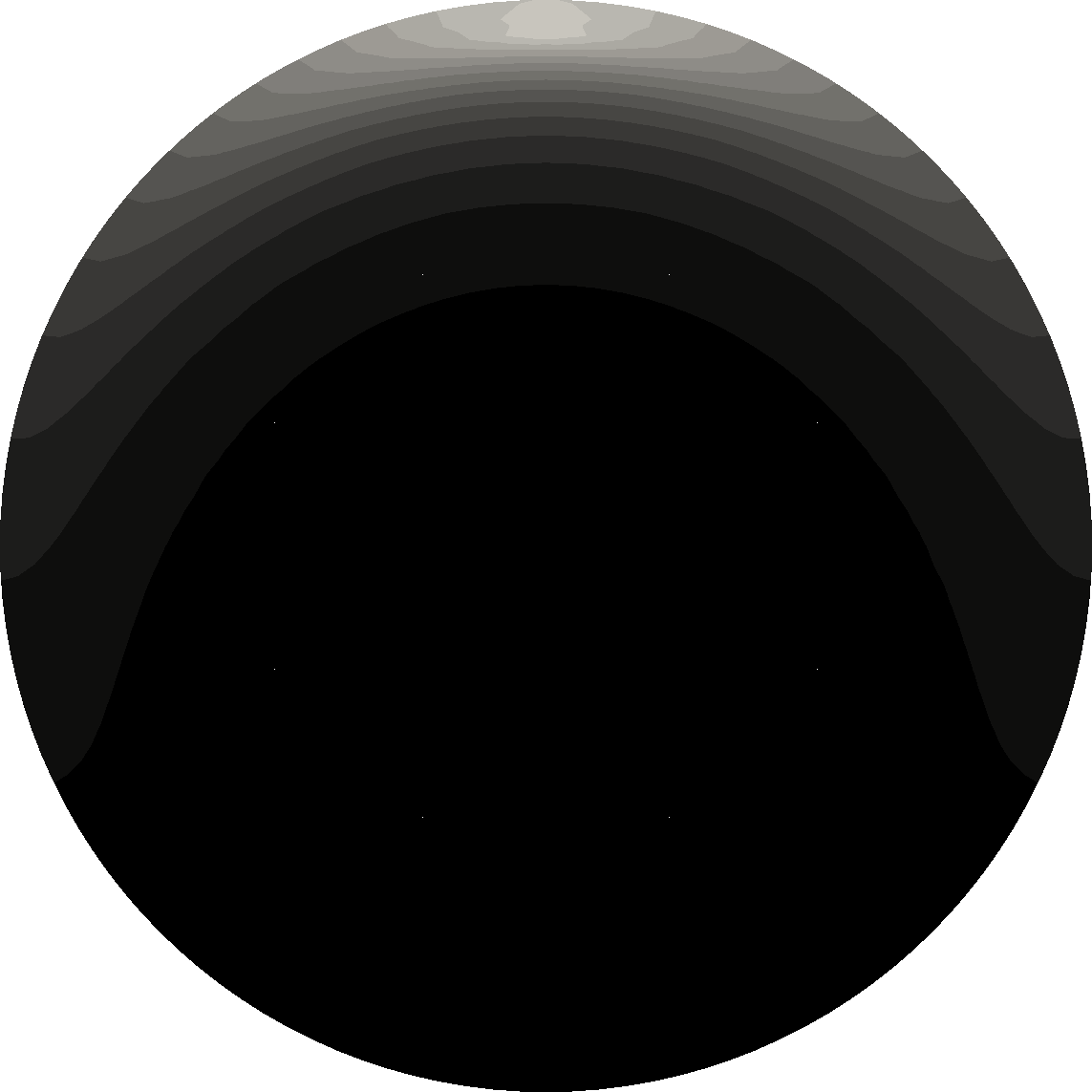} 
  \end{center}
  \caption{Contours of Green's functions at points in different depths
    for the precision operator of our prior $\mc{A}^2$. The contours are
    shown in slices through the Earth that contain the points, and
    larger values of the Green's function correspond to brighter
    shades of gray.  These Green's functions correspond directly to
    the covariance function $c(\x,\y)$ as discussed in
    Section~\ref{sec:bayesian}. Note the anisotropy for points closer
    to the surface.}\label{fig:greens_fns}
\end{figure}

\begin{figure}
  \begin{center}
    \hfill
    \begin{minipage}{0.2\columnwidth}
      \centering
      $\pm 2 \sigma$ fields \\[.5ex]
      \includegraphics[width=.96\columnwidth]{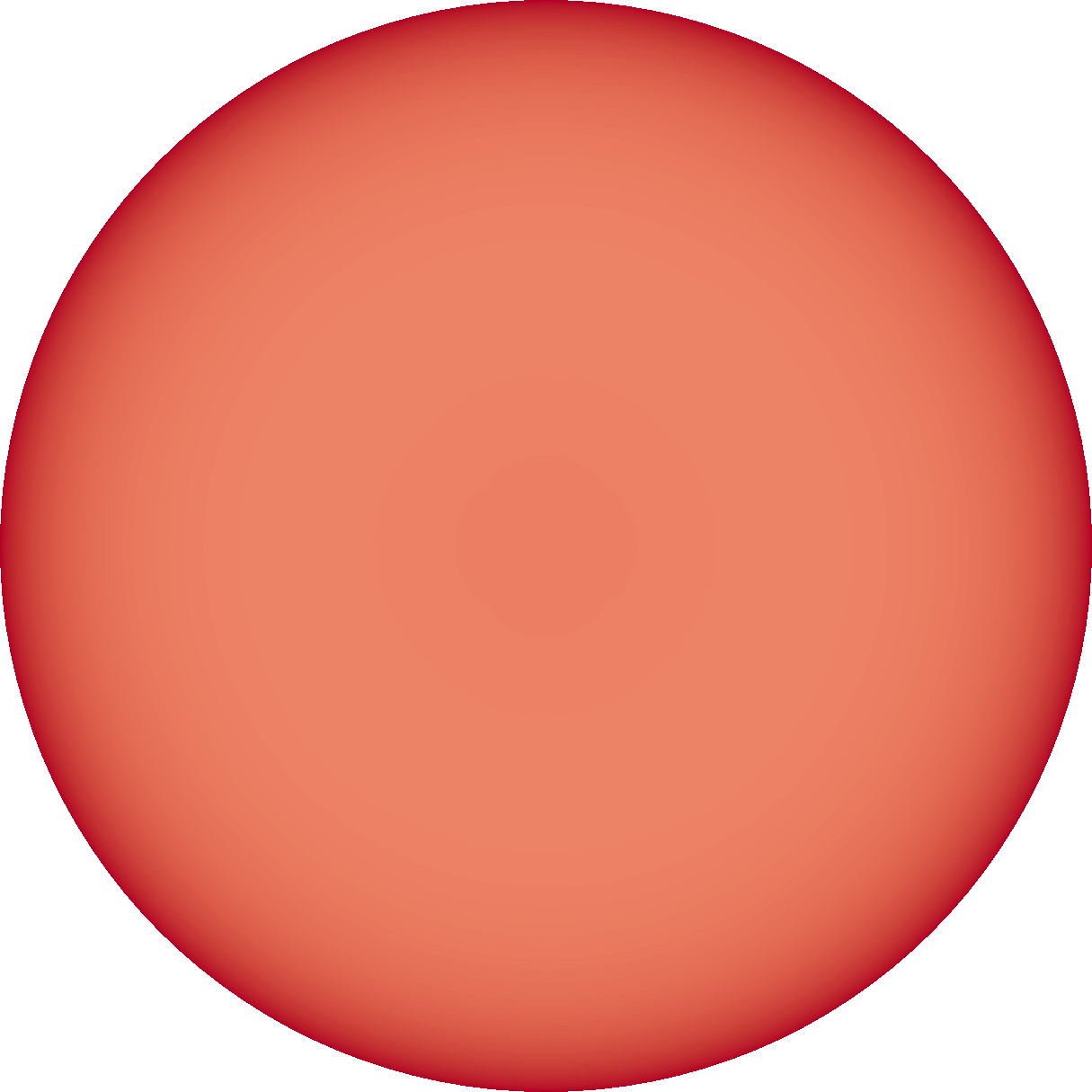} \\
      \includegraphics[width=.96\columnwidth]{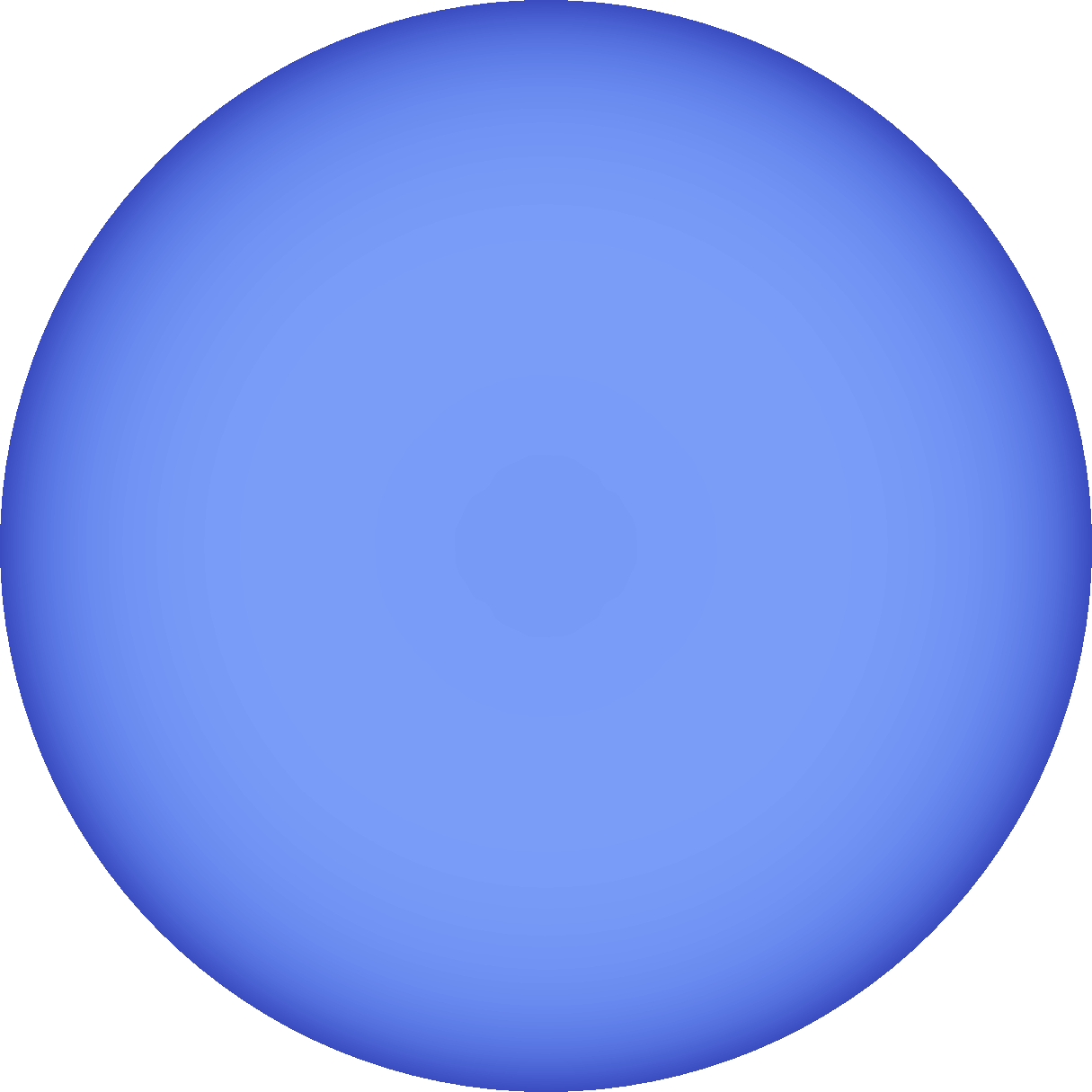}
    \end{minipage}
    \hfill
    \begin{minipage}{0.4\columnwidth}
      \centering
      prior samples \\[.5ex]
      \includegraphics[width=.48\columnwidth]{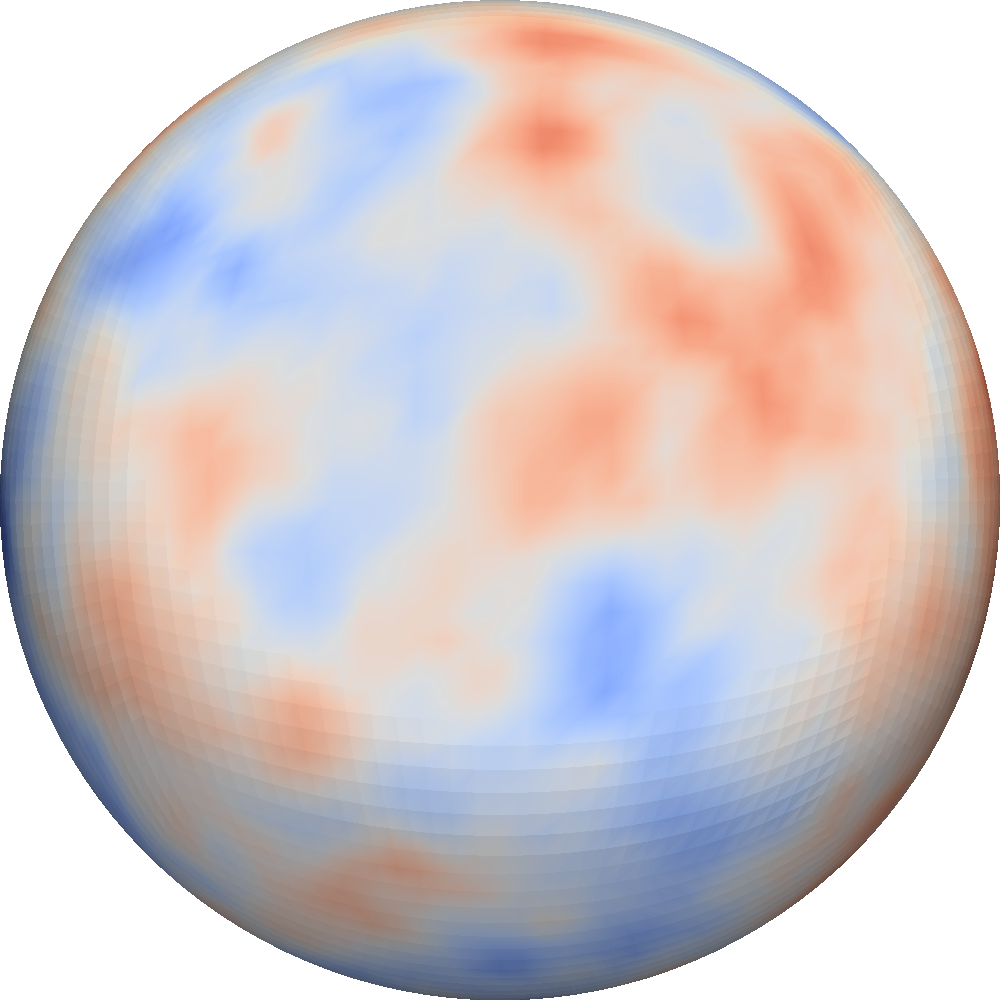}\hspace{.3ex}
      \includegraphics[width=.48\columnwidth]{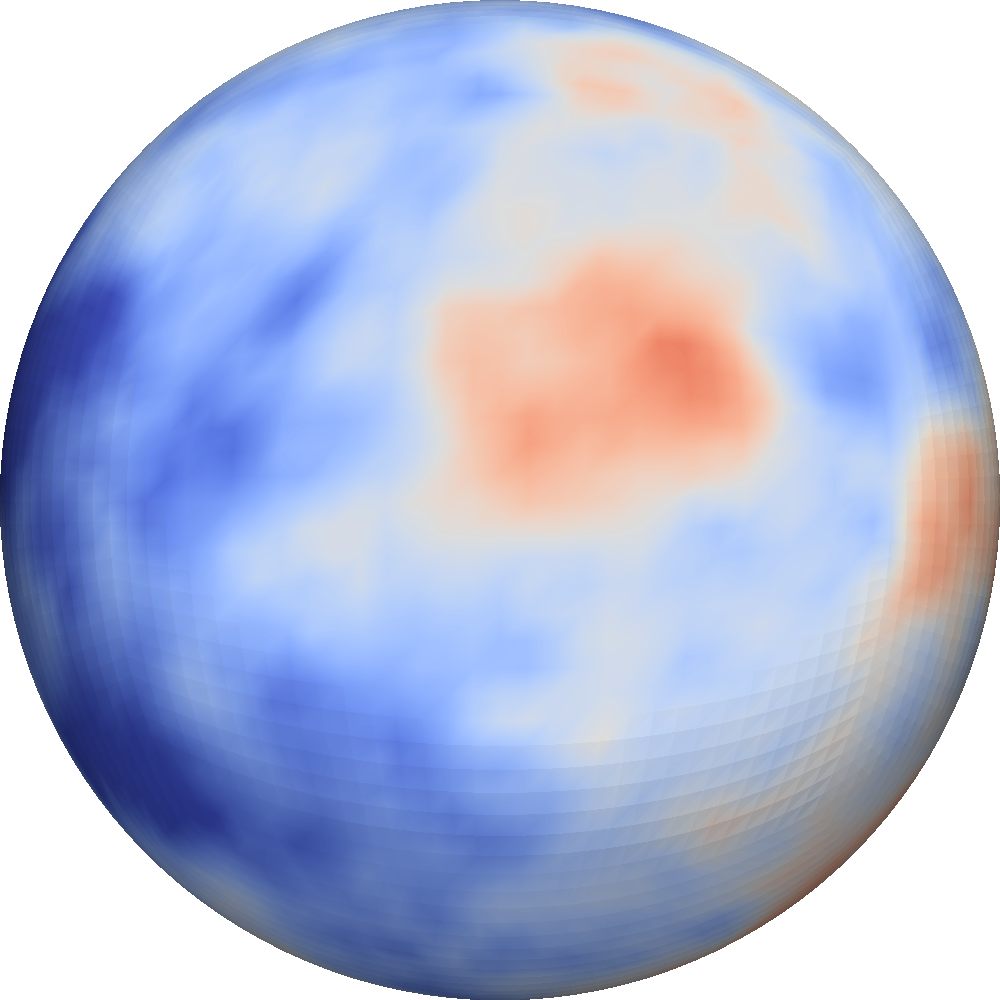}\\
      \includegraphics[width=.48\columnwidth]{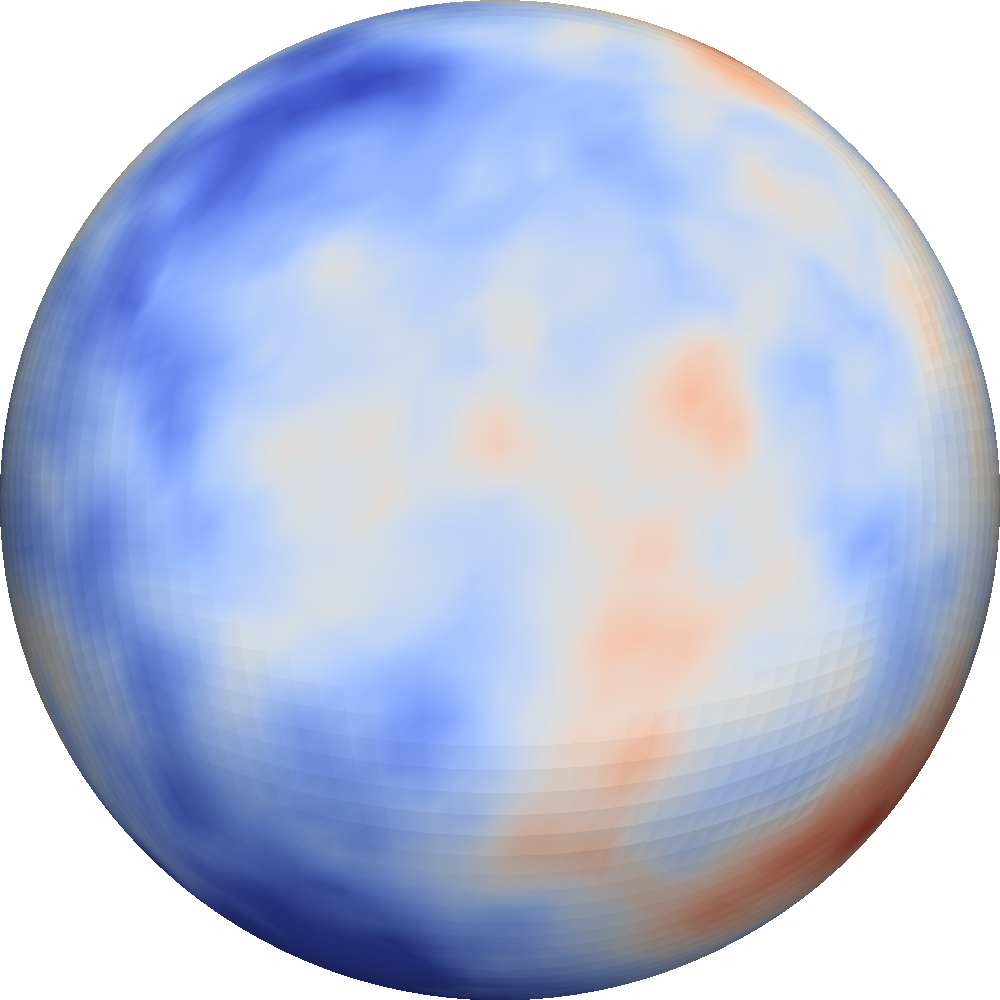}\hspace{.3ex}
      \includegraphics[width=.48\columnwidth]{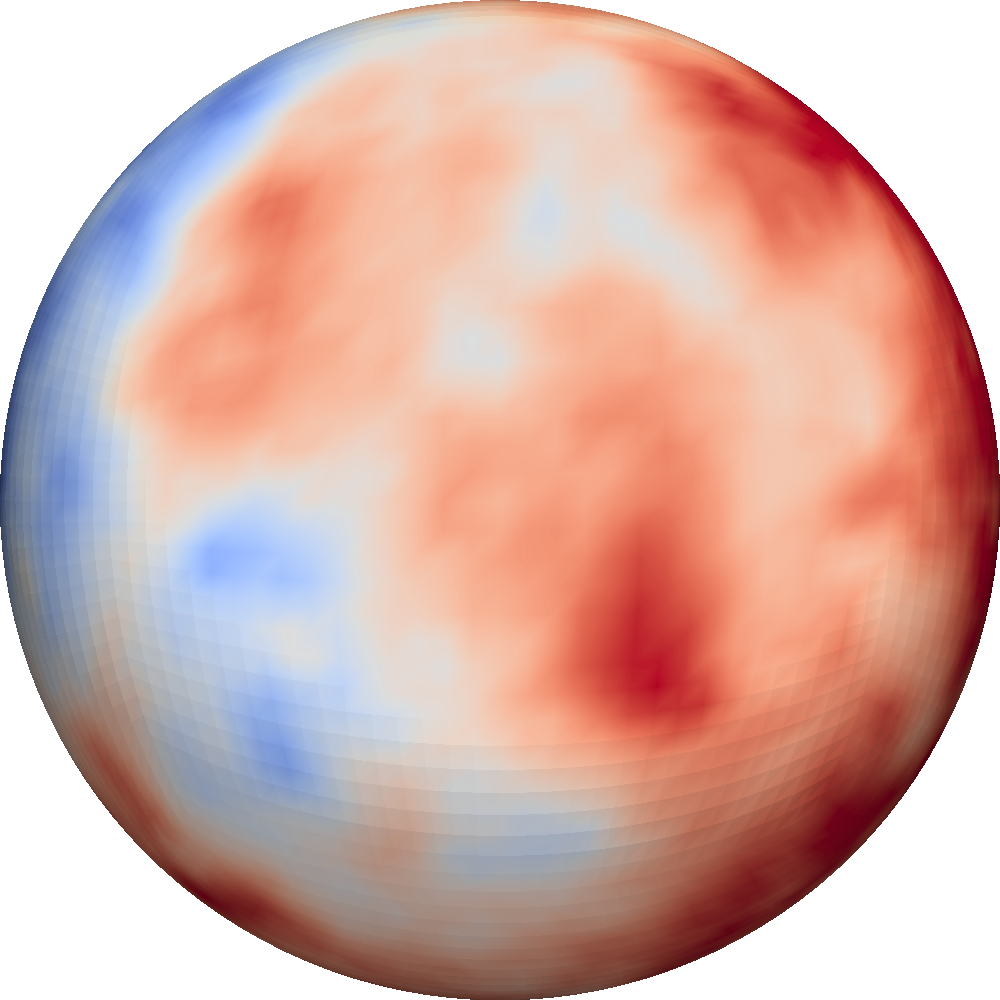}
    \end{minipage}
    \begin{minipage}{0.2\columnwidth}
      \centering
      ground truth \\[.5ex]
      \includegraphics[width=.96\columnwidth]{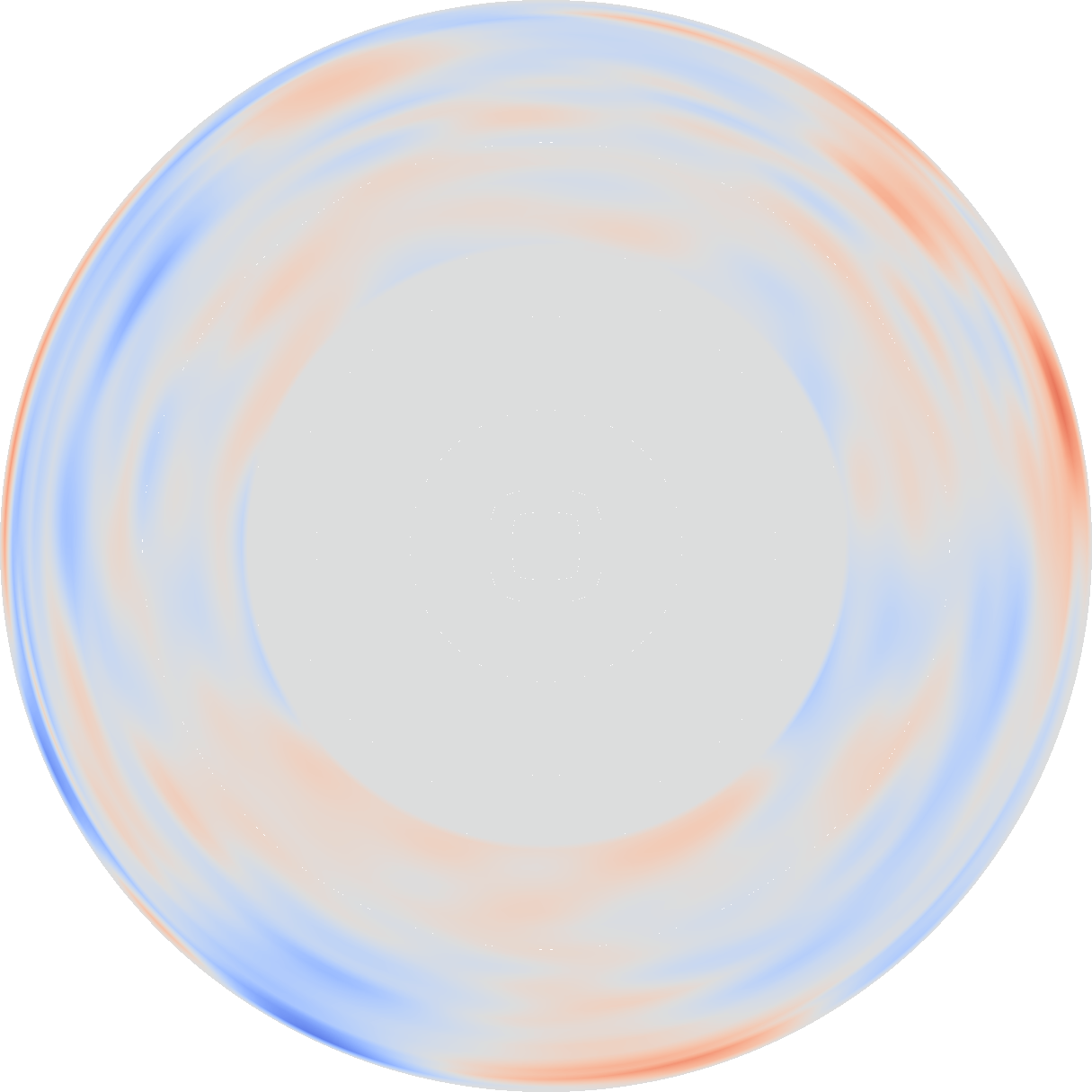}
    \end{minipage}
    \hfill
    \begin{minipage}{0.16\columnwidth}
      \includegraphics[width=.8\columnwidth]{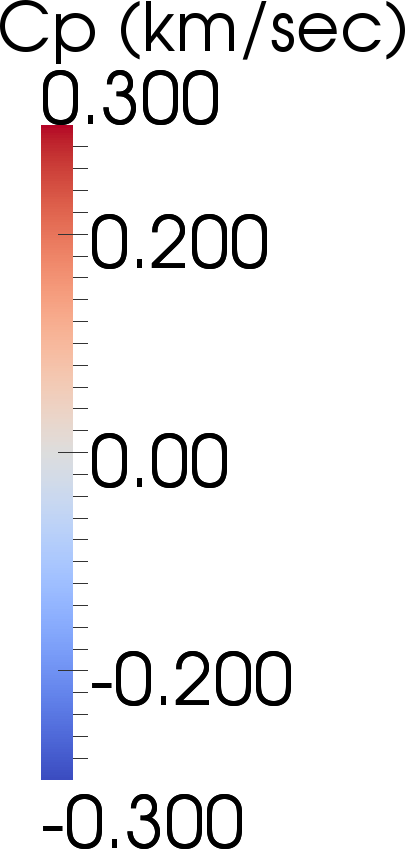}
    \end{minipage}
  \end{center}
  \caption{For illustration, we visualize
    several depictions of the prior using a common color scale.  The
    images on the far left show slices through the pointwise positive and negative
    $2\sigma$-deviation fields, which bound the pointwise $95\%$
    credible interval.  The second and third columns show samples
    drawn from the prior distribution, while the fourth column depicts
    the ``ground truth'' parameter field.  The prior has been chosen so that
    samples display similar qualitative features to the ``ground
    truth'' medium; they exhibit anisotropy in the outer layers of the
    mantle with larger correlation lengths in the lateral directions,
    and become more isotropic with higher overall correlation at
    depth.}\label{fig:prior}
\end{figure}

From the preceding observations and discussion, we choose $\alpha=1.5\cdot
10^{-2}$, while $\Theta$
is chosen to have the following form
\begin{equation}\label{eq:Theta}
\Theta = \beta\left(\bs I_3 - \theta(\xx)\xx\xx^T\right) \quad \text{with } \theta(\xx):=
\begin{cases}
  \frac{1-\theta}{r\|\xx\|^2}\left(2\|\xx\| - \frac 1r\|\xx\|^2\right) & \text{if } \|\xx\|\not=0\\
  0 & \text{if } \|\xx\|=0,
  \end{cases}
\end{equation}
where $\bs I_3$ is the $3\times 3$ identity matrix, $r=6,371$km is the earth radius,
$\beta = 7.5\cdot 10^{-3}r^2$, and
$\theta=4\cdot 10^{-2}$. The choice
$0<\theta<1$ introduces anisotropy in $\Theta$ such that the prior
assumes longer correlation lengths in tangential than in radial
directions, and the anisotropy decreases smoothly towards the center of
the sphere. In Figure~\ref{fig:greens_fns} we show several Green's functions
for the precision operator $\mc{A}^2$, which illustrate this anisotropy.
Figure \ref{fig:prior} shows a slice through the $\pm2\sigma$ fields,
through samples from the prior and through the ground truth model,
which is used to generate the synthetic seismograms.  Note that close
to the boundary, the standard deviation of the prior becomes
larger. This is partly a results of the anisotropy in the differential
operator used in the construction of the prior, but mainly an effect
of the homogeneous Neumann boundary condition used in the construction of the
square root of the prior.  This larger variance close to the boundary
is also reflected in the prior samples, which have their largest
values close to the boundary.  Note that these samples have a larger
correlation length in tangential than in normal directions, as intended
by the choice of the anisotropy in \eqref{eq:Theta}.  The ground truth
model, which is also shown in Figure \ref{fig:prior}, is comparable to
realizations of the prior in terms of magnitude as well as
correlation.

\subsection{The likelihood}
\label{sec:wave_map}

In this section, we construct the likelihood \eqref{eq:likelihoodInf} 
for the inverse acoustic wave
problem. In order to do this, we need to construct
 the parameter-to-observable map
$\bs f\LRp{\ipar}$ and the observations $\yobs$. 
Let us start by
considering the  acoustic wave equation written in the
first-order form,
\begin{subequations}\label{eq:acoustic}
  \begin{align}
    \rho  \frac{\partial\vv}{\partial t} - \nabla (\rho c^2 e) &= \g,\label{eq:acoustic1}\\
    \frac{\partial e}{\partial t} - \nabla \cdot \vv &=0,\label{eq:acoustic2}
  \end{align}
  where $\rho=\rho(\xx)$ denotes the mass density, $c=c(\xx)$ the
  local acoustic wave speed, $\g(\xx,t)$ a (smoothed) point source
  $\xx\in\Omega$, $\vv(\xx,t)$ the velocity, and $e(\xx,t)$ the trace of the strain
  tensor, i.e., the dilatation. We equip \eqref{eq:acoustic} with the
  initial conditions
  \begin{align}\label{eq:acoustic3}
    e(\xx,0) = e_0(\xx), \text{ and } \:\vv(\xx,0) = \vv_0(\xx), \quad  \xx\in \Omega,
  \end{align}
  together with the  boundary condition
  \begin{align}
    e(\xx,t) = 0,\quad \xx\in \Gamma=\partial\Omega, t\in (0,T).
  \end{align}
\end{subequations}
Here, the acoustic wave initial-boundary value problem
\eqref{eq:acoustic} is a simplified mathematical model for seismic
waves propagation in the earth \cite{AkiRichards02}.  The choice of
strain dilatation $e$ together with velocity $\vv$ in the first order
system formulation is motivated from the strain-velocity formulation
for the elastic wave equation used in
\cite{WilcoxStadlerBursteddeEtAl10}. 

Our goal is to quantify the uncertainty in inferring the spatially
varying wave speed $\ipar = c(\xx)$ from waveforms observed at
receiver locations. To define the parameter-to-observable map for a
given wave speed $c(\xx)$, we first solve the acoustic wave equation
\eqref{eq:acoustic} given $c$, 
and record the velocity $\vv$ at a finite
number of receivers in the time interval $(0,T)$. Finally, we 
compute the Fourier coefficients of the seismograms and truncate them;
the truncated coefficients are the observables in the map  $\bs
f\LRp{\ipar}$. A similar procedure is used to generate synthetic
seismograms to define $\yobs$. The noise covariance matrix 
$\matrix\Gamma^{-1}_\text{noise}$ is prescribed as a diagonal matrix
with constant variance. 

\subsection{Gradient and Hessian of the negative log posterior}
\label{sec:wave_gradhessian}
Our proposed method for uncertainty quantification in
Section~\ref{sec:finite_bayesian} requires the computation of
derivatives of the negative log posterior, which in turn requires
the gradient and Hessian of the likelihood and the prior. These
derivatives can be computed efficiently using an adjoint method, as we
now show. For clarity, we derive the gradient and action of the
Hessian  in an infinite-dimensional setting.
Let us begin by denoting  $\vv(c)$ as the space-time solution of the wave equation given the
wave speed $c=c(\xx)$, and
$\mathcal B$ as the observation operator. The parameter-to-observable map $\bs
 f(c)$ can be written as $\mathcal B\vv(c)$. Thus, the negative log posterior
 is (compare with \eqref{eq:optfun})
\begin{equation}
\mc{J}\LRp{c} := \half
\nor{{\mathcal B}\vv(c) - \yobs}_{\matrix\Gamma^{-1}_\text{noise}}^2 +
\half \nor{\mc{A}(c-c_0)}_{L^2\LRp{\D}}^2,
\end{equation}
where $\matrix\Gamma_\text{noise}$ is specified in Section~\ref{sec:wave_problemsetup}.
The dependence on the wave speed $c$ of the velocity $\vv$ and
dilatation $e$ is given by the solution of the forward wave
propagation equation \eqref{eq:acoustic}.
The adjoint approach \cite{EpanomeritakisAkccelikGhattasEtAl08} allows
us to write $\mathcal G(c)$, the gradient of $\mc{J}$ at a point $c$
in parameter space, as
\begin{equation}\label{eq:grad}
\mathcal{G}(c) := 2\rho c\int^T_0  e (\nabla\cdot \w) \, dt + \Acal^2(c-c_0),
\end{equation}
where the adjoint velocity $\w$ and adjoint strain dilatation
$d$ satisfy the {\em adjoint wave propagation terminal-boundary value problem}
\begin{subequations}\label{eq:adjoint}
\begin{alignat}{2}
-\rho \frac{\partial\w}{\partial t} + \bs{\nabla}(c^2\rho
d) &= - \mathcal{B}^*\matrix\Gamma^{-1}_{\text{noise}}(\mathcal{B}\vv-\yobs)  &&\quad\text{in }
\Omega \times (0,T), \\
-\frac{\partial d}{\partial t} + {\nabla}\cdot\w
&= 0 && \quad\text{in } \Omega \times (0,T), \\
\rho \w = \bs{0}, d &= 0  && \quad\text{in } \Omega \times
\left\{t=T\right\}, \\
d &= 0  && \quad\text{on } \Gamma \times
(0,T).
\end{alignat}
\end{subequations}
Here, $\mathcal{B}^*$, an operator from $\R^q$ to the space-time
cylinder $\Omega \times (0,T)$, is the adjoint of $\mathcal{B}$.  Note
that the adjoint wave equations must be solved backward in time (due
to final time data) and have the data misfit as a source term, but
otherwise resemble the forward wave equations.

Similar to the computation of the gradient, the
Hessian operator of $\mc{J}$ at $c$ acting on an arbitrary variation $\tCC$
is given by
\begin{align}\label{eq:wave_hessian}
\mathcal{H}(\CC)\tCC &:= 2\rho \int^T_0 ce(\nabla\cdot\tilde \w) +
c\tilde e(\nabla\cdot \w) + \tilde c e(\nabla\cdot \w)\,dt + \Acal^2\tilde c,
\end{align}
where  $\tilde{\bs{v}}$ and  $\tilde{e}$ satisfy the
{\em incremental forward wave propagation initial-boundary value problem}
\begin{alignat*}{2}
\rho \frac{\partial\tilde\vv}{\partial t} - \nabla(\rho c^2 \tilde e) &= \nabla(2 \rho c \tilde c e)   &&\quad\text{in }
\Omega \times (0,T), \\
\frac{\partial \tilde e}{\partial t} - \nabla\cdot\tilde\vv
&= 0 && \quad\text{in } \Omega \times (0,T), \\
\rho \tilde\vv =\bs{0}, \tilde e &= 0  && \quad\text{in } \Omega \times
\left\{t=0\right\}, \\
\tilde e &= 0  && \quad\text{on } \Gamma \times
(0,T).
\end{alignat*}
On the other hand, $\tilde\w$ and $\tilde d$
satisfy the {\em incremental adjoint wave propagation
  terminal-boundary value problem}
\begin{alignat*}{2}
-\rho \frac{\partial\tilde \w}{\partial t} + {\nabla}(c^2\rho
\tilde d) &=  -{\nabla}(2\tilde c c\rho d)-
\mathcal{B}^*\matrix\Gamma^{-1}_{\text{noise}}\mathcal{B}\tilde \vv
  &&\quad\text{in }
\Omega \times (0,T), \\
-\frac{\partial \tilde d}{\partial t} + {\nabla}\cdot\tilde \w
&= 0 && \quad\text{in } \Omega \times (0,T), \\
\rho \tilde\w = \bs{0}, \tilde d &= 0  && \quad\text{in } \Omega \times
\left\{t=T\right\}, \\
\tilde d &= 0  && \quad\text{on } \Gamma \times
(0,T).
\end{alignat*}
As can be seen, the incremental forward and incremental adjoint wave
equations are linearizations of their forward and adjoint
counterparts, and thus differ only in the source terms.  Moreover, we
observe that the computation of the gradient and the Hessian action
amounts to solving a pair of forward/adjoint and a pair of
incremental-forward/incremental-adjoint wave equations, respectively.
For our computations, we use the Gauss-Newton approximation of the
Hessian, which is guaranteed to be positive. This amounts to
neglecting the terms that contain $\nabla\cdot \w$ in
\eqref{eq:wave_hessian}, and neglecting the term that includes $d$ in
the incremental adjoint wave equations. 

\subsection{Discretization of the wave equation and implementation details}
\label{sec:wave_implementation}

We use the same hexahedral mesh to compute the wave solution
$\LRp{\vv,e}$ as is used for the parameter $c$. While the parameter is
discretized using trilinear finite elements, the wave equation, and
its three variants (the adjoint, the incremental forward, and the
incremental adjoint), are solved using a high-order discontinuous
Galerkin (dG) method. The method, for which details are provided in
\cite{WilcoxStadlerBursteddeEtAl10, Bui-ThanhGhattas12b}, supports
$hp$-non-conforming discretization, but only $h$-non-conformity is
used in our implementation. For efficiency and scalability, a tensor
product of Lagrange polynomials of degree $N$ (we use $N\in\{2,3,4\}$
for the examples in the next section) is employed together with a
collocation method based on Legendre-Gauss-Lobatto (LGL) nodes. As a
result, the mass matrix is diagonal, which facilitates time
integration using the classical four-stage fourth-order Runge Kutta
method. We equip our dG method with exact Riemann numerical fluxes at
element faces. To treat the non-conformity, we use the mortar
approach of Kopriva \cite{Kopriva96,KoprivaWoodruffHussaini02} to
replace non-conforming faces by mortars that connect pairs of
contributing elements. The actual computations are performed on the
mortars instead of the non-conforming faces, and the results are then
projected onto the contributing element faces. The method has been
shown to be consistent, stable, convergent with optimal order, and
highly scalable
\cite{Bui-ThanhGhattas12b,WilcoxStadlerBursteddeEtAl10}.

It should be pointed out that the discretizations of the gradient and
Hessian action given in Section~\ref{sec:wave_gradhessian} are not
consistent with the discrete gradient and Hessian-vector product
obtained by first discretizing the negative log posterior and then
differentiating it. Here, inconsistency means that the former are
equivalent to the latter only in the limit as the mesh size approaches
zero (see also \cite{Gunzburger03,HinzePinnauUlbrichEtAl09}). The
reason is that additional jump terms due to numerical fluxes at
element interfaces are introduced in the discontinuous Galerkin
discretization of the wave equation. In our implementation, we include
these terms to ensure consistency, and this is verified by comparing
the discretized gradient and Hessian action expressions with their
finite difference approximations.

Moreover, since we use a continuous Galerkin finite element method for
the parameter, but a discontinuous Galerkin method for the wave solution, it is
necessary to prolongate the parameter to the solution space before solving the
forward wave equation, and its variants (adjoint, incremental state,
incremental adjoint). Conversely, the gradient and the
Hessian-vector application are computed in the wave solution space,
and then restricted to
the parameter space to provide the correct derivatives for the
optimization solver. To ensure the symmetry of the Hessian, we
construct these restriction and prolongation operations such that they
are adjoint of each other.

Our discretization approach for the Bayesian inverse problem in
\secref{finite_bayesian} requires the repeated application of $\bs
A^{-1}$, each amounting to an elliptic PDE solve on the
finite-dimensional parameter space. To accomplish this task
efficiently, we use the parallel algebraic multigrid (AMG) solver {\em
  ML} from the Trilinos project \cite{GeeSiefertHuEtAl06}. The cost of
this elliptic solve is negligible compared to that of solving the
time-dependent seismic wave equations, which employ high order
discretization in contrast to the trilinear discretization of the
anisotropic Poisson operator, $\bs A$.

The adjoint equation has to be solved backwards-in-time (as shown in
Section~\ref{sec:wave_gradhessian}); computation of the gradient
\eqref{eq:grad} requires combinations of the state and adjoint
solutions corresponding to same time.  Thus the gradient computation
requires the complete time history of the forward solve, which cannot
be stored due to the large-scale nature of our problem. A similar, but
slightly more challenging storage problem occurs in the Hessian-vector
application. Here, solving the incremental state equation requires the
solution of the state equation, and the incremental adjoint solution
requires the solution of the incremental state equation.  We avoid
storage of the time history of these wavefields by using a
checkpointing method as employed in
\cite{EpanomeritakisAkccelikGhattasEtAl08}. This scheme reduces the
necessary storage at the expense of increasing the number of wave
propagation solves. 

Between 1200 and 4096 processor cores\footnote{These computations were
  performed on the Texas Advanced Computing Center's Lonestar 4
  system, which has 22,656 Westmere processor cores with 2GB memory
  per core.} for 10-20 hours are needed to solve the seismic inverse
problems presented in the next section. The vast majority of the
runtime is spent on computing solutions of the forward, adjoint and
incremental wave equations either for the computation of the MAP point
(see Section~\ref{sec:findMAP}) or the Lanczos iterations for
computing the low rank approximation of the misfit Hessian (see
Section~\ref{sec:wave_lowrank}). Due to the large number of required
wave propagation solves, good strong scalability of the wave
propagation solver is important for rapid turnaround.  We refer to the
discussion in \cite{Bui-ThanhBursteddeGhattasEtAl12} on the
scalability of the wave propagation solver, as well to the overall
scalability of our Bayesian inversion approach applied to seismic
inverse problems of up to one million parameters.

\subsection{Setup of model problems}
\label{sec:wave_problemsetup}

\begin{figure}
\begin{center}
\includegraphics[width=.35\columnwidth]{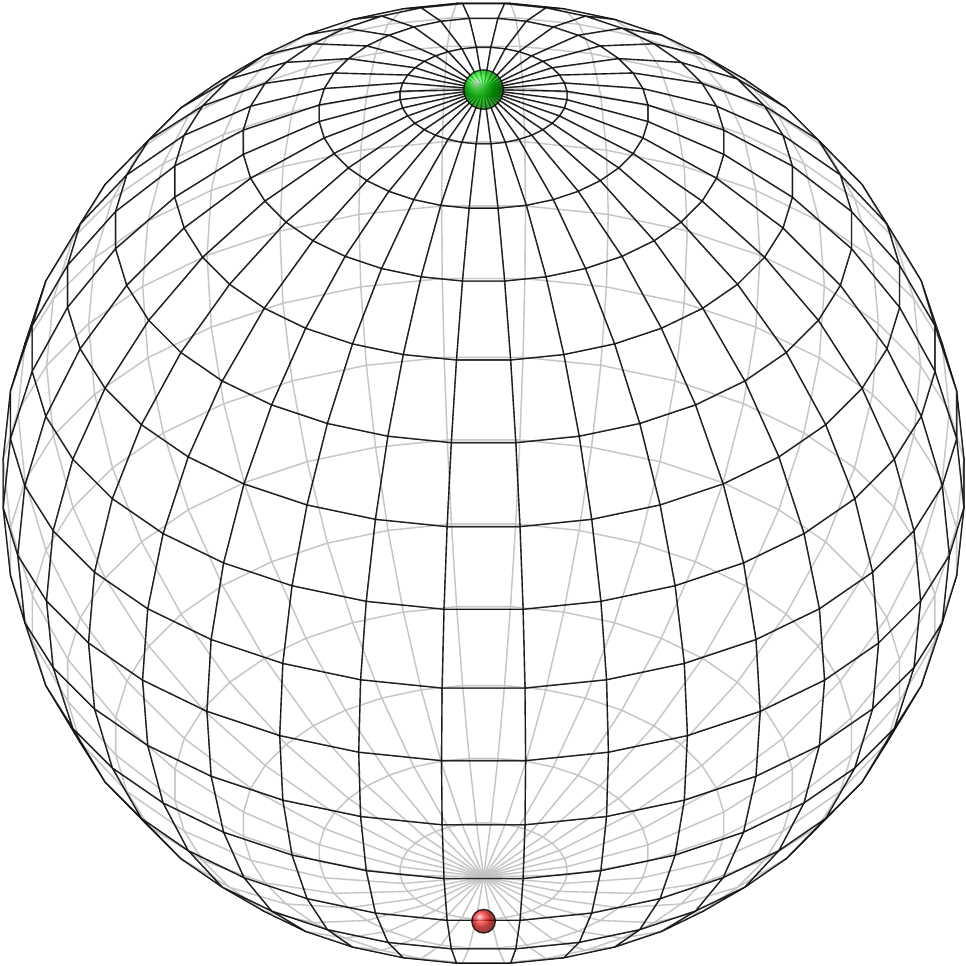}\hspace{4ex}
\includegraphics[width=.35\columnwidth]{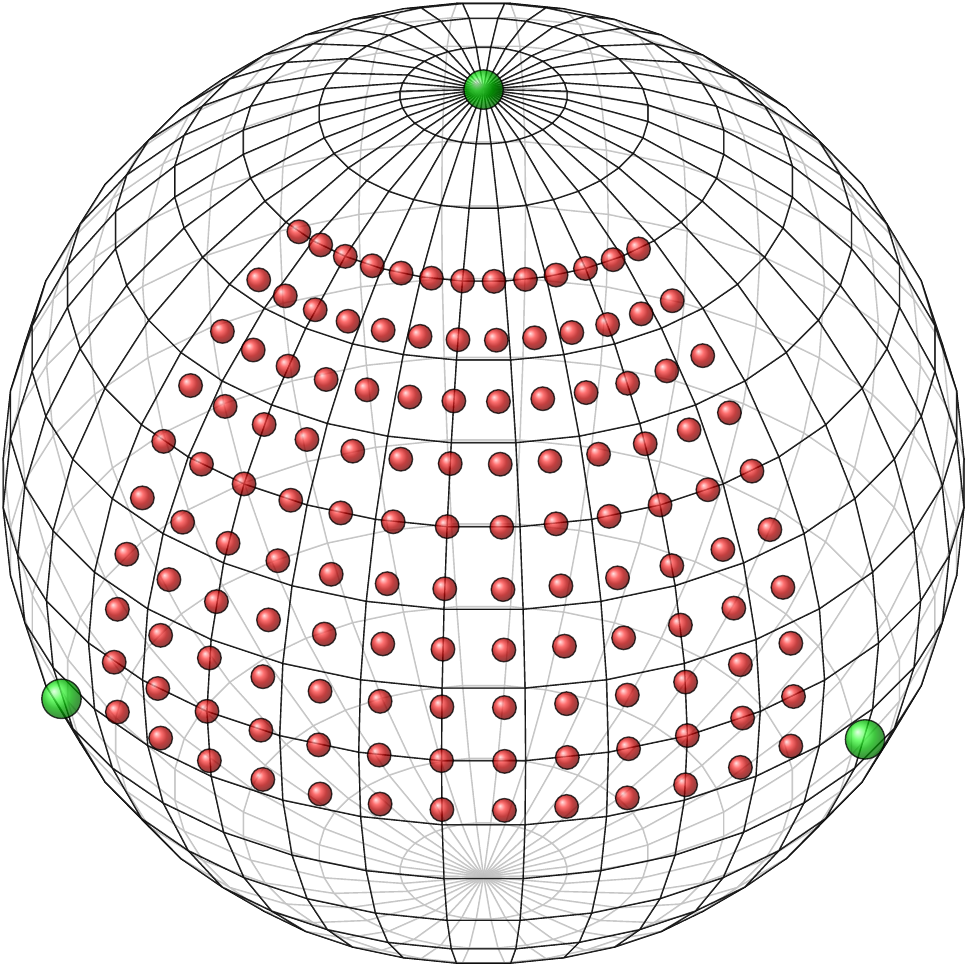}
\end{center}
\caption{Location of sources (green) and receivers (red) for Problem I
  (left) and Problem II (right).}
\label{fig:sourcesreceivers}
\end{figure}%

Synthetic observations $\yobs$ are generated from solution of the wave
equation using the S20RTS earth model. To mitigate the inverse crime
\cite{KaipioSomersalo05}, the local wave speed on LGL nodes of the
wave propagation mesh is used to generate the observations, which
implies that a higher order approximation of the earth model is used
to generate the synthetic data, but the inversion is carried out on a
lower order mesh. Both sources and receivers are located at 10km depth
from the earth surface. For the source term $\g$ in
\eqref{eq:acoustic}, we use a delta function point source in the
$z$-direction convolved with a narrow Gaussian in space. In time, we
employ a Gaussian
with standard deviation of $\sigma=20$s centered at $60$s.  The wave
propagation mesh (i.e., the discretization of velocity and dilatation)
is chosen fine enough to accurately resolve frequencies below
$0.05$Hz. We Fourier transform the (synthetic) observed velocity
waveforms at each receiver location and retain only the first 101
Fourier modes to define the observations $\yobs$. In our problems,
the Fourier coefficients $\yobs$ vary between $10^{-5}$ and
$10^{-1}$, and we choose for the noise covariance
a diagonal matrix with a standard deviation of $0.002$.

We consider the following two model problems:
\begin{itemize}
\item {\bf Problem I:} The first problem has a single source at the
  North pole and a single receiver at 45$^\circ$ south of the equator,
  as illustrated in the left image of Figure
  \ref{fig:sourcesreceivers}. The wave propagation time is 1800s. The
  wave speed (i.e., unknown material parameter) field is discretized
  on a mesh of trilinear hexahedra with 78,558 nodes, representing the
  unknowns in the inverse problem. The forward problem is discretized
  on the same mesh with 3rd-order dG elements, resulting in about 21.4
  million spatial wave propagation unknowns, and in 2100 four-stage,
  fourth-order Runge Kutta time steps.

\item {\bf Problem II:} The second problem uses 130 receivers
  distributed on a quarter of the Northern hemisphere along zonal
  lines with 7.5$^\circ$ spacing and 3 simultaneous sources as shown
  on the right of \figref{sourcesreceivers}. The wave propagation time
  is 1200s. The wave speed is discretized on three different trilinear
  hexahedral meshes with 40,842, 67,770 and 431,749 wave speed
  parameters, which represent the unknowns in the inverse
  problem. These meshes corresponding to discretizations with 4th, 3rd
  and 2nd order discontinuous elements for the wave propagation
  variables (velocity and dilatation).  The results in the next
  section were computed with 67,770 wave speed parameters and the 3rd
  order dG discretization for velocity and dilatation. This amounts to
  18.7 million spatial wave propagation unknowns, and 1248 Runge Kutta
  time steps.

\end{itemize}

\subsection{Low rank approximation of the prior-preconditioned misfit Hessian}
\label{sec:wave_lowrank}

\begin{figure}\centering
    \includegraphics[width=0.48\columnwidth]{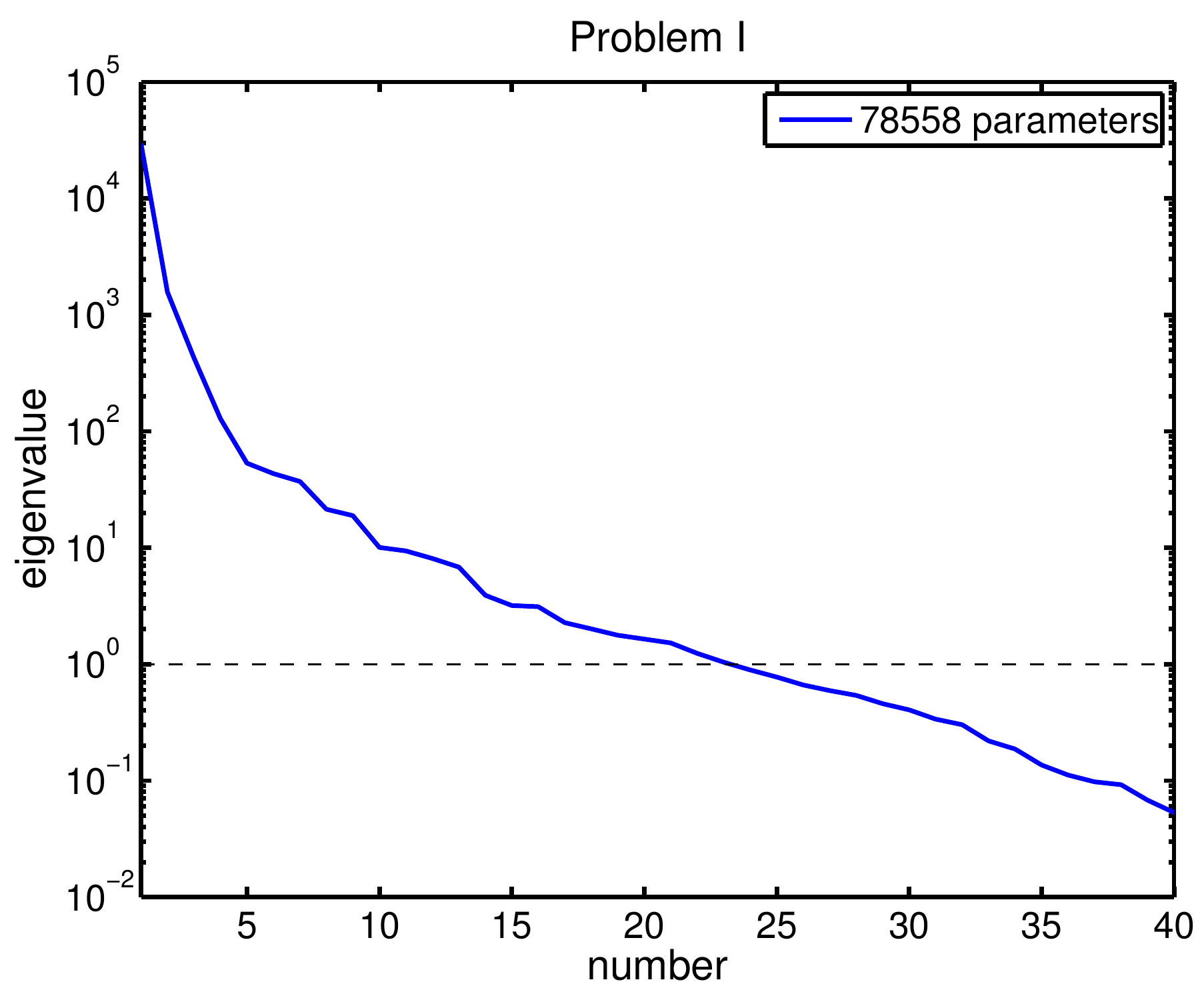}
    \hfill
    \includegraphics[width=0.48\columnwidth]{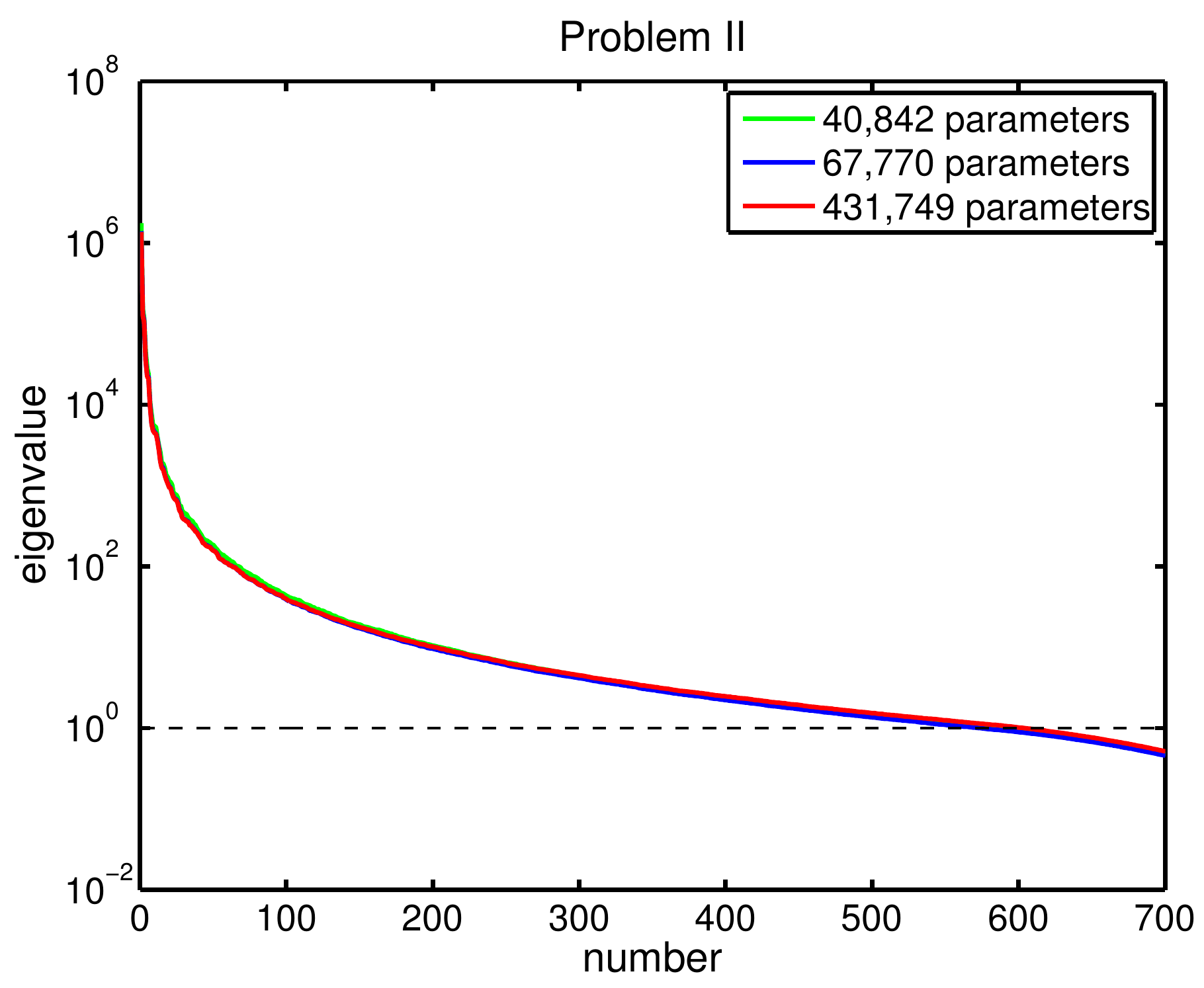}
    \caption{Logarithmic plot of the spectrum of prior-preconditioned
      data misfit Hessian for Problem I (left) and Problem II
      (right). The computation of the approximate spectrum for Problem
      I uses a discretization with $78558$ wave speed parameters,
      third-order dG finite elements for the wave propagation
      solution, and 50 Lanczos iterations. The spectrum for Problem II
      is computed on different discretizations of the parameter mesh
      using 900 Lanczos iterations. The eigenvalues for the three
      discretizations essentially lie on top of each other, which
      illustrates that the underlying infinite-dimensional statistical
      inverse problems is properly approximated. The horizontal line
      $\lambda=1$ shows the reference value for the truncation of the
      spectrum of the misfit Hessian.  For an accurate approximation
      of the posterior covariance matrix (i.e., the inverse of the
      Hessian), eigenvalues that are small compared to $1$ can be
      neglected as discussed in Section \ref{lowrank}, and in
      particular as shown in \eqref{eqn:eigen_error}.
      \label{fig:wave_spectra}}
\end{figure}

Before discussing the results for the quantification of the
uncertainty in the solution of our inverse problems, we numerically study the
spectrum of the prior-preconditioned misfit Hessian.
In Figure \ref{fig:wave_spectra}, we show the dominant spectrum of the
prior-preconditioned Hessian evaluated at the MAP estimate for Problem
I (left) and Problem II (right).  As can be observed, the eigenvalues
decay faster in the former than in the latter.  That is, the former is
more ill-posed than the latter. The reason is that the three
simultaneous source and 130 receivers of Problem II provide more
information on the earth model. 
This implies that retaining more eigenvalues is necessary to
accurately approximate the prior-preconditioned Hessian of the data
misfit for Problem II compared to Problem I. In particular, we need at
least 700 eigenvalues for Problem II as compared to about 40 for
Problem I to obtain a sensible low-rank approximation of the Hessian,
and this constitutes the bulk of computation time (since each
Hessian-vector product in the Lanczos solver requires incremental
forward and adjoint wave propagation solutions). These numbers compare
with a total number of parameters of 78,558 (Problem I) and 67,770
(Problem II), which amounts to a reduction of between  two and three
orders of magnitude. This directly translates into two to three 
orders of magnitude reduction in cost of solving the statistical
inverse problem. 

Figure \ref{fig:wave_spectra} presents the spectra for Problem II for
three different discretization of the wave speed parameter field. The
figure suggests that the dominant spectrum is essentially
mesh-independent and that all three parameter meshes are sufficiently
fine to resolve the dominant eigenvectors of the prior-preconditioned Hessian.
Consequently, the Hessian low-rank approximation, particularly the
number of Lanczos iterations, is independent of the number of discrete
parameters.  Thus, in this example, the number of wave propagation
solutions required by the low-rank approximation does not depend on
the parameter dimension.

\begin{figure}\centering
  \includegraphics[width=.18\columnwidth]
                  {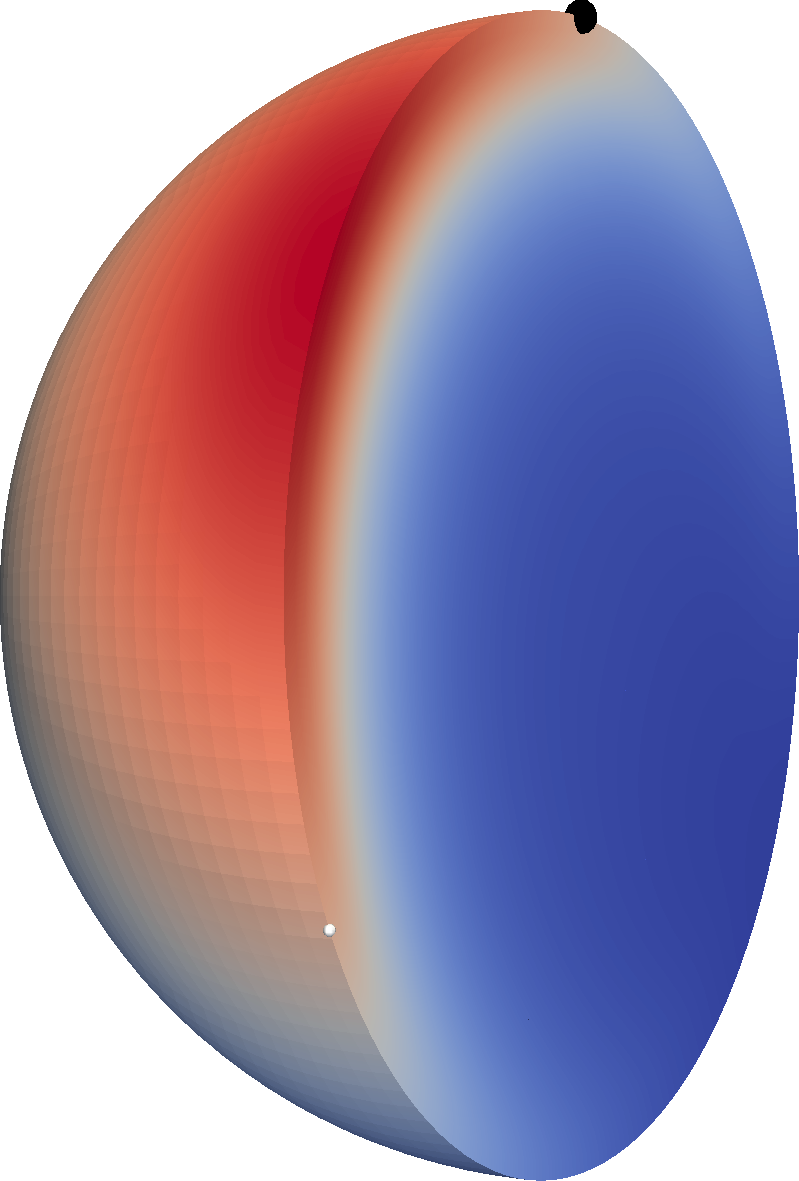}\hfill
  \includegraphics[width=.18\columnwidth]
                  {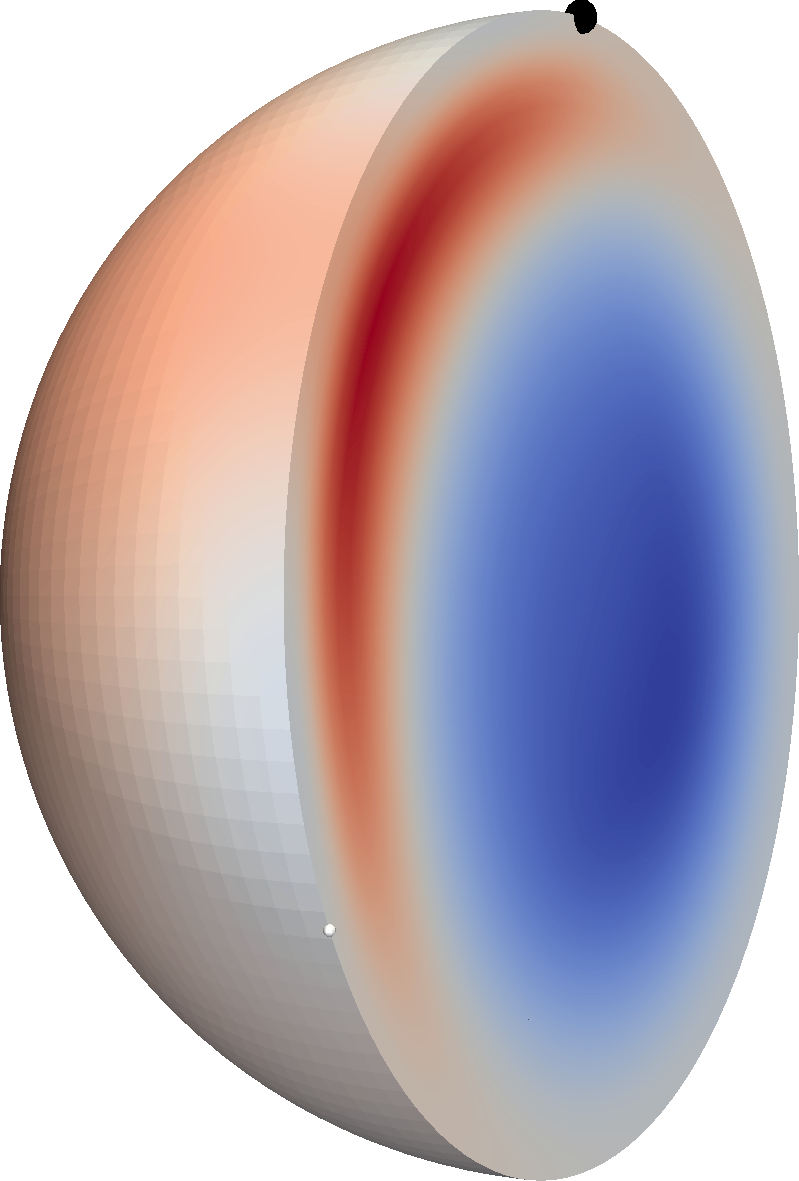}\hfill
  \includegraphics[width=.18\columnwidth]
                  {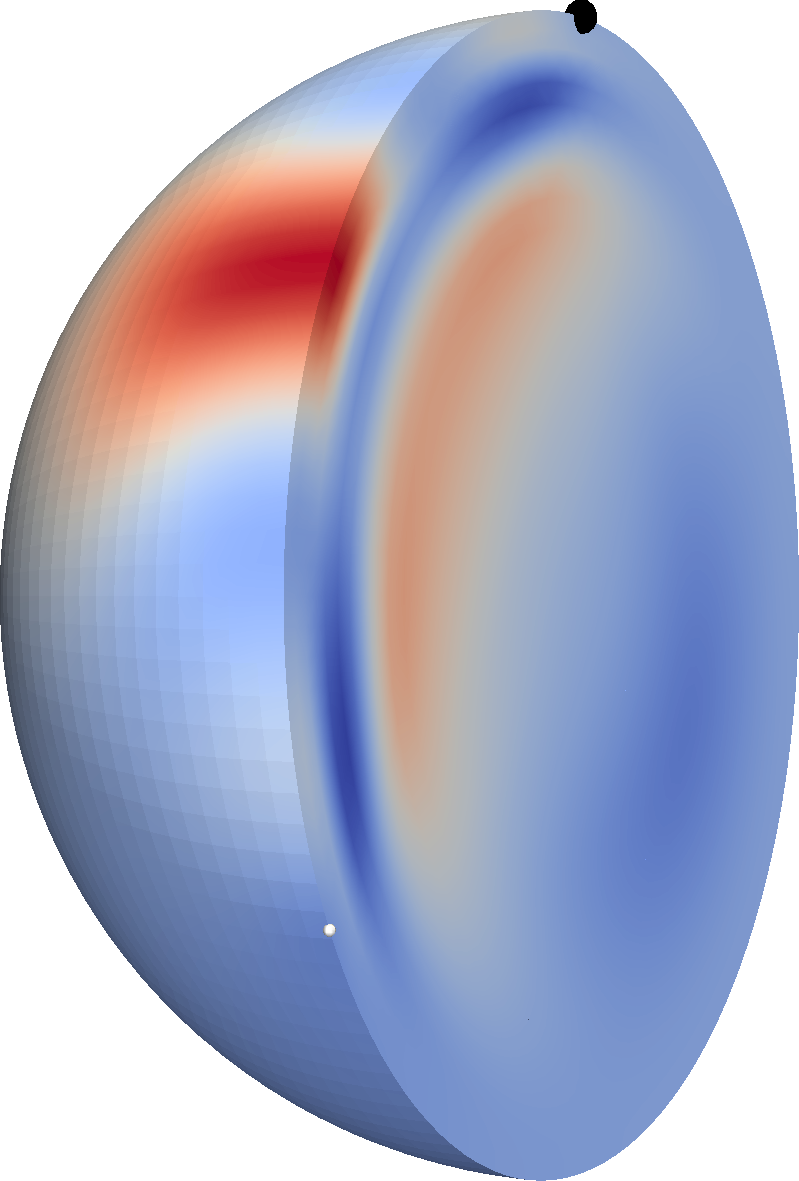}\hfill
  \includegraphics[width=.18\columnwidth]
                  {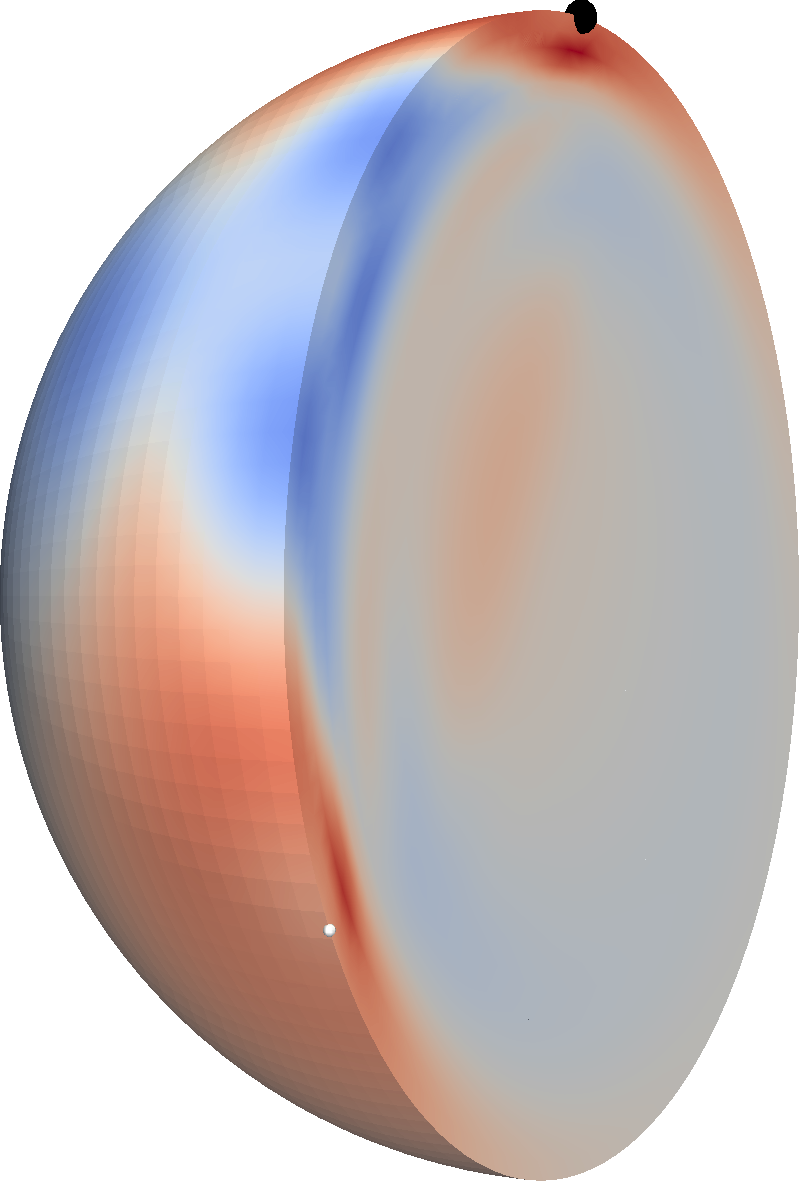}\hfill
  \includegraphics[width=.18\columnwidth]
                  {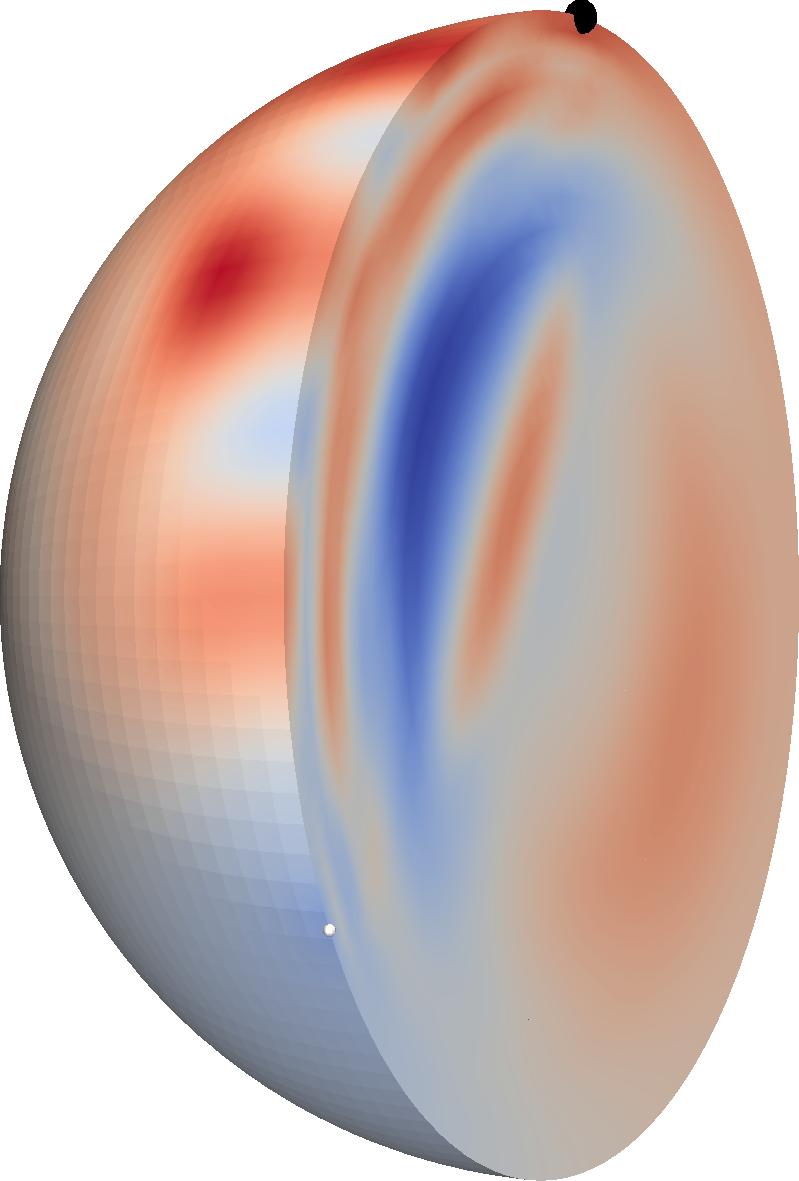}
  \caption{Problem I: Eigenvectors of the prior-preconditioned misfit
    Hessian corresponding to the first (i.e., the largest), the 3rd,
    the 5th, 8th and 13th eigenvalues (from left to right). The
    visualization employs a slice through the source and receiver
    locations.\label{fig:eigenvectorsI}}
\end{figure}

\begin{figure}
\centering
 \includegraphics[width=.19\columnwidth]
                 {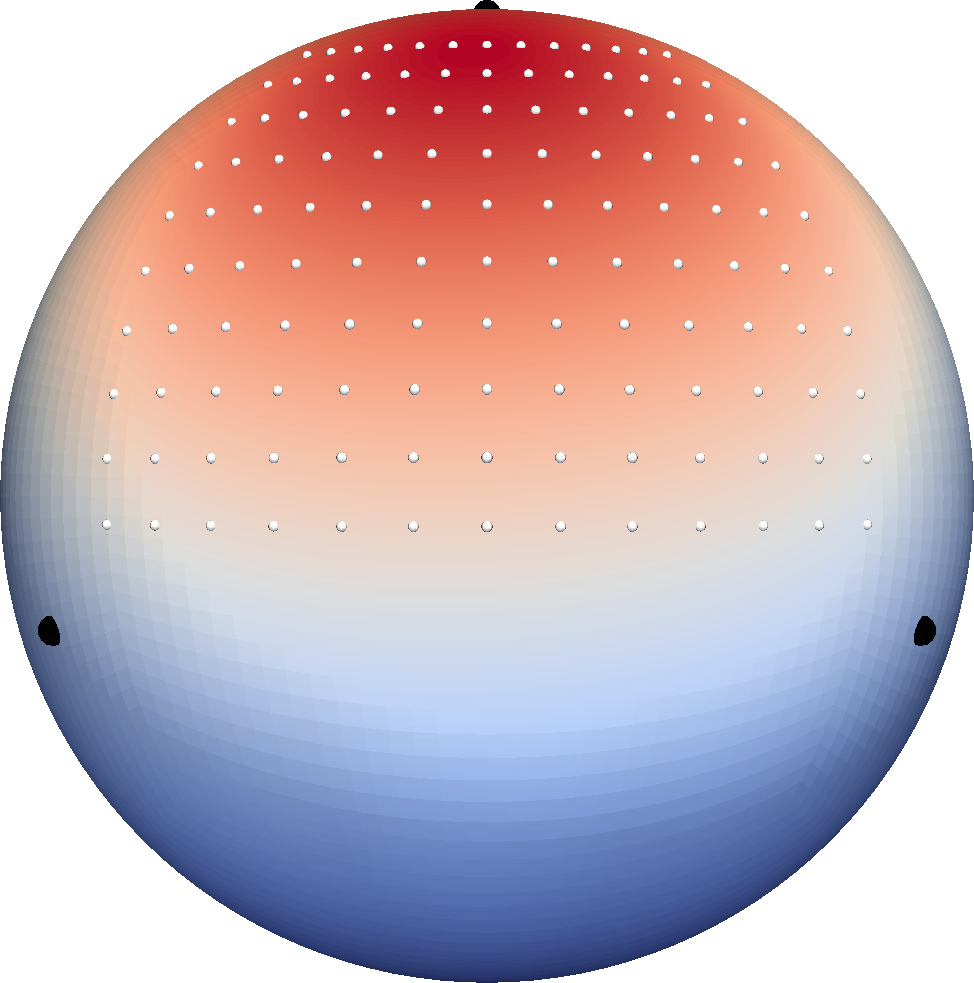} \hfill
 \includegraphics[width=.19\columnwidth]
                 {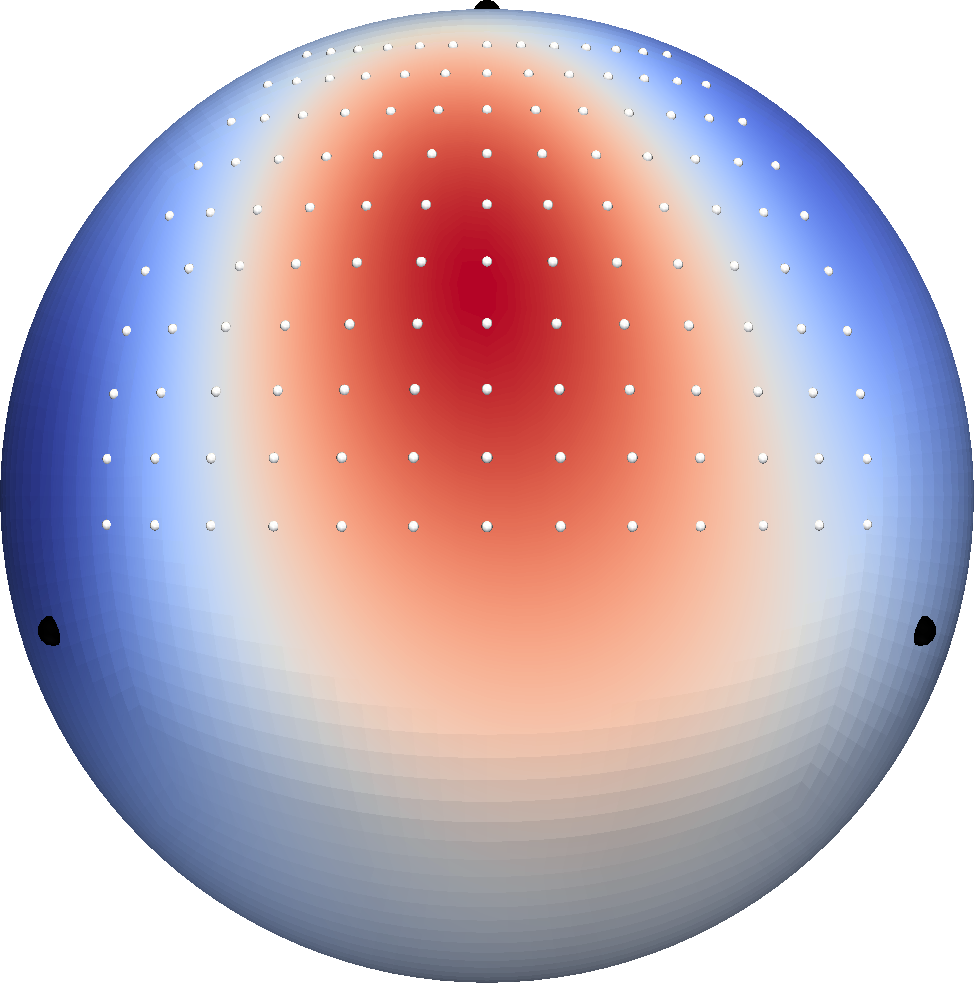} \hfill
 \includegraphics[width=.19\columnwidth]
                 {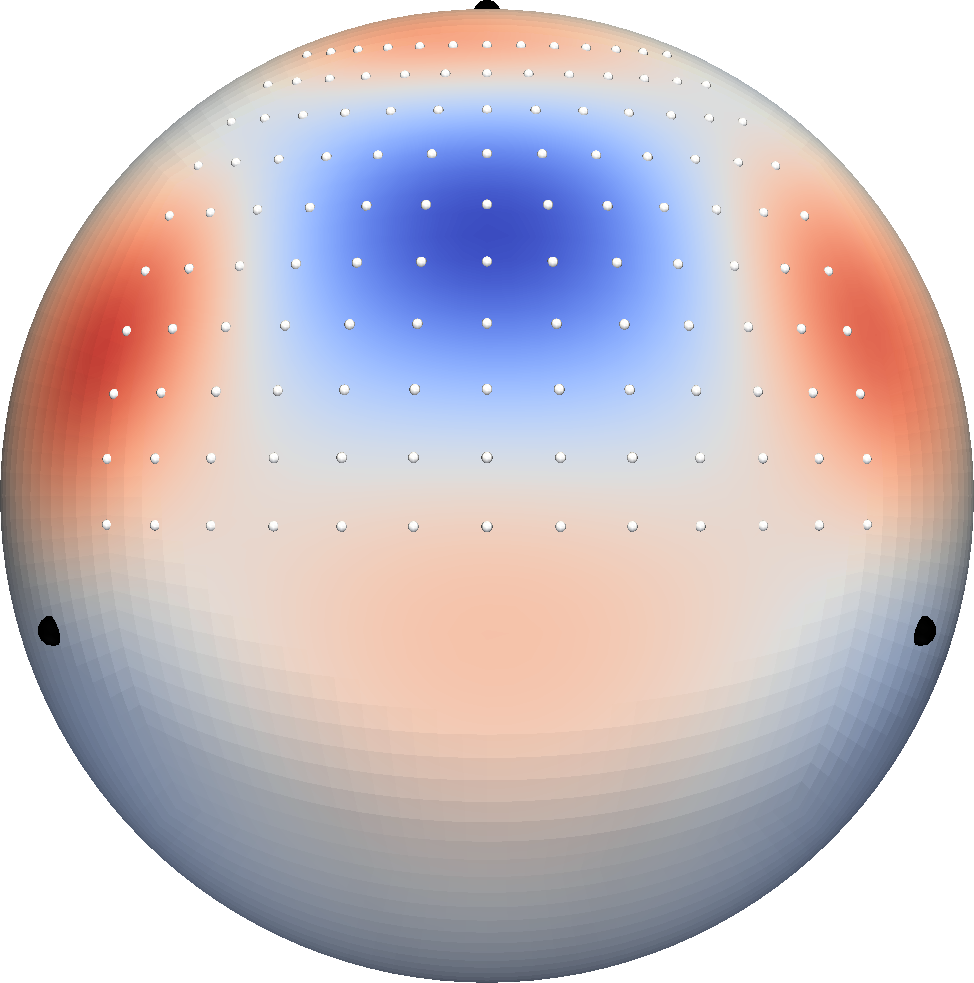} \hfill
 \includegraphics[width=.19\columnwidth]
                 {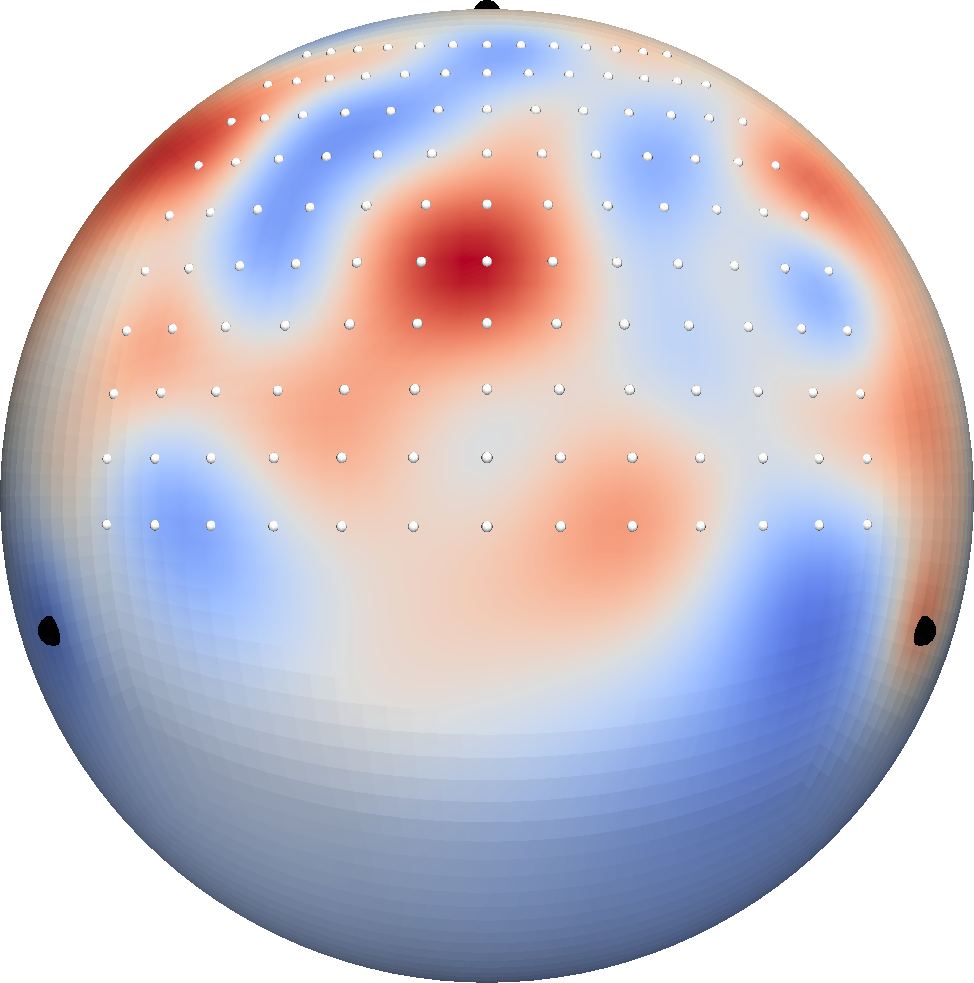} \hfill
 \includegraphics[width=.19\columnwidth]
                 {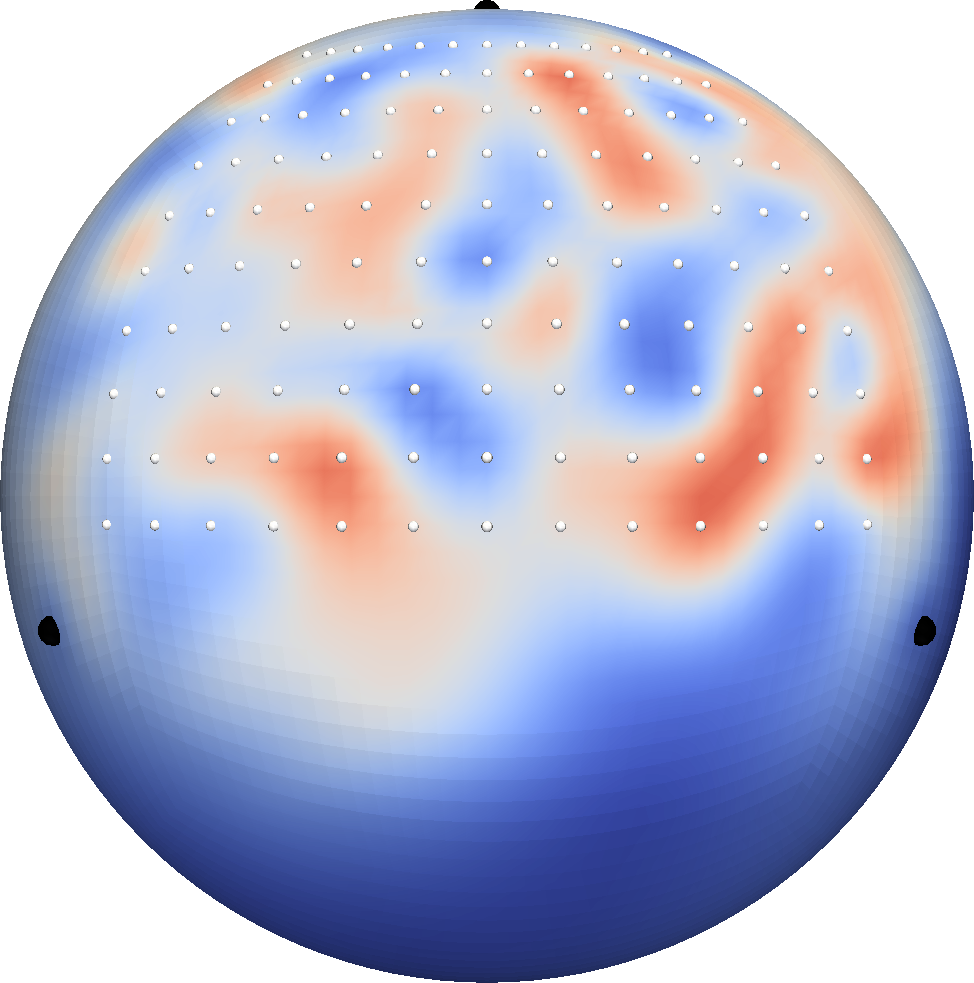}
\caption{Problem II: Eigenvectors of the misfit Hessian
  corresponding to eigenvalues 1, 5, 20, 100 and 350 respectively.
  Note that the lower  modes are smoothest and become
  more oscillatory with increasing mode number.}
\label{fig:eigenvectorsII}
\end{figure}

Figures \ref{fig:eigenvectorsI} and \ref{fig:eigenvectorsII} show
several eigenvectors of the prior preconditioned data misfit Hessian
\eqref{Hmisfit} (corresponding to several dominant eigenvalues) for
Problems I and II.
Eigenvectors corresponding to dominant eigenvalues represent the earth
modes that are ``most observable'' from the data, given the
configuration of sources and receivers. As can be seen in these
figures, the largest eigenvalues produce the smoothest modes, and as
the eigenvalues decrease, the associated eigenvectors become more
oscillatory, due to the reduced ability to infer smaller length scales
from the observations. 


\subsection{Interpretation of the uncertainty in the solution of the
  inverse problem} 
\label{sec:wave_resultsUQ}

\begin{figure}\centering
    \framebox{
      \begin{minipage}{.85\columnwidth}
        \begin{minipage}{0.21\columnwidth}
          \includegraphics[width=1.0\columnwidth]
                          {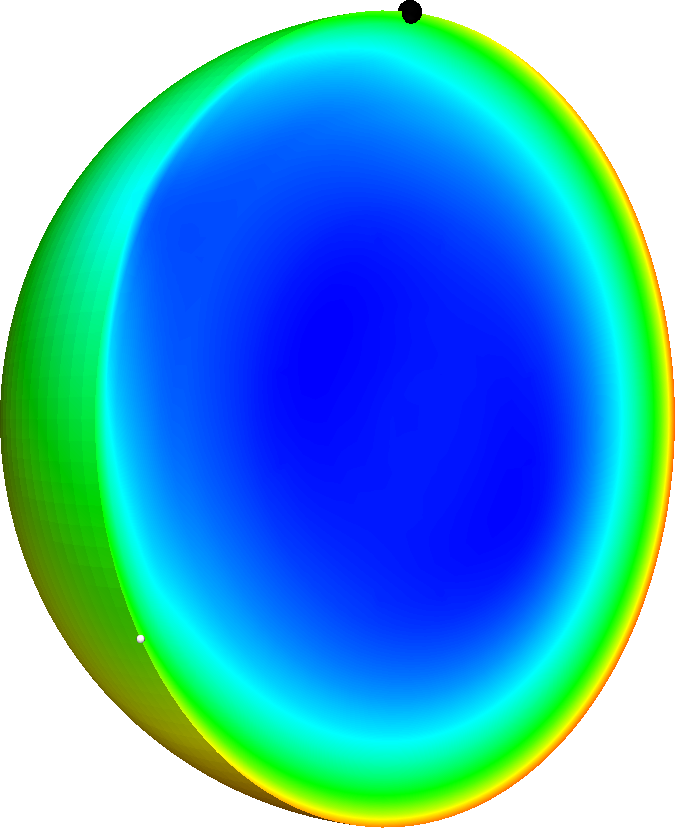}
        \end{minipage}
        \hfill
        \begin{minipage}{0.21\columnwidth}
          \includegraphics[width=1.0\columnwidth]
                          {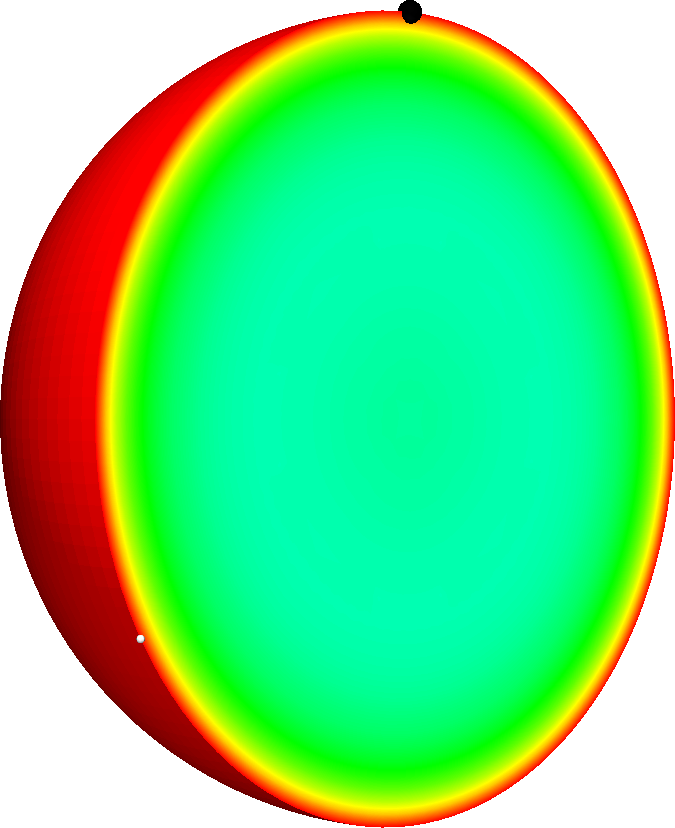}
        \end{minipage}
        \hfill
        \begin{minipage}{0.21\columnwidth}
          \includegraphics[width=1.0\columnwidth]
                          {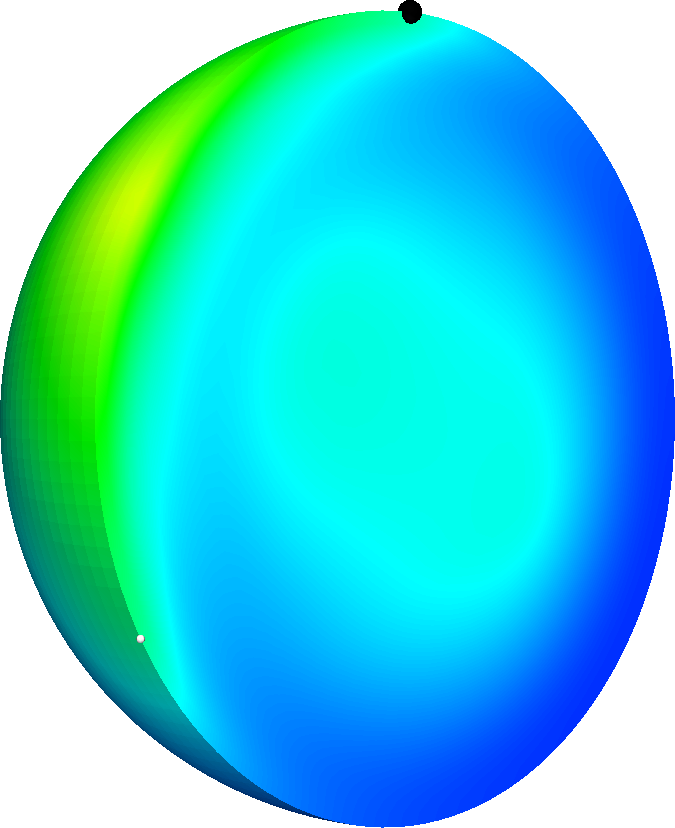}
        \end{minipage}
        \hspace{1ex}
        \begin{minipage}{0.1\columnwidth}
          \includegraphics[width=1.0\columnwidth]
                          {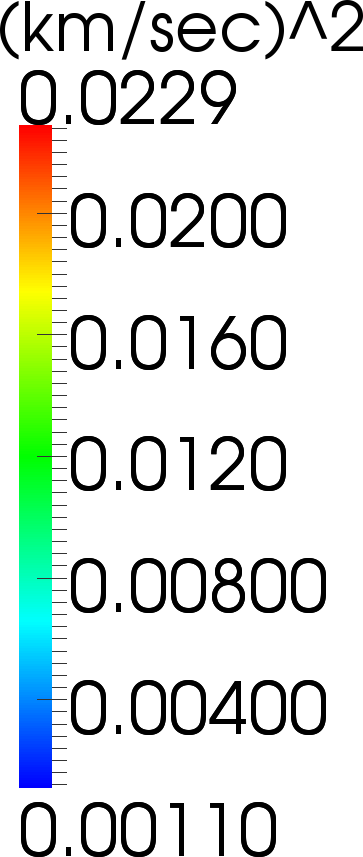}
        \end{minipage}
        \\
        \begin{minipage}{1.0\columnwidth}
          \[
          \phantom{df}\hspace{.5cm}
          \matrix{\Gamma}_{\text{post}}
          \hspace{1.1cm}
          \approx
          \hspace{1cm}
          \matrix{\Gamma}_{\text{prior}}
          \hspace{1cm}
          -
          \hspace{1cm}
          \matrix{\Gamma}_{\text{prior}}^{1/2} \matrix{V}_r \matrix{D}_r \matrix{V}_r^T \M
          \matrix{\Gamma}_{\text{prior}}^{1/2}
          \hspace{0.5cm}
          \]
        \end{minipage}
    \end{minipage}
    }
  \caption{Problem I: The left image depicts the pointwise posterior
    variance field, which is represented as the difference between the
    original prior variance field (middle), and the reduction in
    variance due to data (right).  The locations of the single source
    and single receiver is shown by the black and white dot,
    respectively.  Colorscale is common to all three images.}\label{fig:var_setupI}
\end{figure}

We first study Problem I, i.e., the single source, single receiver
problem. Since the data are very sparse, it is expected that we
can reconstruct only very limited information from the truth
solution; this is reflected in the smoothness of the dominant
eigenmodes shown in Figure \ref{fig:eigenvectorsI}. 
To assess the uncertainty, Figure \ref{fig:var_setupI} shows prior
variance, knowledge gained from the data (i.e., reduction in the
variance), and posterior variance, which are computed from 
\eqref{eq:postvarfield}. 
As discussed in Section~\ref{lowrank}, the posterior is the
combination of the prior information and the knowledge gained from the
data, so that the posterior uncertainty is decreased relative to the
prior uncertainty.  That is, the inference has less uncertainty in
regions for which the data are more informative.  In particular, the
region of lowest uncertainty is at the surface half-way between source
and receiver, as Figure \ref{fig:var_setupI} shows.
Note that the data are also informative about the core-mantle
boundary, where the strong material contrast results in stronger
reflected energy back to the surface receivers, allowing greater
confidence in the properties of that interface.

Next, we study the results for Problem II. The comparison between the
MAP estimate and the ground truth earth model (S20RTS) at different
depths is displayed in \figref{targetmap}. As can be seen, we are able
to recover accurately the wave speed in the portion of the Northern
hemisphere covered by sources and receivers. 
We plot the prior and posterior pointwise standard deviations in
\figref{variance}.
One observes that the uncertainty reduction is
greatest along the wave paths between sources and receivers,
particularly in the quarter of the Northern hemisphere surface where
the receivers are distributed. 

In \figref{NA_samples}, we show a comparison between samples from the
prior distribution and from the posterior.  We observe that in the
quarter of the Northern hemisphere where the data are more informative
about the medium, we have a higher degree of confidence about the wave
speed, which is manifested in the common large scale features across
the posterior samples.  The fine-scale features (about which the data
are least informative) are qualitatively similar to those of the prior
distribution, and vary from sample to sample in the posterior.  We
note that the samples shown here are computed by approximating
$\matrix M^{-1/2}$ in expression \eqnref{Lmatrix} using the (diagonal)
lumped mass matrix to avoid computing a factorization of $\matrix
M$. If desired, this mass lumping can be avoided by 
applying $\matrix M^{-1/2}$ to a vector using an iterative scheme
based on polynomial approximations to the matrix function $f(t) =
t^{-1/2}$, as in \cite{ChenAnitescuSaad11}.

\begin{figure}
  \begin{center}
    \begin{minipage}{0.2\columnwidth}
      \includegraphics[width=\columnwidth]{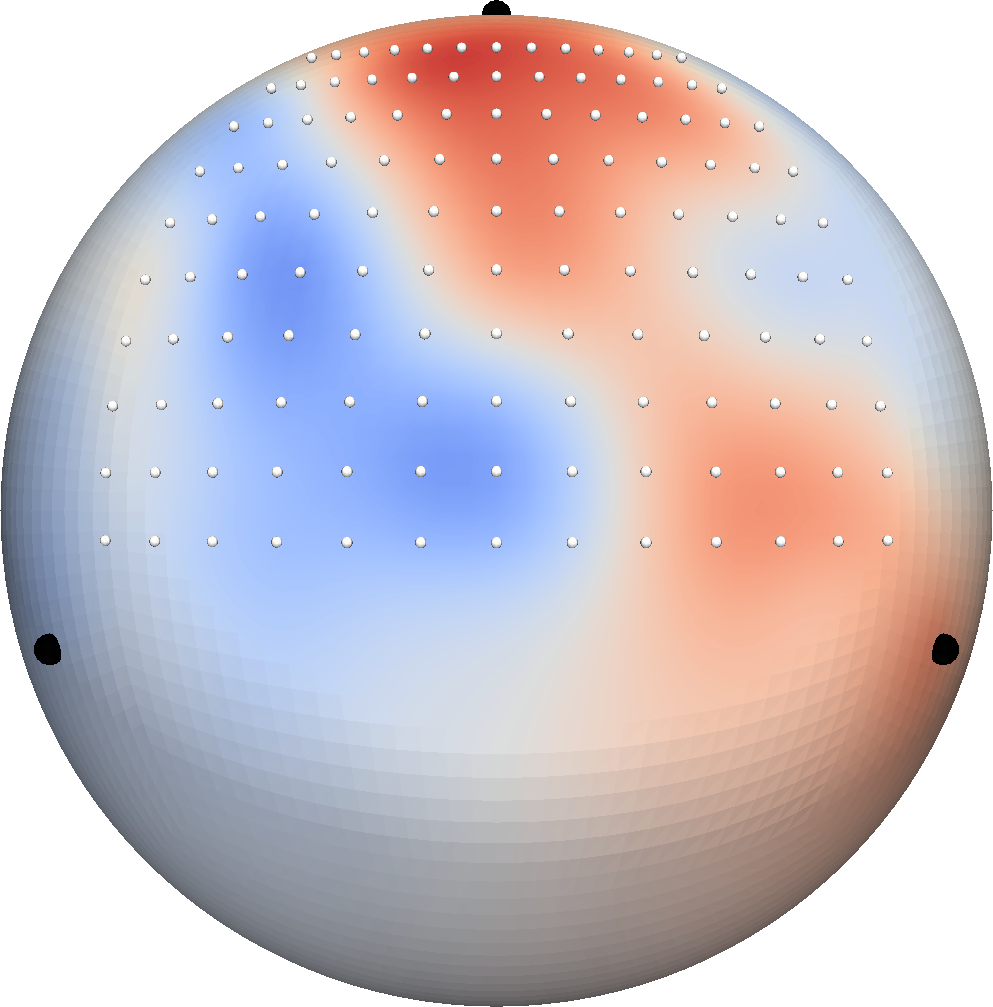} \\
      \includegraphics[width=\columnwidth]{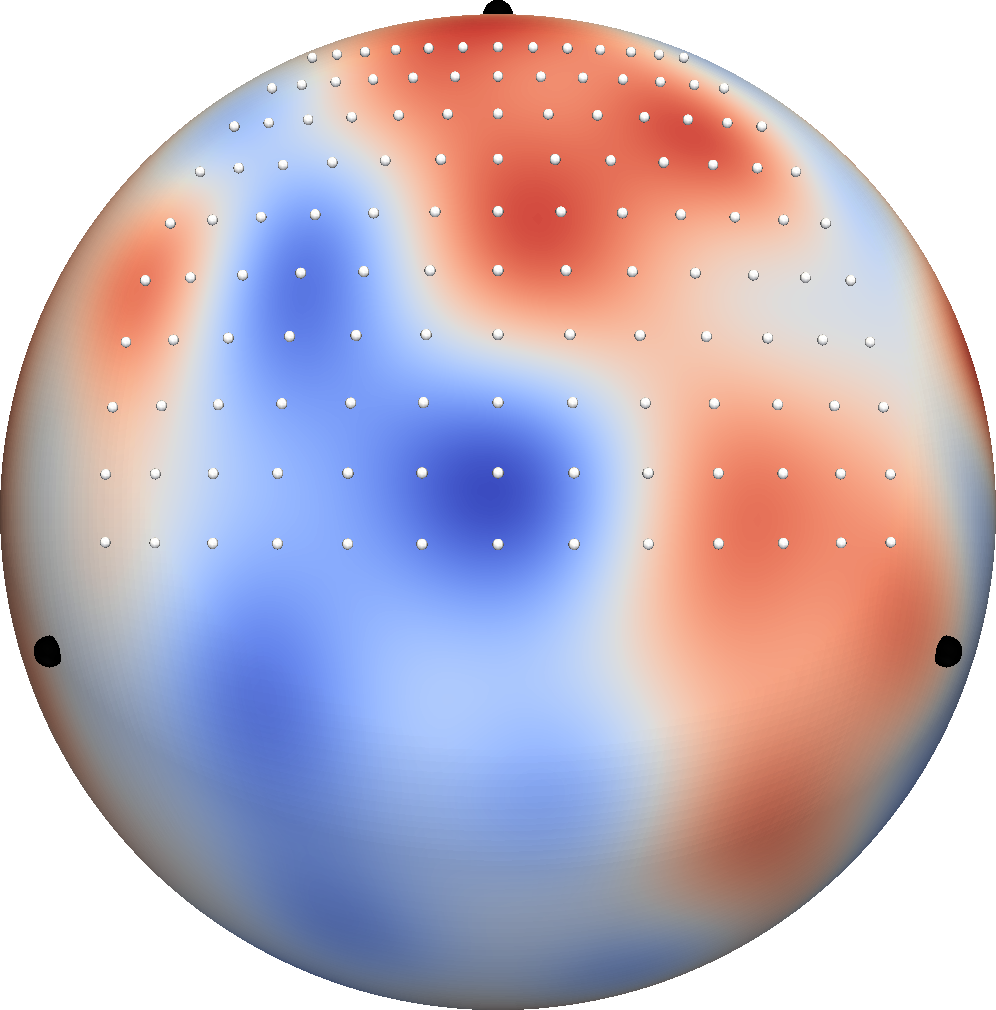}
    \end{minipage} \hfill
    \begin{minipage}{0.1\columnwidth}
    \includegraphics[width=\columnwidth]{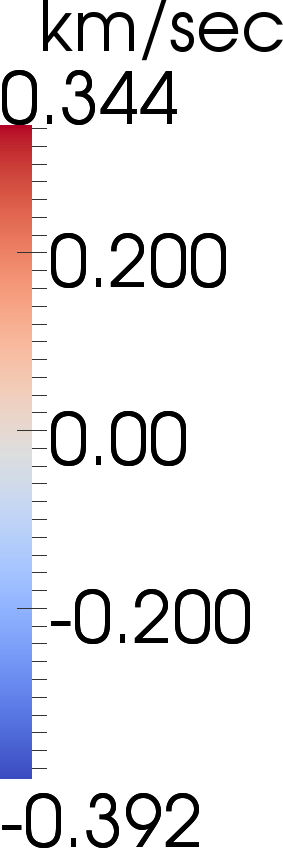} \hfill
    \end{minipage}
    \begin{minipage}{0.2\columnwidth}
      \includegraphics[width=\columnwidth]{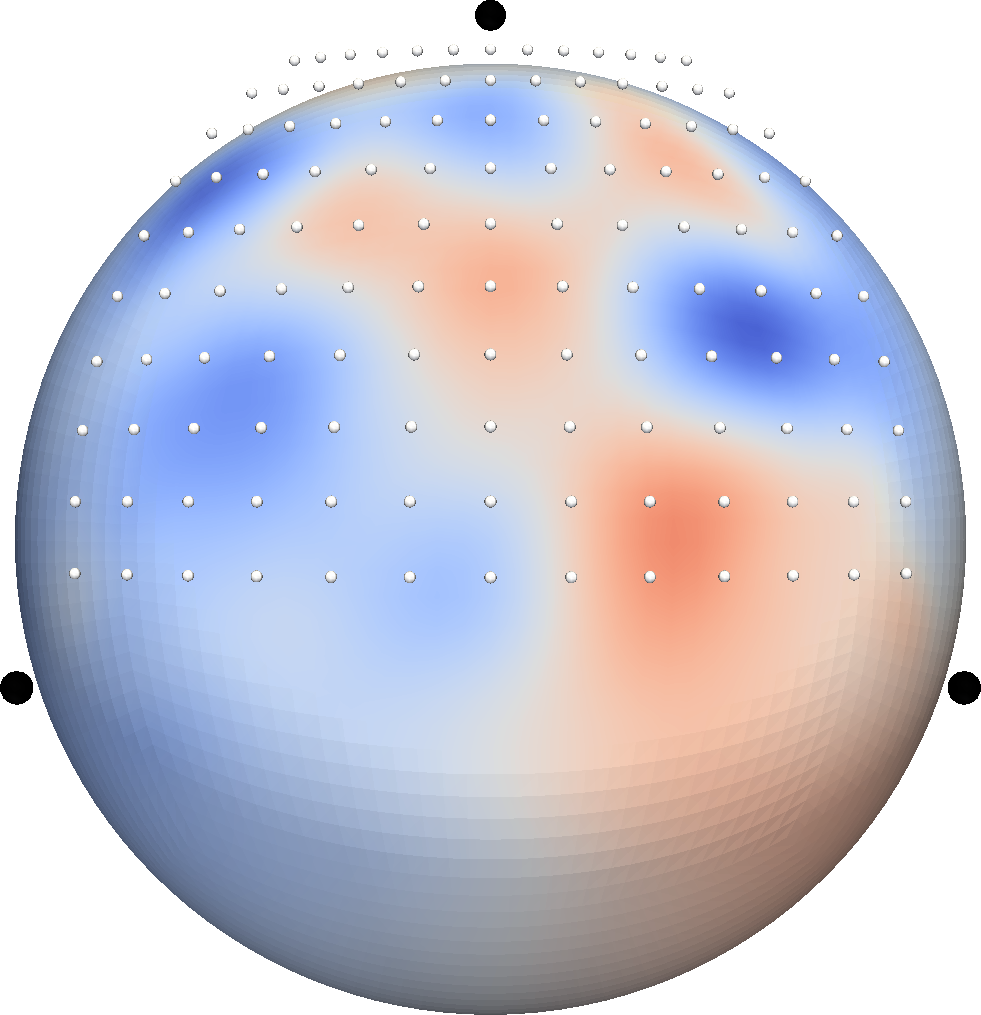} \\
      \includegraphics[width=\columnwidth]{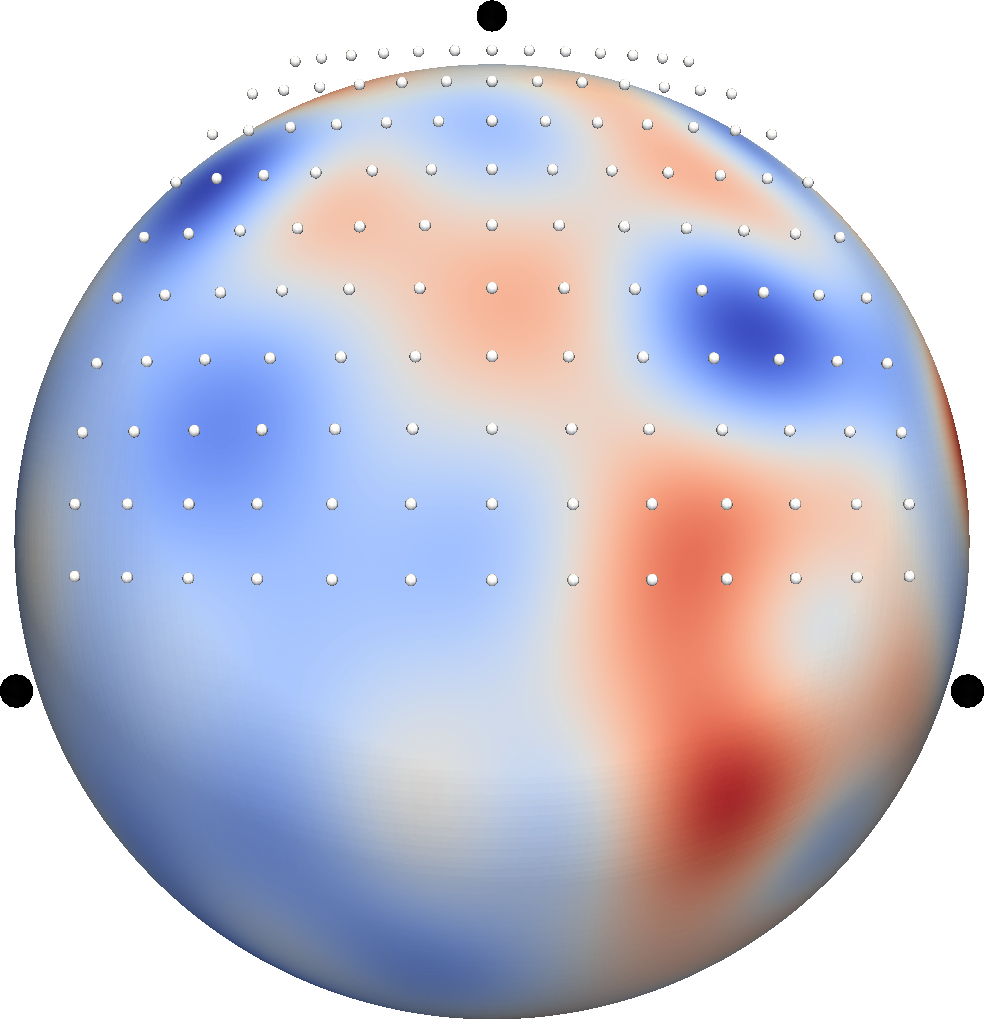}
    \end{minipage} \hfill
    \begin{minipage}{0.1\columnwidth}
    \includegraphics[width=\columnwidth]{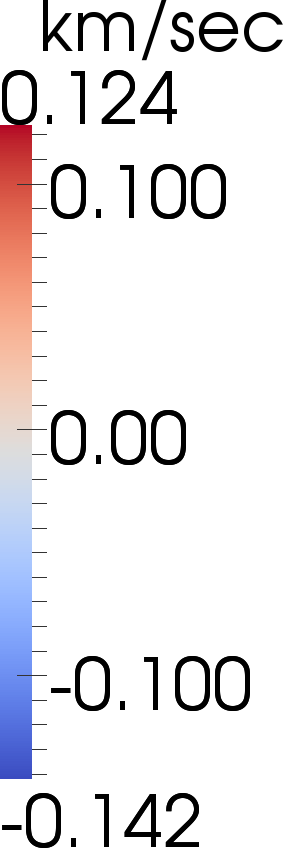} \hfill
    \end{minipage}
    \begin{minipage}{0.2\columnwidth}
      \includegraphics[width=\columnwidth]{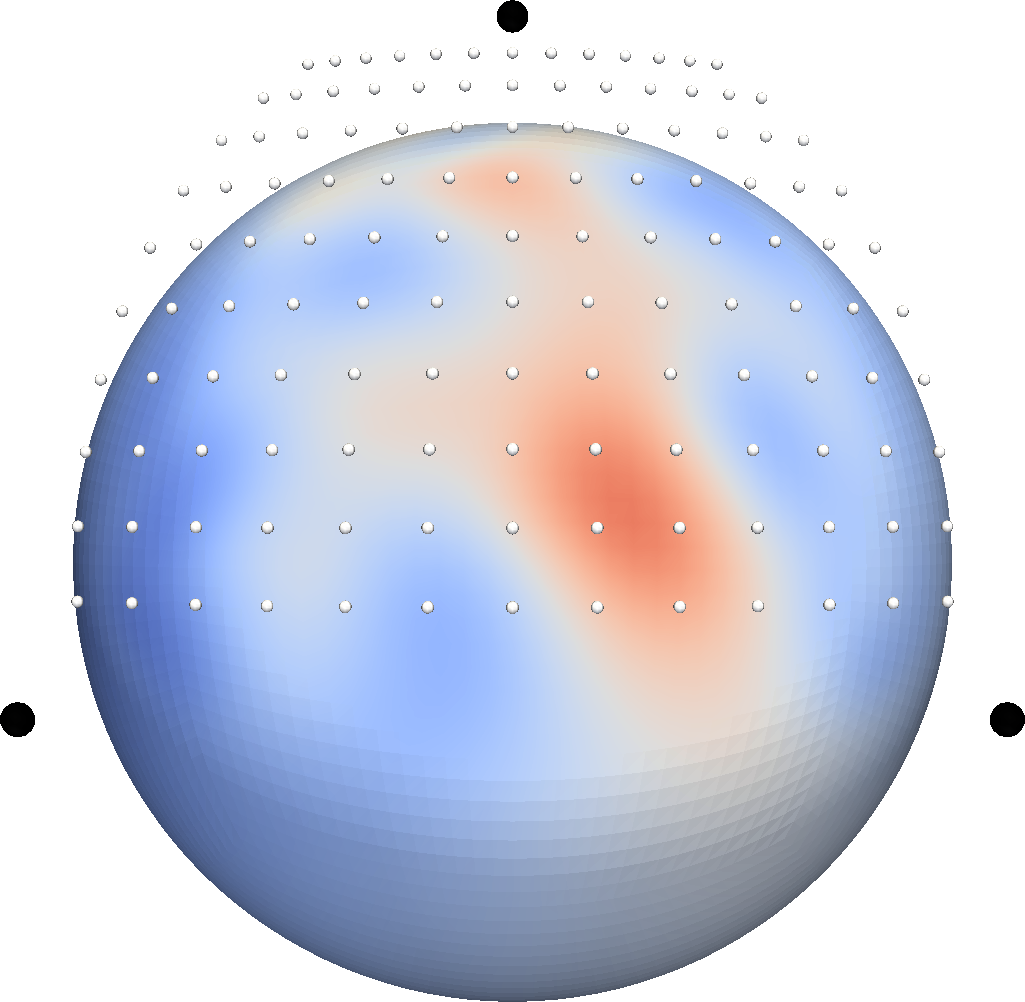} \\
      \includegraphics[width=\columnwidth]{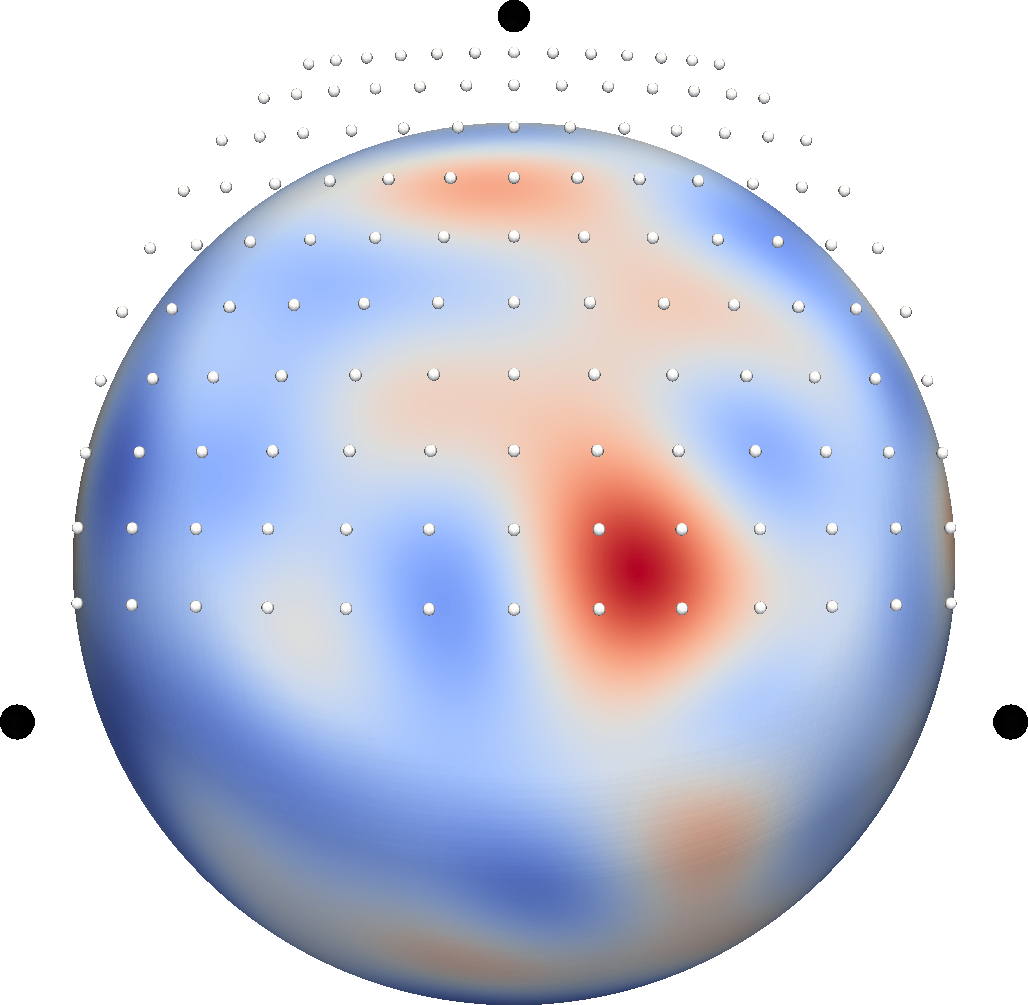}
    \end{minipage} \hfill
    \begin{minipage}{0.1\columnwidth}
    \includegraphics[width=\columnwidth]{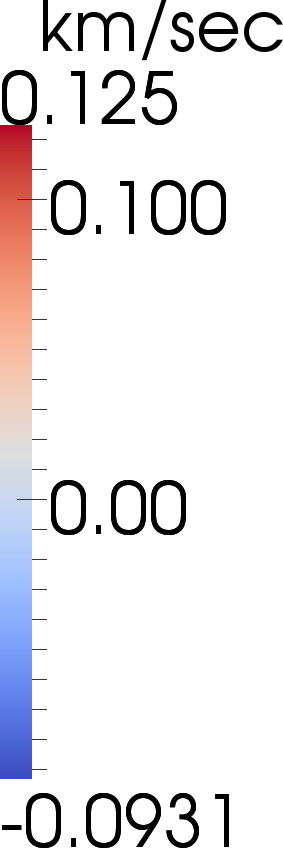} \hfill
    \end{minipage}
  \end{center}
  \caption{Comparison of MAP of posterior pdf (upper row) with the
    ``truth'' earth model (lower row) in a depth of $67km$
    (left image), $670km$ (middle image) and $1340km$ (right image).
    The colormap varies with depth, but is held
    constant between the MAP and ``truth'' images at each depth.}
  \figlab{targetmap}%
\end{figure}%

\begin{figure}\centering
    \includegraphics[height=0.23\columnwidth]
{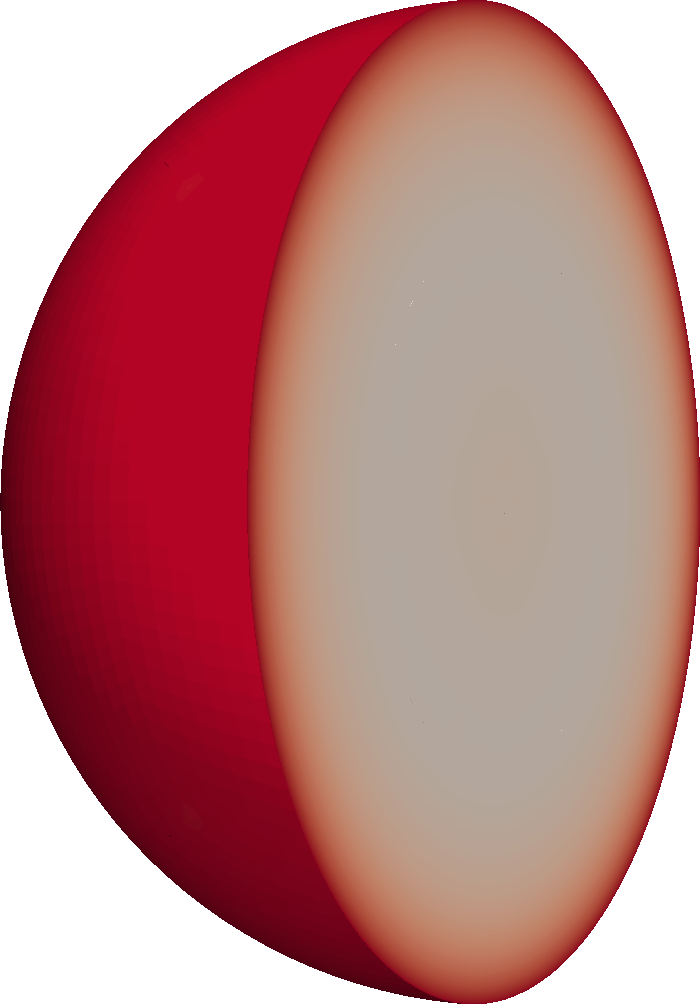}\hfill
    \includegraphics[height=0.23\columnwidth]
{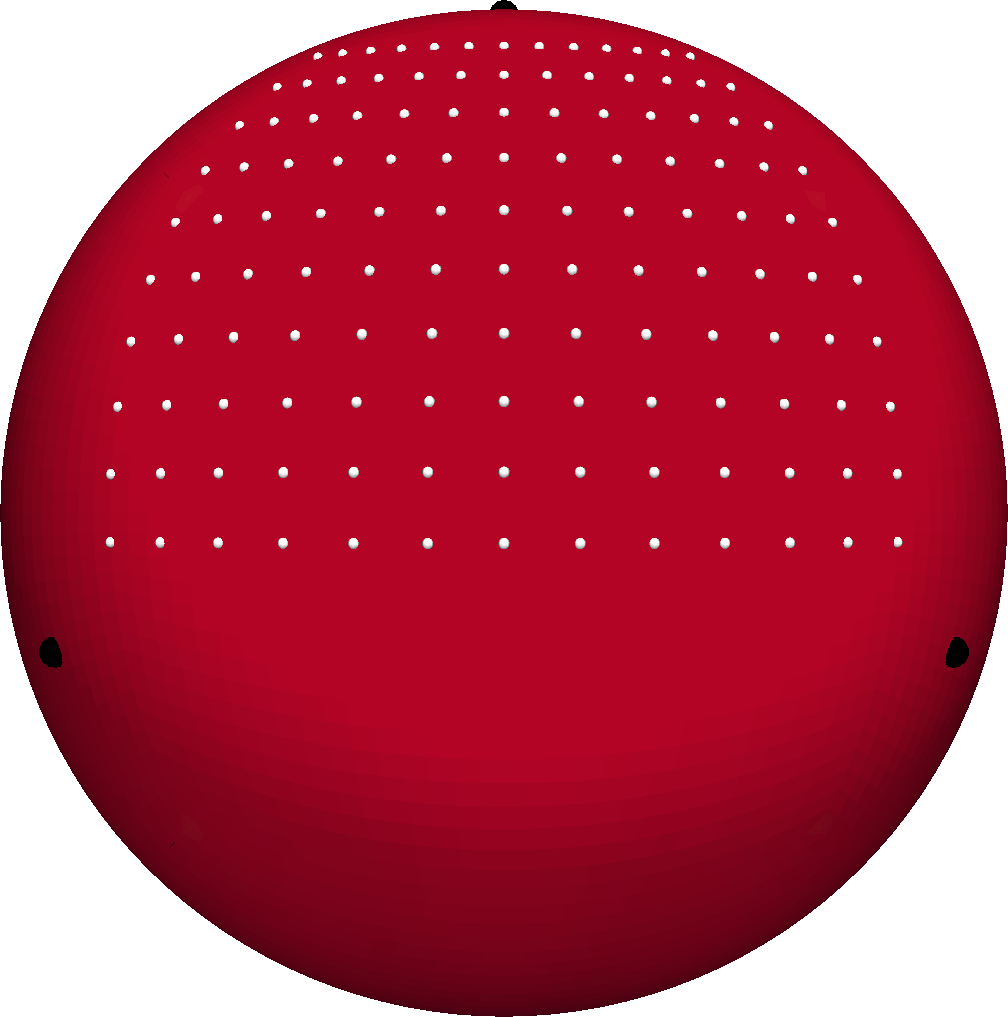}\hfill\hfill
    \includegraphics[height=0.23\columnwidth]
{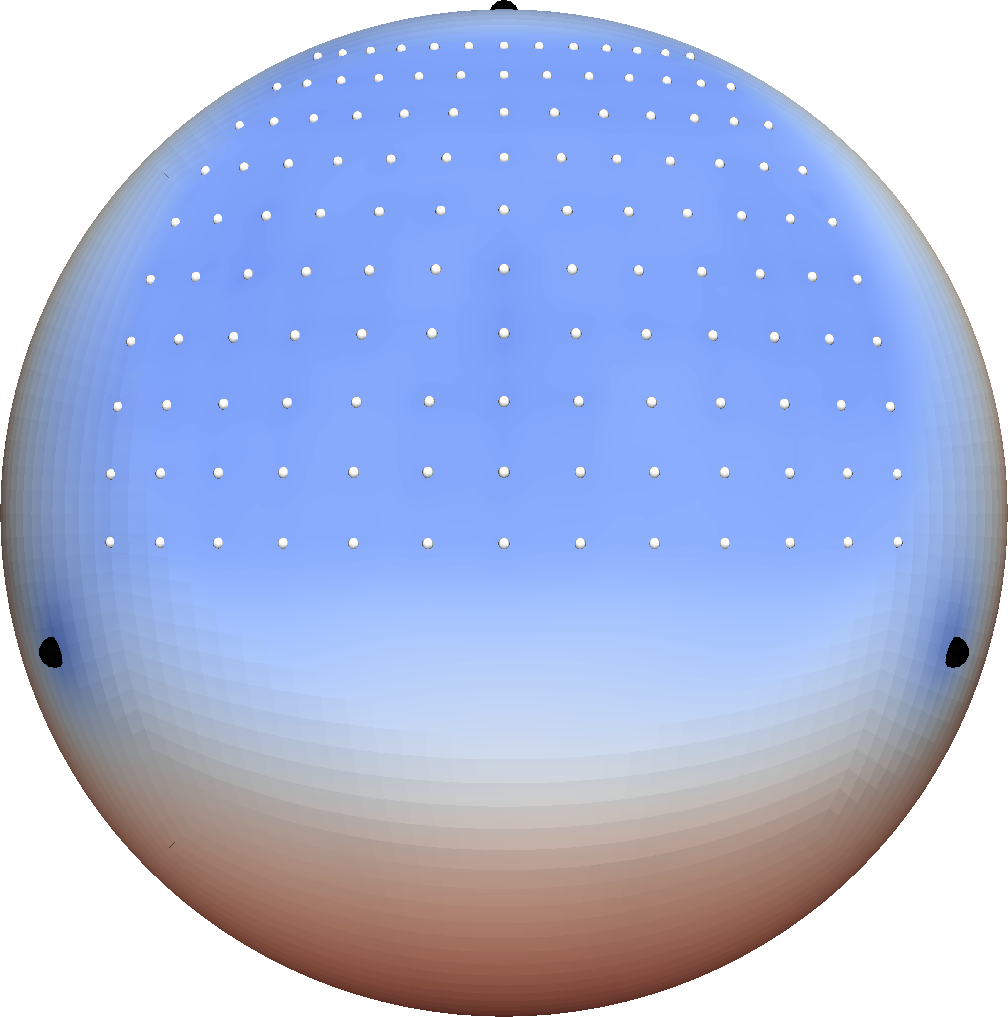}\hfill
    \includegraphics[height=0.23\columnwidth]
{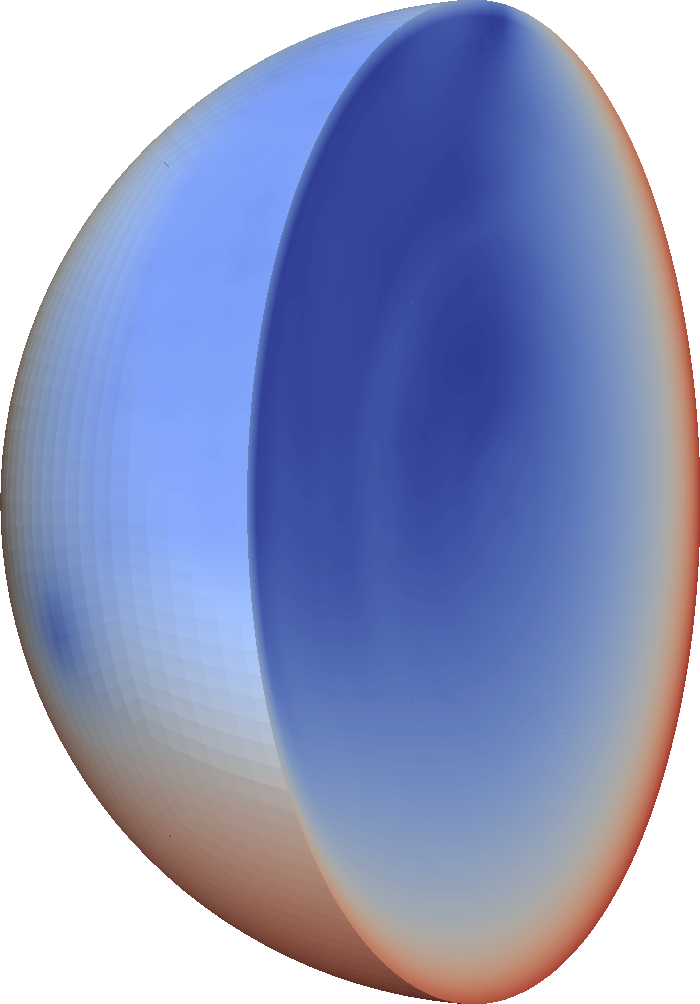}\hfill
    \includegraphics[height=.23\columnwidth]
{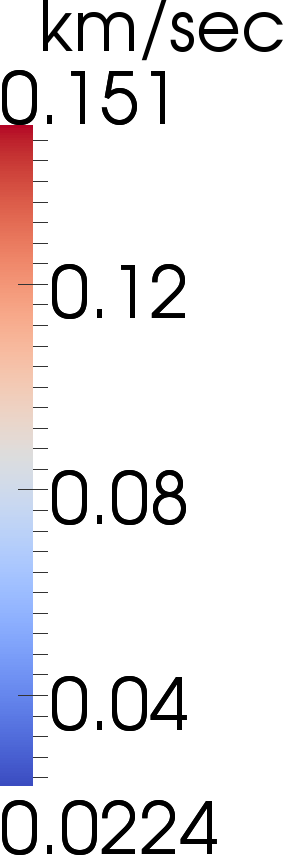}
  \caption{The figure compares the pointwise standard deviation 
    for the prior (left) and posterior (right) distributions at a
    depth of $67km$.  The color indicates one standard deviation, and
    the scale is common to both prior and posterior images.  We
    observe that the most reduction in variance due to data occurs in
    the region near sources and receivers, whereas the least reduction
    occurs on the opposite side of the Earth.}  \figlab{variance}
\end{figure}%

\begin{figure}
\centering
\begin{minipage}{0.7\columnwidth}
  \includegraphics[width=.24\columnwidth]{prior_sample_3} \hfill
  \includegraphics[width=.24\columnwidth]{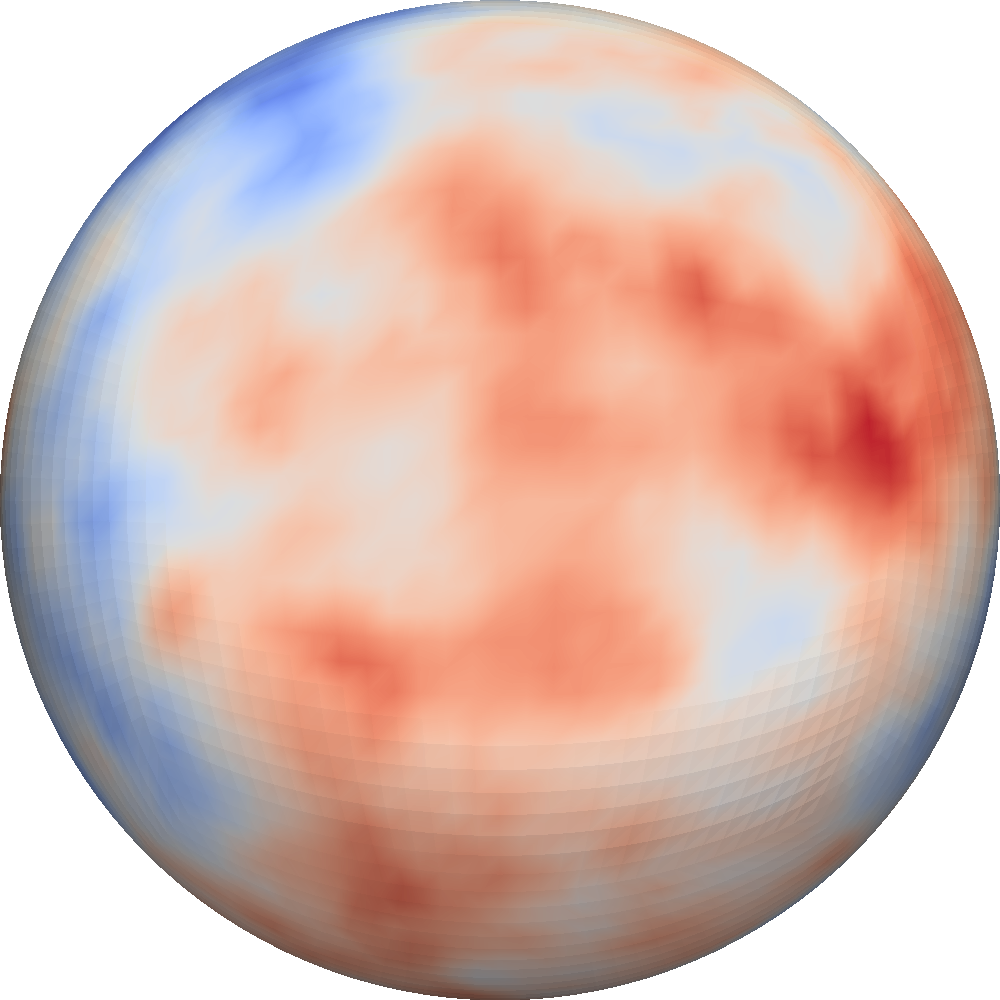} \hfill
  \includegraphics[width=.24\columnwidth]{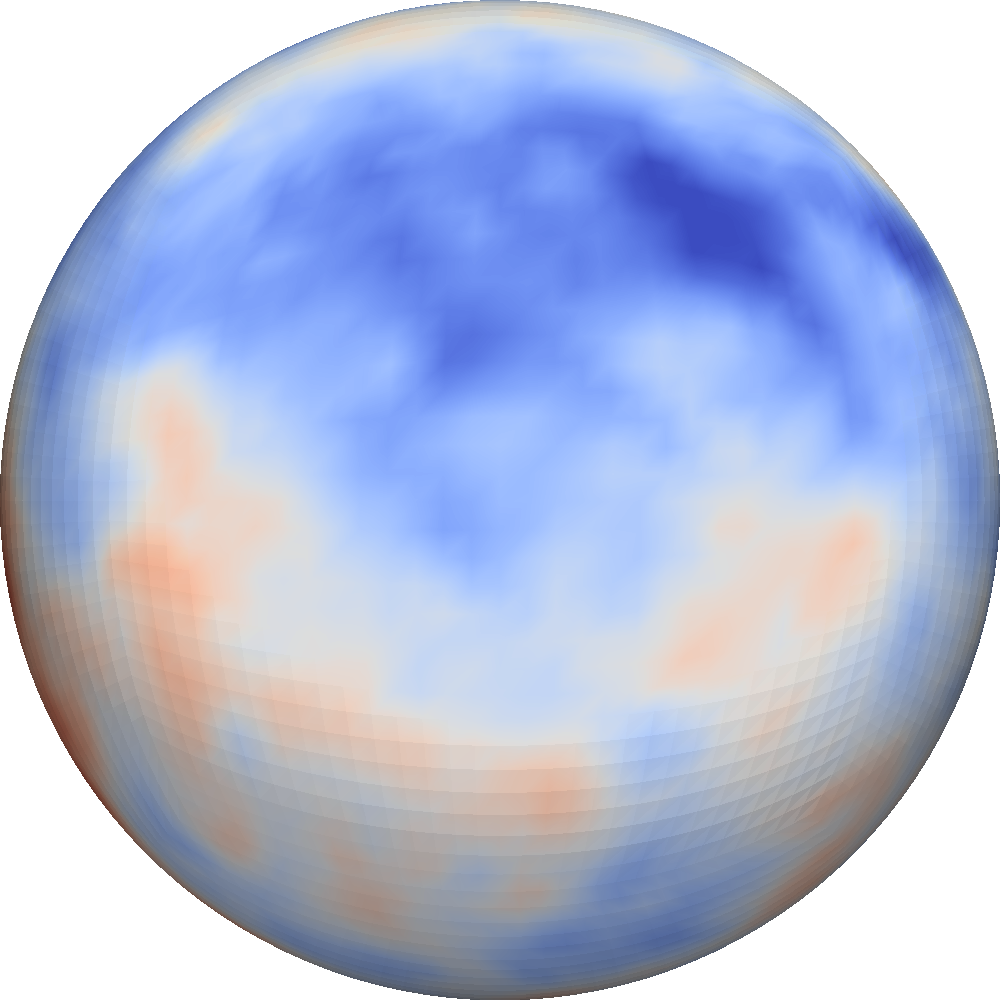} \hfill
  \includegraphics[width=.24\columnwidth]{prior_sample_6} \\
  \includegraphics[width=.24\columnwidth]{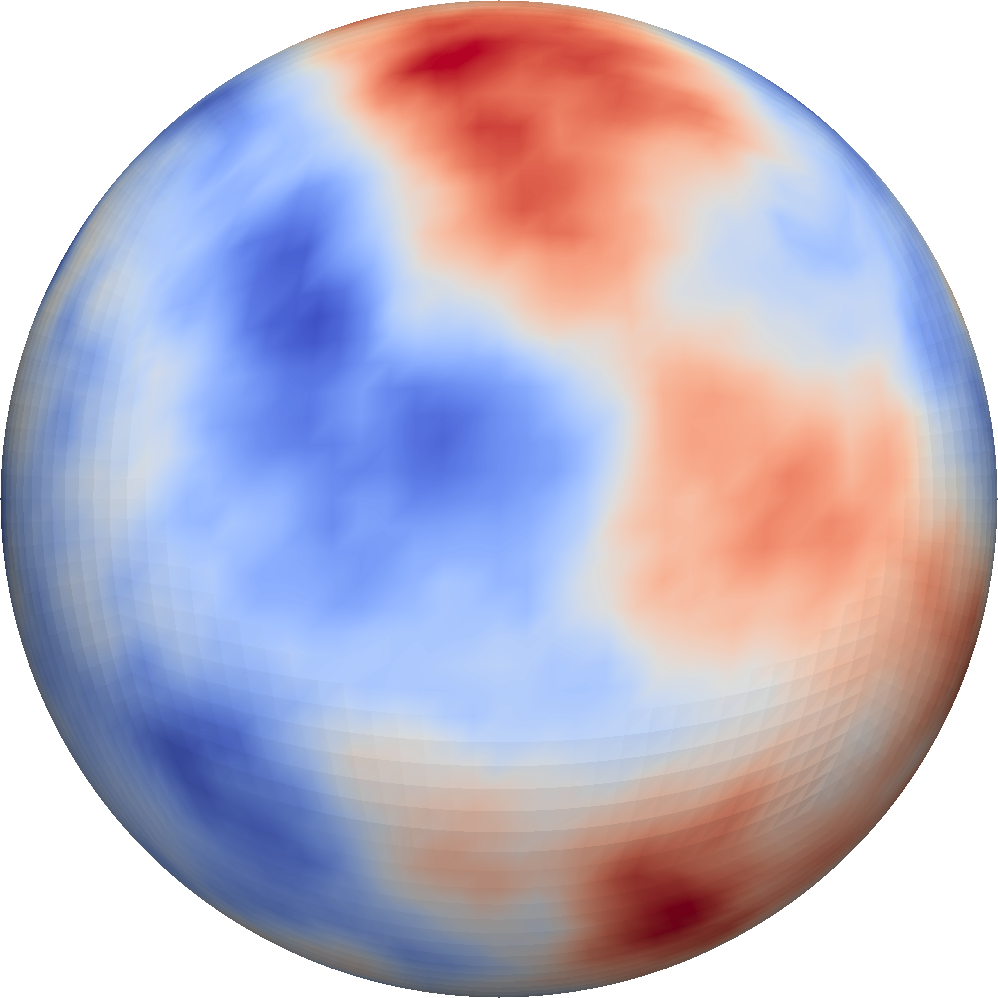} \hfill
  \includegraphics[width=.24\columnwidth]{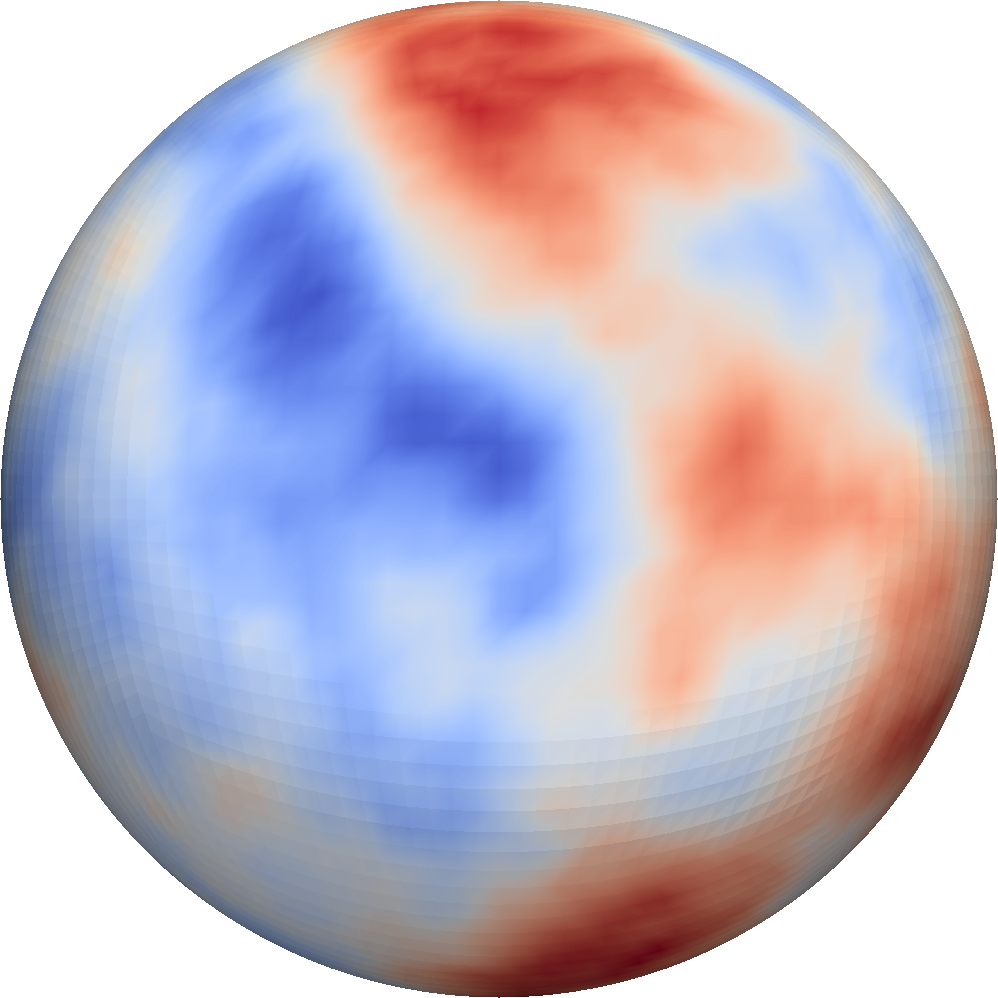} \hfill
  \includegraphics[width=.24\columnwidth]{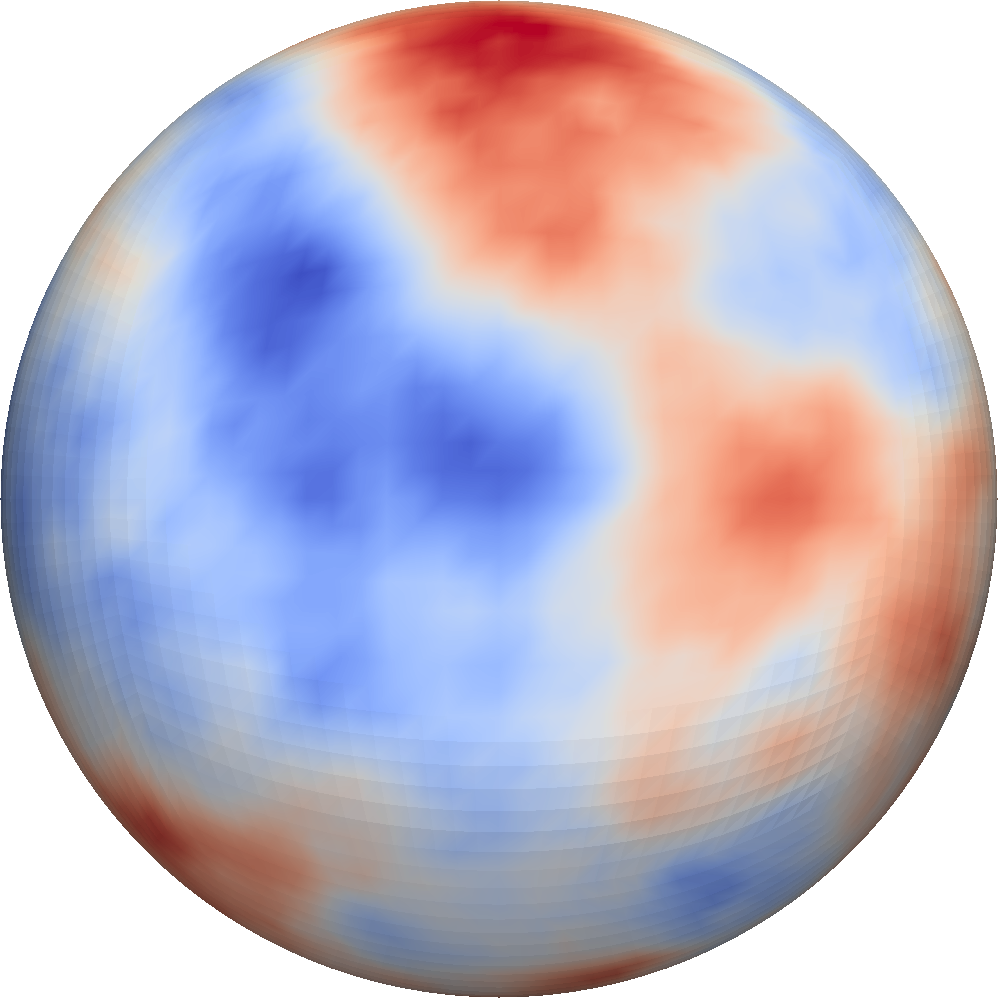} \hfill
  \includegraphics[width=.24\columnwidth]{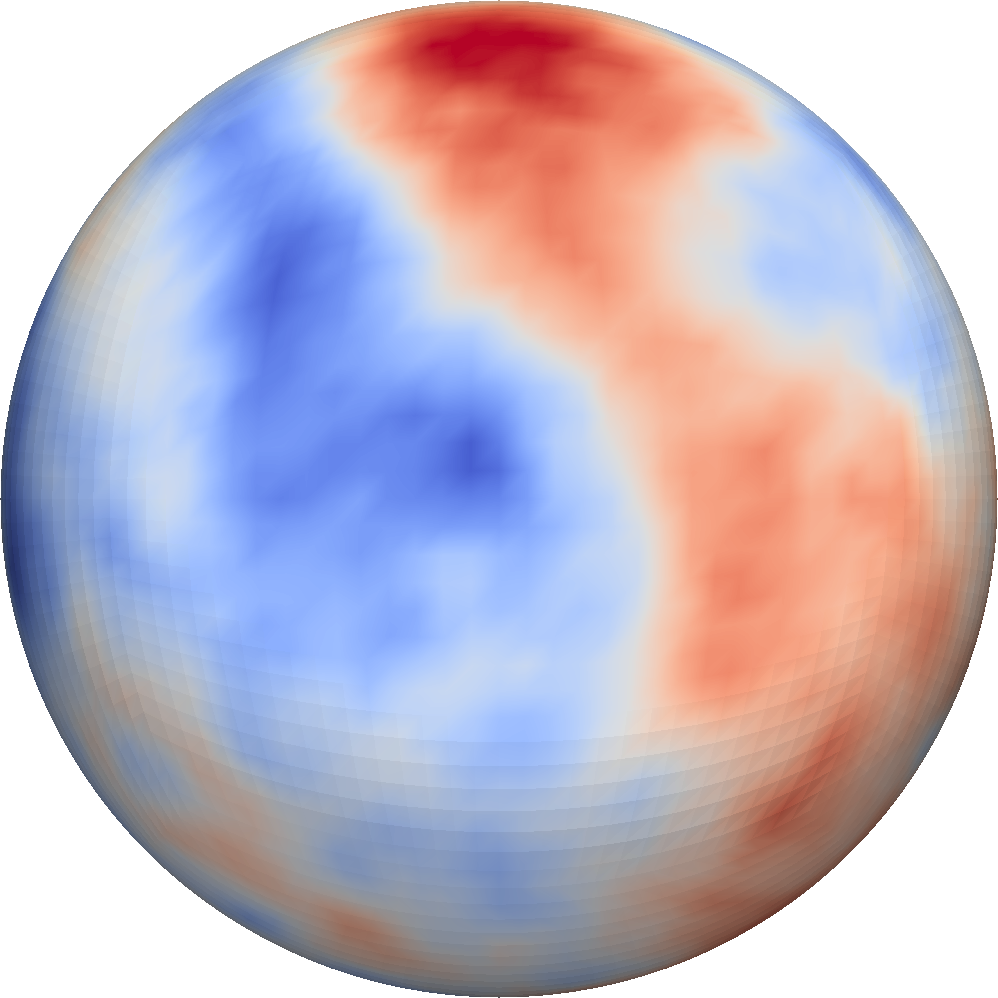}
\end{minipage}
\begin{minipage}{0.18\columnwidth}
  \includegraphics[width=\columnwidth]{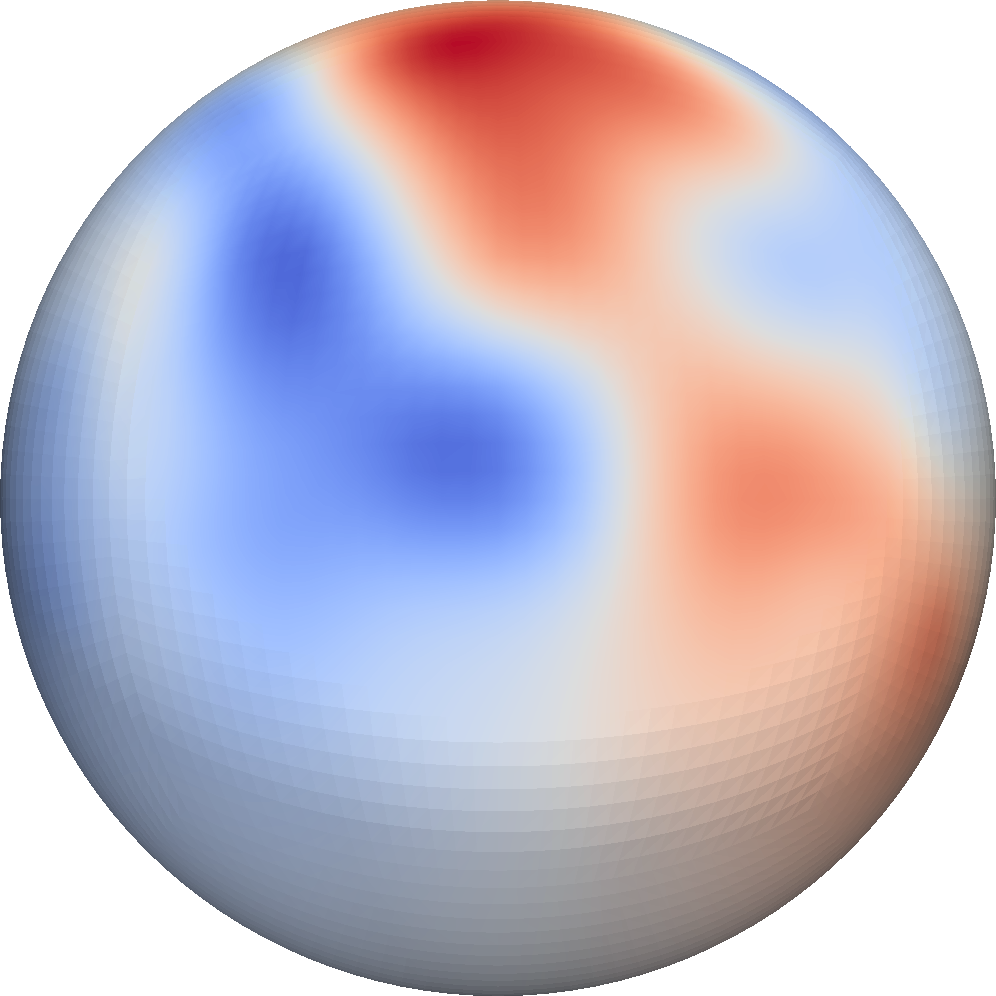}
\end{minipage}
\begin{minipage}{0.09\columnwidth}
  \includegraphics[width=\columnwidth]{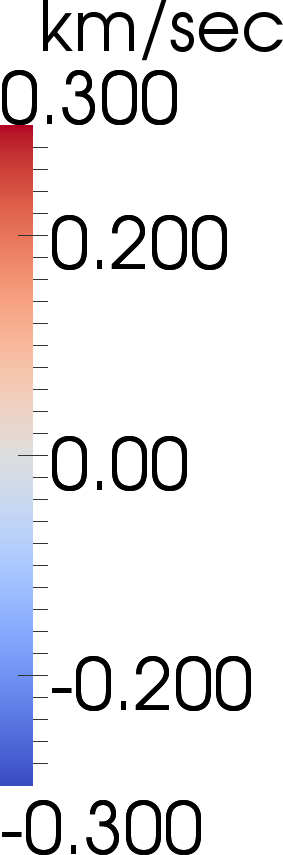}
\end{minipage}

\caption{Samples from the prior (top row) and posterior (bottom row)
  distributions.  The prior scaling was chosen such that the ``ground
  truth'' S20RTS would be a qualitatively reasonable sample from the
  prior distribution.  For comparison purposes, the MAP estimate is
  shown on the far right.}
\figlab{NA_samples}
\end{figure}

\section{Conclusions}

A computational framework for estimating the uncertainty in the
numerical solution of linearized infinite-dimensional statistical
inverse problems is presented. We adopt the Bayesian inference
formulation: given observational data and their uncertainty, the
governing forward problem and its uncertainty, and a prior probability
distribution describing uncertainty in the parameter field, find the
posterior probability distribution over the parameter field.
The framework, which builds on the infinite-dimensional formulation
proposed by Stuart \cite{Stuart10}, incorporates a number of
components aimed at ensuring a convergent discretization of the
underlying infinite-dimensional inverse problem. It additionally
incorporates algorithms for manipulating the prior, constructing a low
rank approximation of the data-informed component of the posterior
covariance operator, and exploring the posterior, that together ensure
scalability of the entire framework to very high parameter
dimensions. Since the data are typically informative about only a low
dimensional subspace of the parameter space, the Hessian is sparse
with respect to some basis. We have exploited this fact to construct a
low rank approximation of the Hessian and its inverse using a parallel
matrix-free Lanczos method. Overall, our method requires a
dimension-independent number of forward PDE solves to approximate the
local covariance. Uncertainty quantification for the linearized
inverse problem thus reduces to solving a fixed number of forward and
adjoint PDEs (which resemble the original forward problem),
independent of the problem dimension. The entire process is thus
scalable with respect to the forward problem dimension, uncertain
parameter dimension, and observational data dimension. We applied this
method to the Bayesian solution of an inverse problem in 3D global
seismic wave propagation with up to 430,000 parameters, for which we
observe 2--3 orders of magnitude dimension reduction, making UQ for
large-scale inverse problems tractable.

\section*{Appendix}
\seclab{appendix} In the following, we provide a constructive
derivation of $\L$ in \eqnref{Lmatrix} such that it satisfies
$\matrix{\Gamma}_{\text{post}} = \L \Ladj$. Our goal is to draw
posterior Gaussian random sample with covariance matrix
$\matrix{\Gamma}_{\text{post}}$ in $\R^n_\M$. To accomplish this, a
standard approach is first to find a factorization
$\matrix{\Gamma}_{\text{post}} = \tilde{\L} \adjMacroMM{\tilde{\L}}$,
where $\tilde{\L}$ is a linear map from $\R^n_\M$ to $\R^n_\M$. Then,
any random sample from the posterior can be written as
\begin{equation}
\eqnlab{postSample}
\bs{\nu}^{\text{post}} = \dparmap + \tilde{\L}\tilde{\bs{n}},
\end{equation}
where $\tilde{\bs{n}}$ is a Gaussian random
sample with zero mean and identity covariance matrix in
$\R^n_\M$. It follows that 
\[
\tilde{\bs{n}} = \M^{-1/2}\bs{n},
\] 
where $\bs{n}$ is the standard Gaussian random sample with zero mean
and identity covariance matrix in $\R^n$, i.e. $\bs{n} \sim \mathcal
N(0,\matrix{I})$, and $\M^{-1/2}$ a linear map from $\R^n$ to
$\R^n_\M$. 

Therefore, what remains to be done is to construct
$\tilde{\L}$. To begin the construction, we rewrite
\eqref{eqn:reduced_postfinal} as
\[
\matrix{\Gamma}_{\text{post}} \approx
\matrix{\Gamma}_{\text{prior}}^{1/2}\underbrace{\LRp{\matrix{I}
-{\V_r} \matrix{D}_r {\Vadj_r}}}_{\matrix{B}}\matrix{\Gamma}_{\text{prior}}^{1/2}.
\]
The simple structure of $\matrix{B}$ allows us to write its spectral
decomposition as
\[
\matrix{B} = \sum_{i=1}^n\vv_i\adjMacroEM{\vv_i} - \sum_{i=1}^r\frac{\lambda_i}{\lambda_i+1}\vv_i\adjMacroEM{\vv_i} = \sum_{i=1}^r\frac{1}{\lambda_i+1}\vv_i\adjMacroEM{\vv_i} + \sum_{i=r+1}^n\vv_i\adjMacroEM{\vv_i},
\]
which, together with the standard definition of the square root of
postive self-adjoint operators \cite{ArbogastBona08}, immediately gives
\[
\matrix{B}^{1/2} = \sum_{i=1}^r\frac{1}{\sqrt{\lambda_i+1}}\vv_i\adjMacroEM{\vv_i} + \sum_{i=r+1}^n\vv_i\adjMacroEM{\vv_i} = \sum_{i=1}^r\LRp{\frac{1}{\sqrt{\lambda_i+1}} - 1}\vv_i\adjMacroEM{\vv_i} + \sum_{i=1}^n\vv_i\adjMacroEM{\vv_i} = 
{\V_r} \matrix{P}_r {\Vadj_r} + \matrix{I},
\]
where $\matrix{P}_r = \diag\LRp{1/\sqrt{\lambda_1+1}-1, \hdots,
  1/\sqrt{\lambda_r+1}-1} \in \R^{r\times r}$, and $\matrix{B}^{1/2}$ is self-adjoint in $\R^n_\M$, namely, $\adjMacroMM{\LRp{\matrix{B}^{1/2}}} = \matrix{B}^{1/2}$. Now, we define
\[
\tilde{\L} = \matrix{\Gamma}_{\text{prior}}^{1/2} \matrix{B}^{1/2},
\]
which, by construction, is the desired matrix
 owing to the following trivial identity
\[\tilde{\L}
\adjMacroMM{\tilde{\L}} =  \matrix{\Gamma}_{\text{prior}}^{1/2} \matrix{B}^{1/2} \adjMacroMM{\LRp{\matrix{B}^{1/2}}} \adjMacroMM{\LRp{\matrix{\Gamma}_{\text{prior}}^{1/2}}} =
\matrix{\Gamma}_{\text{prior}}^{1/2} \matrix{B}\matrix{\Gamma}_{\text{prior}}^{1/2} = 
\matrix{\Gamma}_{\text{post}},
\]
where we have used the self-adjointness of
$\matrix{\Gamma}_{\text{prior}}^{1/2}$ and $ \matrix{B}^{1/2}$ in
$\R^n_\M$.

Finally, we can rewrite \eqnref{postSample} in terms of $\bs{n}$ and
$\L = \tilde{\L}\M^{-1/2}$, a linear map from $\R^n$ to $\R^n_\M$, as
\[
\bs{\nu}^{\text{post}} = \dparmap + \L\bs{n},
\]
where $\L$ satisfies the following desired identity
\[
\L\adjMacroEM{\L} = \tilde{\L} \M^{-1/2}
\adjMacroEM{\LRp{\M^{-1/2}}} \adjMacroMM{\tilde{\L}} = 
\tilde{\L} \M^{-1/2}
\M^{-1/2}\M \adjMacroMM{\tilde{\L}} =  \tilde{\L}\adjMacroMM{\tilde{\L}} =   \matrix{\Gamma}_{\text{post}}.
\]



\end{document}